\documentclass[twoside,11pt]{article}

\usepackage{blindtext}

\usepackage{jmlr2e}
\usepackage{comment}
\usepackage{amsmath}
\usepackage{enumerate}
\usepackage{multirow}
\usepackage{xcolor}
\usepackage{natbib}
\usepackage{bm}

\usepackage{lastpage}
\jmlrheading{24}{2024}{1-\pageref{LastPage}}{3/15; Revised ?/??}{?/??}{21-0000}{Yang Chen, Cheng-Der Fuh and Chu-Lan Michael Kao}

\ShortHeadings{Det. the Num. of States in HMM via Marg. L.H.}{Chen, Fuh and Kao}
\firstpageno{1}

\begin{document}

\title{Determine the Number of States in Hidden Markov Models via Marginal Likelihood}

\author{\name Yang Chen \email ychenang@umich.edu \\
       \addr Department of Statistics\\
       University of Michigan\\
      1085 South University
Ann Arbor, MI 48109-1107, USA
       \AND
       \name Cheng-Der Fuh \email cdffuh@gmail.com \\
       \addr Graduate Institute of Statistics\\
       National Central University\\
       No. 300, Zhongda Rd., Zhongli District, Taoyuan City 320317, Taiwan
       \AND
       \name Chu-Lan Michael Kao \email chulankao@gmail.com \\
       \addr Institute of Statistics\\
       National Yang-Ming Chiao-Tung University \\
       1001 University Road, Hsinchu 30010, Taiwan
       }

\editor{My editor}

\maketitle

\begin{abstract}
Hidden Markov models (HMM) have been widely used by scientists to model stochastic systems: the underlying process is a discrete Markov chain, and the observations are noisy realizations of the underlying process. Determining the number of hidden states for an HMM is a model selection problem which is yet to be satisfactorily solved, especially for the popular Gaussian HMM with heterogeneous covariance. In this paper, we propose a consistent method for determining the number of hidden states of HMM based on the marginal likelihood, which is obtained by integrating out both the parameters and hidden states. Moreover, we show that the model selection problem of HMM includes the order selection problem of finite mixture models as a special case. We give rigorous proof of the consistency of the proposed marginal likelihood method and provide an efficient computation method for practical implementation. We numerically compare the proposed method with the Bayesian information criterion (BIC), demonstrating the effectiveness of the proposed marginal likelihood method.
\end{abstract}


\section{Introduction}

It is well recognized that hidden Markov models (HMM) and general state space models provide useful frameworks for describing noisy observations from an underlying stochastic process. They are popular for processing time series data and widely used in fields like speech recognition, signal processing, and computational molecular biology. 

The basic components of a hidden Markov model include the observations $\{Y_i = y_i, 1\leq i\leq n\}$ and the corresponding hidden states $\{X_i = x_i,1\leq i\leq n\}$, which is a Markov chain. Throughout the paper, we use upper cases $\{Y, X\}$ to denote the random variables and the corresponding lower cases $\{y, x\}$ to denote the realizations (observations). In this article, we consider discrete state space hidden Markov models, i.e., the hidden states $X_i$ have a finite support, observed at discrete time points $\{t_1,\ldots, t_n\}$, or $\{1,\ldots, n\}$ for notation simplicity. The size of the support of the hidden states, denoted by $K$, is the number of hidden states of an HMM. In most real-world problems, the number of hidden states is not known beforehand but conveys important information of the underlying process. For example, in molecular biology, $K$ could be the number of distinct conformations of a protein; in chemistry, $K$ could be the number of distinct chemical species in a biochemical reaction. Existing methods to estimate $K$ either suffer from a lack of theoretical guarantee or unfeasible/impractical implementation, which we review in detail in Section~\ref{subsec:literature_review_model_selection_HMM}. The goal of this article is to provide a consistent method, \textit{the marginal likelihood method}, to determine the number of hidden states $K$ based on the observations $\{y_1,\ldots, y_n\}$ of an HMM, which is computationally feasible for practitioners with minimal tuning.

\subsection{Recap of HMM and Notations}
\label{subsec:recaphmm}

Consider the following hidden Markov model (HMM): let ${\bm X} = \{X_i, i \geq 0\}$ be an ergodic (irreducible, aperiodic, and positive recurrent) Markov chain on a finite state space $\mathcal{X}_K = \{1, \cdots, K \}$ with transition matrix 
$Q_K = \{ q_{k\ell}, 1 \leq k, \ell \leq K\} \in \mathcal{Q}_K$, i.e., $q_{k\ell} = P(X_{i+1}=\ell | X_i = k)$ for all $i \geq 0$, where $\mathcal{Q}_K := \left\{ Q_K: q_{k\ell} > 0, \sum_{k'=1}^K q_{kk'} = 1, \ \forall 1 \leq k,\ell \leq K \right\}$ is the collection of transition matrices with all positive entries. Note that a Markov chain governed by a $Q_K \in \mathcal{Q}_K$ is irreducible, aperiodic, and positive recurrent, so it has an invariant measure $\mu(Q_K)=(\mu_1(Q_K), \cdots, \mu_K(Q_K))$; we further assume that $X_0$ follows $\mu(Q_k)$ so that ${\bm X}$ is stationary. Conditioning on ${\bm X}$, ${\bm Y} = \{Y_i, i \geq 1\}$ are independent random variables on a measurable space $\mathcal{Y}$, and for all $i \geq 1$, when given $X_i = k$, $Y_i$ is assumed to have probability density function $f(\cdot|\boldsymbol{\theta}_k)$ (which is independent to $i$) with respect to some $\sigma$-finite measure $\lambda$ on $\mathcal{Y}$, where $\boldsymbol{\theta}_k \in \Theta$, and $\Theta$ is a subspace of $\mathbb{R}^d$, the $d$-dimensional Euclidean space. We assume that $f$ is distinguishable on $\Theta$, i.e., for all $1 \leq k < \ell \leq K$, $\lambda \{y: f(y|\boldsymbol{\theta}_k) \neq f(y|\boldsymbol{\theta}_\ell)\} > 0$. We denote the model parameters by $\boldsymbol{\phi}_{K} = (Q_K; \boldsymbol{\theta}_1,\ldots, \boldsymbol{\theta}_K) \in \mathcal{Q}_K \times \Theta^K := \Phi_K$.

Suppose we observe ${\bm y}_{1:n} = \{y_1, y_2, \cdots, y_n\} \in \mathcal{Y}^n$, but the underlying process ${\bm x}_{1:n} = \{x_1, x_2, \cdots, x_n\}$ remains hidden (unobserved). The joint likelihood of $({\bm y}_{1:n}, {\bm x}_{1:n})$ given the parameters $\boldsymbol{\phi}_K$ is
\begin{align}\label{eqn_lik}
p({\bm y}_{1:n}, {\bm x}_{1:n}| \boldsymbol{\phi}_{K}) 
:= 
\sum_{k=1}^K \mu_k(Q_K) q_{kx_1} f(y_1|\boldsymbol{\theta}_{x_1}) \times \prod_{i=2}^n q_{x_{i-1}x_i} f(y_i|\boldsymbol{\theta}_{x_i}).
\end{align}
The likelihood after integrating out the hidden states is
\begin{equation}
\label{eqn_full_lik}
p({\bm y}_{1:n}| \boldsymbol{\phi}_K) = \sum_{{\bm x}_{1:n} \in \mathcal{X}_K^{n}} p({\bm y}_{1:n}, {\bm x}_{1:n}|\boldsymbol{\phi}_K),
\end{equation}
where $\mathcal{X}_K^{n}$ denotes the product space of $n$ copies of $\mathcal{X}_K$.

The maximum likelihood estimator (MLE) of a hidden Markov model given $K$, the number of hidden states, can be obtained through the Baum-Welch/ Expectation-Maximization (EM) algorithm { \color{blue}\citep{BW:1966, Baum:1970, DLR:1977}.} Under certain regularity conditions, the consistency and asymptotic normality of the maximum likelihood estimator of HMMs are established in~\cite{LerouxM:1992} and \cite{Bickel:1998}, respectively, when the correct $K$ is specified.

\subsection{Brief Literature Review}
\label{subsec:literature_review_model_selection_HMM}

It has been recognized that the model parameters of an HMM are not identifiable when the number of hidden states is over-estimated~\cite[Chapter~22 of][]{Hamilton:1994,Ferguson:1980,Ryden:1998}. Thus, determining the number of hidden states, also known as the order selection in the machine learning literature, is an important problem for conducting valid inferences on model parameters of hidden Markov models. There is a vast literature on the model selection for hidden Markov models. We briefly review a fraction of the most widely adopted methods here to place our work in the context of literature. 

A special case of HMMs is finite mixture models, where the rows of the transition matrix are identical to each other. The model selection of finite mixture models is mostly based on penalized likelihood, also known as information-theoretic approaches, such as the Akaike Information Criterion (AIC). {A rich literature has been developed for finite mixture models, including \cite{ChenK:1996},\cite{ ChenK:2012}, \cite{Chen:2009}, \cite{ChenT:2009}, \cite{ChenT:2008}, \cite{Tao:2013}, \cite{Hui:2015},  \cite{Jeffries:2003}, \cite{Lo:2001}, \cite{Rousseau:2011}, and many others.}


When the observations $\{y_1,\ldots, y_n\}$ are supported on a finite set (i.e., when they are discrete-valued), we call it a finite-alphabet hidden Markov process~\citep{MacDonald:1997}. Information-theoretic approaches (based on maximum likelihood estimation and penalization terms) for the order estimation of finite-alphabet hidden Markov processes are widely used.~\cite{Finesso:1990} proposes a penalized likelihood method, which is proved to be strongly consistent for finite-alphabet HMMs under certain regularity conditions.~\cite{Ziv:1992} derives the estimator by minimizing the under-estimation probability, which is shown to be not consistent~\citep{Kieffer:1993,LiuN:1994}. ~\cite{LiuN:1994} gives a modified version, which is shown to be consistent given an upper bound of the order of a finite-alphabet HMMs.~\cite{Kieffer:1993} gives a strongly consistent estimator that resembles the Bayesian information criterion (BIC)~\citep{Schwarz:1978} for finite-alphabet HMMs.~\cite{GassiatB:2003} proves strong consistency of these penalized maximum likelihood estimations without assuming any upper bound on the order for finite-alphabet HMMs, with smaller penalties than previous works. See \cite{Ephraim:2002} and \cite{Ryden:1995} for more detailed discussions about the literature on order selection of finite-alphabet HMMs.

However, when the observations $\{y_1,\ldots, y_n\}$ are supported on the real line, as in the Gaussian HMM, where each observation follows a Gaussian distribution conditioning on its hidden state, the problem becomes more difficult. The major difficulty comes from the fact that the overly-fitted mixture models are not identifiable and that the likelihood ratio statistics becomes unbounded, see~\cite{Gassiat2014}. The majority of the methodologies proposed in the literature rely on the idea of penalized likelihood, the consistency of which remains to be satisfactorily solved. {\color{blue}For instance, both AIC and BIC have been proposed to Gaussian HMMs in \cite{Leroux2:1992}, but AIC has been shown to be inconsistent (\cite{fuh2024kullback})\footnote{{\color{blue}It is known that AIC is good for forecasting other than estimation in typical model selection problem.}}, and the existing consistency results for BIC do not apply to hetergoeneous Gaussian HMMs with unequal variance (\cite{Leroux:1992}, \cite{yonekura2021asymptotic}). }

{\color{blue}Other methods using penalized likelihood exist in the literature, including the minimum description length (MDL) in \cite{MDL2}, \cite{MDL3} and \cite {MDL}.  \cite{hung2013hidden} gives a consistent estimator of the number of hidden states using double penalization when assuming that the maximum likelihood estimators are consistent. ~\cite{Ryden:1995} introduces an estimator that does not asymptotically under-estimate the order, given an upper bound for the order. ~\cite{Ryden:1998} applies the bootstrap technique to perform a likelihood ratio test for the order estimation of hidden Markov models for a real-data example.~\cite{GassiatK:2000} investigates the likelihood ratio test for testing a single population i.i.d.~model against a mixture of two populations with Markov regime.~\cite{MacKAY:2002} estimates the order and the parameters together by minimizing a penalized distance function of the empirical distribution with all finite mixture models. Information-theoretic approaches make it possible to add heavier penalties as opposed to that of the BIC; see, for instance, \cite{MDL3}, \cite{gassiat2002likelihood}, and \cite{GassiatB:2003}. All of the above faces the aforementioned obstacles of non-identifiability and unbounded likelihood.}

Bayesian methods, which do not depend on the maximum likelihood estimator, also play an important role in the HMM model selection literature. Reversible jump methods {\color{blue}proposed by \cite{fan2011reversible} and  \cite{green2009reversible} have been successfully adopted in practice by \cite{boys2004bayesian}, \cite{green2002hidden}, \cite{robert2000bayesian}, and \cite{spezia2010bayesian},} with a lack of theoretical justification. \cite{Gassiat2014} provides a frequentist asymptotic evaluation of Bayesian analysis methods purely from a theoretical perspective: under certain conditions on the prior, the posterior concentration rates and a consistent Bayesian estimation of the number of hidden states are given; practical implementation, guidance for tuning of the algorithm, and numerical results are not provided therein. 

\textcolor{blue}{In this work, we adopt the Bayesian approach and show both theoretically and numerically that the marginal likelihood method can give a consistent selection of the number of hidden states.} Some authors have studied approaches that are related to our marginal likelihood method. \cite{chambaz2005nonasymptotic} uses marginal likelihood ratio for the order estimation of mixture models and obtains similar results for the marginal likelihood ratio: $O(e^{-cn})$ for underestimation, and $O(n^{-1/2+\delta})$ for overestimation. \cite{wang2015likelihood} adapts the penalty approach to stochastic block models. Though the aforementioned studies share similarities with the results in this paper, in these studies on mixture models, the hidden state variables are assumed to be independent and identically distributed (i.i.d.), which is not the case for HMMs.

\subsection{Gaussian Hidden Markov Models}
\label{sec:gaussianhmmproblem}

In this section, we discuss the difficulties of the order selection of HMMs using a concrete example widely adopted in applications, the heterogeneous Gaussian HMM.

In a heterogeneous Gaussian HMM, let $X_0=x_0$, and given $X_i = k$, $Y_i$ follows a Gaussian distribution with mean $\mu_k$ and variance $\sigma_k^2$. Thus $\boldsymbol{\phi}_K = (Q_K; \{\mu_k,\sigma_k^2\}_{1\leq k\leq K})$ and the joint likelihood is
\begin{align*}
p({\bm y}_{1:n}, {\bm x}_{1:n}| \boldsymbol{\phi}_{K}) \propto & \prod_{k=1}^K \left\{  \prod_{i: x_i = k} \frac{1}{\sigma_k}\exp\left(-\frac{(y_i - \mu_k)^2}{2\sigma_k^2}\right) \right\}  \times \left\{ \prod_{i=1}^n q_{x_{i-1}{x_i}}\right\}.
\end{align*}
Note that this likelihood is unbounded for some paths ${\bm x}_{1:n}$. For example, consider a path with $x_1=1$ and $x_i \neq 1$ for all $i \geq 2$. If one takes $\mu_1 = y_1$, then as $\sigma_1 \rightarrow 0$, $p({\bm y}_{1:n}, {\bm x}_{1:n}| \boldsymbol{\phi}_{K}) \rightarrow\infty$. Since the full likelihood  $p({\bm y}_{1:n} | \boldsymbol{\phi}_{K})$ sums over all possible paths, such unbounded path always exists, and thus, {\color{blue}the full likelihood becomes unbounded.} 
This can be a serious issue when one overfits an HMM -- the extra component could concentrate on only one single observation with zero variance, which blows up the likelihood. Therefore, methods of model selection for Gaussian HMM based on penalized likelihoods, {\color{blue}such as the Bayesian Information Criterion defined as
\begin{equation*}
\mbox{BIC} := -2\log p({\bm y}_{1:n}|\hat{\boldsymbol{\phi}}_{K}) + K(K+d-1) \log n,
\end{equation*}
}
which requires the consistency of the maximum likelihood estimator, becomes problematic. General consistency results of model selection based on penalized likelihoods have to exclude this case in their required regularity conditions~\cite{Leroux:1992}. Therefore, the BIC, though widely adopted in practice, is theoretically questionable for its validity as a model selection criterion for HMM, {\color{blue}as discussed by \cite{MacDonald:1997} and \cite{Gassiat2014}.} This is the same issue as the unbounded likelihood for heterogeneous Gaussian mixture models~\citep{ChenK:2012}. In fact, Gaussian mixture models can be obtained by setting $q_{ij} = s_j$ for all $i, j\in \{1,2,\ldots, K\}$, where $\{s_j\}_{1\leq j\leq K}$ are the proportions of the mixture components, satisfying $\sum_{j=1}^K s_j = 1$.

Furthermore, as noted in~\cite{Gassiat2014}, for overly fitted HMMs or other finite mixture models, the model parameters become non-identifiable. In an overly fitted HMM, the neighborhood
of the true transition matrix contains transition matrices arbitrarily close to
non-ergodic transition matrices. Adding hard thresholds to entries in the transition matrix does not satisfactorily solve the problem.

\subsection{Outline}

The remainder of the paper has five sections. We first propose the marginal likelihood method for general HMM order selection in Section~\ref{sec:marginal_likelihood_def}. The consistency of the method is presented in Section~\ref{sec:asymptotic_study}, for which we use a special case to illustrate the proof strategy.
In Section~\ref{section:computation}, we describe the computational method, demonstrate the effectiveness of the marginal likelihood method using numerical experiments, and conclude with discussions on choices of hyper-parameters. Section~\ref{subsec:realdata} applies our proposed marginal likelihood method to real data from single-molecule experiments on protein transportation. Section~\ref{section:conclusions} concludes the paper with a summary. The proof of the consistency theory is provided in the Appendix, with additional details presented in the Online Supplement.

\subsection{Contributions}
\label{subsec:contributions}

The major contributions of the paper are as follows. (1) We investigate the marginal likelihood method for HMM order selection, which resolves the difficulties caused by unbounded likelihood by incorporating the prior distribution. (2) The theoretical result on the consistency of our estimator for the number of states is established. (3) The computational algorithm is efficient, robust, and has been tested to work very well. (4) An easy-to-use R package, \texttt{HMMmlselect}, which implements our algorithm, is provided and publicly available at CRAN (https://cran.r-project.org/package=HMMmlselect).

\section{Model Selection via Marginal Likelihood}
\label{sec:marginal_likelihood_def}

As discussed in Section~\ref{subsec:literature_review_model_selection_HMM}, the existing model selection methods for HMM either have no theoretical guarantee or are theoretically justified only for a very restricted family of HMMs that excludes the popular heterogeneous Gaussian HMM. We propose a marginal likelihood method, which directly compares the probability of obtaining the observations under HMMs with different numbers of hidden states, after integrating out both the model parameters and the hidden states. This method, as we will see, is consistent under weak regularity conditions that are satisfied by a wide range of HMMs, including the heterogeneous Gaussian HMM. 

\subsection{Marginal Likelihood Method}
\label{subsec:marginalLM}

Given the number of states $K$, we assume that each $\boldsymbol{\theta}_k$ is independently drawn from the prior distribution $\pi(\boldsymbol{\theta}|\boldsymbol{\alpha})$ and that the transition probabilities $Q_K$ are drawn from the prior distribution $\nu_K(Q_K|\boldsymbol{\beta}_K)$, independent of each $\boldsymbol{\theta}_k$; here $\boldsymbol{\alpha}$ and $\boldsymbol{\beta}_K$ are the hyper-parameters that are assumed to be fixed constants. We let $p_0(\boldsymbol{\phi}_K)$ denote the joint prior: $p_0(\boldsymbol{\phi}_K) = p_0(\boldsymbol{\phi}_K|\boldsymbol{\alpha}, \boldsymbol{\beta}_K) = \nu_K(Q_K|\boldsymbol{\beta}_K) \prod_{k=1}^K \pi(\boldsymbol{\theta}_k|\boldsymbol{\alpha}).$ The marginal likelihood under a $K$-state HMM is then defined as
\begin{equation}
\label{eqn_marg_likelihood}
p_K({\bm y}_{1:n}) = \int_{\Phi_K}  p({\bm y}_{1:n}|\boldsymbol{\phi}_K)\  p_0(\boldsymbol{\phi}_K)\  d\boldsymbol{\phi}_K,
\end{equation}
where $\boldsymbol{\phi}_K\in \Phi_K = \mathcal{Q}_K \times \Theta^K$ is defined in Section~\ref{subsec:recaphmm}.

Given a sufficiently large $\overline{K}$, we choose $K \in \{ 1, 2, \cdots, \overline{K} \}$ that maximizes the marginal likelihood as the estimator of the number of states, i.e.,
\begin{equation}\label{k-hat}
\hat{K}_n := \mbox{arg max}_{1 \leq K \leq \overline{K}} p_K({\bm y}_{1:n}).
\end{equation}
We show that, if the true number of states is $K^*$ and $\overline{K} \geq K^*$, then, under mild conditions, $\hat{K}_n$ is a consistent estimator of $K^*$.

\subsection{Discussions of Marginal Likelihood Methods}

The marginal likelihood has been used in the model selection literature. The ratio of marginal likelihoods is known as the Bayes factor {\color{blue}\citep{BF:1995}}, a popular model selection criterion. The BIC is in fact an approximation of the marginal likelihood using the Laplace method. \cite{Ghahramani:2001} discusses the practical applicability and calculation of the Bayes factor.~\cite{bauwens2014marginal} applies the marginal likelihood method for model selection of Markov-switching GARCH and change-point GARCH models.~\cite{stepwise} uses the marginal likelihood method to determine the number and locations of change points of a stepwise signal.

As discussed in Section~\ref{sec:gaussianhmmproblem}, the heterogeneous Gaussian HMM suffers from having an unbounded likelihood surface. Adding a prior for the variance parameters in Gaussian HMMs can, in fact, fix the issue of unbounded likelihood. Therefore, the proposed marginal likelihood method, which integrates out the parameters and hidden states, does not suffer from irregularity of the likelihood surface. 

{\color{blue}
It is worth mentioning that there are other methods involving the use of prior. For example, \cite{Gassiat2014} sample parameter $\boldsymbol{\phi}_K$ from $\Phi_K$ with $K \geq K^*$, and apply a function on the sampled parameter to obtain a consistent estimator of $K^*$. Since their approach requires sampling from posterior distribution, their algorithm also involve MCMC. The difference between these two methods is that we provide computational approximation that makes the method efficiently enough for practical usage.

Finally, the conditional marginal likelihood has also been proposed in the literature. For instance, \cite{lotfi2022bayesian} has observed the overfitting and underfitting issues for marginal likelihood under a machine learning setting, and therefore proposed the use of conditional marginal likelihood. However, these issues do not apply to our setting as we will prove the consistency of our estimator, which ensures that the probabilities of overfitting and underfitting both go to zero asymptotically.
}

\section[Theoretical Study]{Theoretical Study of the Marginal Likelihood Estimator}
\label{sec:asymptotic_study}
 
We first present the consistency result of the proposed marginal likelihood estimator for HMM order selection, Theorem \ref{theorem:asymptotic_consistency_hmm}, in Section~\ref{subsec:consistencytheorem}. A brief illustration of the proof concept under a special case is then presented in Section~\ref{subsec:illustrate}. Finally, in Section~\ref{subsec:connectwithGM}, we point out the connections between the order selection of HMMs and the model selection of finite mixture models.
The proofs of these results are deferred to the Appendix, starting with a short explanation of how to extend the proof concept to more general cases. An online supplement is further provided for theoretical details.
 
Throughout the paper, we use $\xrightarrow{P}$ to denote convergence in probability under probability law $P$. For a decreasing sequence $\{\epsilon_n, n > 0\}$ that converges to $0$ as $n\rightarrow\infty$, we denote it by $\epsilon_n\downarrow 0$. For any set $\Omega$, we use $\Omega^c$, $\overline{\Omega}$, $\partial \Omega$ and $1_{\Omega}$ to denote its complement, closure, boundary, and corresponding indicator function, respectively. For any vector or matrix $A$, let $A^t$ be its transpose. We use $||\cdot||$ to denote the $L_2$-norm under the corresponding space. Finally, we use $D_{\boldsymbol{\theta}} f$ and $D_{\boldsymbol{\theta}}^2 f$ to denote the gradient vector and Hessian matrix of a function $f$ with respect to $\boldsymbol{\theta}$, respectively.

\subsection{Consistency and Rate of Convergence}
\label{subsec:consistencytheorem}

Before we state $\hat{K}_n$ is a consistent estimator of $K^*$, we need first to define what is the ``true" number of states $K^*$ for a given HMM. Since a $K^*$-state HMM can always be embedded as a $K$-state HMM with $K>K^*$, we basically define $K^*$ as the smallest possible number of states that can characterize the HMM (which turns out to also be the number of states providing the identifiability {\color{blue} under the weakly identifiable condition 3) to be defined later; see \cite{LerouxM:1992}, Lemma 2 for related results.}) That is, for a $K$ state HMM, recall that $\Phi_K$ is the corresponding parameter space, and let $\mathcal{M}_K := \{ p(\cdot|\boldsymbol{\phi}_{K}): \boldsymbol{\phi}_{K} \in \Phi_K\}$, i.e., the set of probability distributions indexed by parameters $\boldsymbol{\phi}_{K} \in \Phi_K$, where $p(\cdot|\boldsymbol{\phi}_{K})$ is defined in (\ref{eqn_full_lik}). Define $K^*$ as the smallest positive integer $K$ such that the probability distribution of $\{Y_i, i \geq 1\}$ is in $\mathcal{M}_K$.
 
Now, let $K^*$ be the true number of states and $\boldsymbol{\phi}^* = (Q^*; \boldsymbol{\theta}_1^*,\ldots, \boldsymbol{\theta}_{K^*}^*)$ be the true parameters of the HMM considered. Let $P^*$ and $E^*$ denote the probability and expectation under the true parameter $\boldsymbol{\phi}^*$, respectively. In addition, by the definition of $\mathcal{Q}_K$, the closure of $\mathcal{Q}_K$ is 
\begin{equation*}
\overline{\mathcal{Q}}_K := \left\{ Q_K: q_{k\ell} \geq 0, \sum_{k'=1}^K q_{kk'} = 1, \ \forall 1 \leq k,\ell \leq K \right\}.
\end{equation*}
Further, denote
\begin{equation*}
\mathcal{Q}_K^\epsilon := \left\{ Q_K: q_{k\ell} \geq \epsilon, \sum_{k'=1}^K q_{kk'} = 1, \ \forall 1 \leq k,\ell \leq K \right\}
\end{equation*}
for any given $\epsilon > 0$.

We assume the following regularity conditions in our consistency theorems.
\begin{enumerate}[1)]
    \item $q_{k\ell}^* > 0$ for all $1 \leq k, \ell \leq K^*$. This implies that the Markov chain with transition matrix $Q^*$ is irreducible, aperiodic, and positive recurrent.
    \item $\Theta$ is a compact set in $\mathbb{R}^d$, and the true parameters $\{\boldsymbol{\theta}_k^*, k = 1, \cdots, K^*\}$ are distinct interior points of $\Theta$.
    \item  For any $k \in \{1, 2, \cdots, K^*\}$, for any $\boldsymbol{\theta} \neq \boldsymbol{\theta}_k^*$, we have $\lambda \{y: f(y|\boldsymbol{\theta}) \neq f(y|\boldsymbol{\theta}_k^*)\} > 0$ and $E^* | \log f(Y_1|\boldsymbol{\theta}_k^*)| < \infty$. In addition, for any $K>0$, the family of mixtures of at most $K$ elements of $\{ f(\cdot | \boldsymbol{\theta}) : \boldsymbol{\theta} \in \Theta \}$ is weakly identifiable in the sense that, for any positive $q_k$ and $q_k'$ with $\sum_{k=1}^K q_k = \sum_{k=1}^K q_k' = 1$,
    \begin{equation*}
    \lambda \left\{ y: \sum_{k=1}^K q_k f(y| \boldsymbol{\theta}_k) \neq  
        \sum_{k=1}^K q_k' f(y | \boldsymbol{\theta}_k') \right\} = 0
    \end{equation*}
    if and only if
    \begin{equation*}
        \sum_{k=1}^K q_k 1_{\boldsymbol{\theta}_k}(\boldsymbol{\theta}) \equiv 
        \sum_{k=1}^K q_k' 1_{\boldsymbol{\theta}_k'}(\boldsymbol{\theta})
    \end{equation*}
    as a function of $\boldsymbol{\theta}$.

    \item There exists $\delta > 0$ such that, for any $k \in \{1, 2, \cdots, K^*\}$, the function $\boldsymbol{\theta} \rightarrow f(\cdot|\boldsymbol{\theta})$ is twice continuously differentiable in $\mathcal{B}_\delta(\boldsymbol{\theta}_k^*) := \{ \theta: \Vert \boldsymbol{\theta} - \boldsymbol{\theta}_k^* \Vert < \delta \}$. 
    Furthermore, the following holds:
    \begin{equation*}
    E^* \left[\sup_{\boldsymbol{\theta} \in \mathcal{B}_\delta(\boldsymbol{\theta}_k^*) } || D_{\boldsymbol{\theta}} \log f(Y_1|\boldsymbol{\theta}) || \right] < \infty,
    \end{equation*}
    \begin{equation*}
    E^* \left[\sup_{\boldsymbol{\theta} \in \mathcal{B}_\delta(\boldsymbol{\theta}_k^*) } || D_{\boldsymbol{\theta}}^2 \log f(Y_1|\boldsymbol{\theta}) || \right] < \infty,
    \end{equation*}
    \begin{equation*}
    \int_\mathcal{Y} \sup_{\boldsymbol{\theta} \in \mathcal{B}_\delta(\boldsymbol{\theta}_k^*) } || D_{\boldsymbol{\theta}}  f(y|\boldsymbol{\theta}) || \lambda(dy) < \infty,
    \end{equation*}
    \begin{equation*}
    \int_\mathcal{Y} \sup_{\boldsymbol{\theta} \in \mathcal{B}_\delta(\boldsymbol{\theta}_k^*) } || D_{\boldsymbol{\theta}}^2  f(y|\boldsymbol{\theta}) || \lambda(dy) < \infty.
    \end{equation*}
    In addition, for each $k \in \{ 1, 2, \cdots, K^* \}$,
    \begin{equation*}
    P^* \left\{ \sup_{ \boldsymbol{\theta}, \boldsymbol{\theta}' \in \bigcup_{k=1}^{K^*}\mathcal{B}_\delta(\boldsymbol{\theta}_k^*)} \frac{f(Y_1|\boldsymbol{\theta})}{f(Y_1|\boldsymbol{\theta}')} = \infty \Bigg\vert X_1 = k \right\} < 1.
    \end{equation*}

    \item There exists $\delta > 0$ such that, for any $\boldsymbol{\theta} \in \Theta$,
    \begin{equation*}
    E^* \left[ \sup_{\boldsymbol{\theta}' \in \mathcal{B}_\delta(\boldsymbol{\theta})} (\log f(Y_1|\boldsymbol{\theta}'))^+ \right] < \infty.
    \end{equation*}

    \item The prior density $\pi(\boldsymbol{\theta} |\boldsymbol{\alpha})$ is a continuous function of $\boldsymbol{\theta}$ on $\Theta$, and is positive at $\boldsymbol{\theta}_k^*$ for all $1\leqslant k \leqslant K^*$. 
    
    \item For all $K $, the prior density $\nu_K(Q_K|\boldsymbol{\beta}_K)$ is a continuous and positive function of $Q_K$ on $\mathcal{Q}_K$, with a support in $\mathcal{Q}_K^\epsilon$ for some $\epsilon > 0$. In addition, $\nu_{K^*}$ is positive at $Q^*$.
\end{enumerate}

Under these assumptions, we have the following consistency result.

\begin{theorem}\label{theorem:asymptotic_consistency_hmm}
Assume that conditions 1)-7) hold. 
Then, under $P^*$, as $n \rightarrow \infty$,
\begin{itemize}
    \item[1)] there exists $c>0$ such that, for all $K < K^*$,
    \begin{equation}\label{eqn_thm_consist_under}
    \frac{p_K({\bm y}_{1:n})}{p_{K^*}({\bm y}_{1:n})} = O_{P^*} \left( e^{-cn} \right);
\end{equation}
    \item[2)] for all $K > K^*$,
\begin{equation}\label{eqn_thm_consist}
\frac{p_K({\bm y}_{1:n})}{p_{K^*}({\bm y}_{1:n})} = O_{P^*} \left( n^{-\frac{(K-K^*)d}{2}} \right).
\end{equation}
\end{itemize}
Consequently, if $\overline{K} \geq K^*$, i.e., the upper bound is at least $K^*$, then $\hat{K}_n \rightarrow K^*$ in probability as $n \rightarrow \infty$.
\end{theorem}

\begin{remark}
Conditions 1)-7) include the conditions from \cite{HMMCLT}, which ensure the asymptotic normality for the posterior distribution under the true number of states. An additional weakly identifiable condition is enforced in the second part of condition 3) to deal with the case when having misspecified number of states. This weakly identifiable condition is the same as condition (id) in \cite{keribin2000consistent} and implies condition 2 in \cite{LerouxM:1992}.
\end{remark}

\begin{remark}\label{remark:prior}
Note that in Condition 7), we assume that the prior has support in $\mathcal{Q}_K^\epsilon$. This is the strong mixing condition typically used in the literature. See Assumption 3 in \cite{yonekura2021asymptotic}, for instance. However, as suggested by \cite{Gassiat2014},  
it is possible that having a prior vanishes quickly enough as it approaches the boundary of $\mathcal{Q}_K$ would give the same result. The simulation studies in Section \ref{section:computation} show that it is indeed possible. See Online Supplement, Remark \ref{remark:prior-conj} for more details, as well as Remarks \ref{remark-Lemma:L} and \ref{remark-Lemma:L''} for potential issues.
\end{remark}

\begin{remark}\label{remark:Theta}
One can relax the assumption that $\Theta$ is compact to ``if for each $y \in \mathcal{Y}$, $f(y|\cdot)$ vanishes at the infinity, and $\pi$ vanishes at the infinity''. See \cite{LerouxM:1992}, page $130$. {\color{blue}After such relaxation, conditions 1)-7) hold for a wide range of commonly used HMMs, including the aforementioned Gaussian HMMS.} The simulation studies in Section \ref{sec:simulation_study} further confirm this claim.
\end{remark}

{\color{blue}
\begin{remark}\label{remark:Qk}
It should be pointed out that we only assume the compactness on $\Theta$, not $\mathcal{Q}_K \times \Theta^K$. This means that we can have $q_{k\ell}$ arbitrarily close to $0$. However, as we impose condition 7), the prior $\nu_K$ is zero outside of $\mathcal{Q}_K^\epsilon$, and $\mathcal{Q}_K^\epsilon \times \Theta^K$ is compact. In other words, condition 7) is effectively equivalent to assuming the entire parameter space is compact, just as assumed in \cite{douc2009subgeometric}, \cite{douc2011consistency}, and others. We choose to present it in the form of condition 7) since it has the potential to be weakened, as mentioned in Remark \ref{remark:Theta}.
\end{remark}
}

\subsection{Illustration of the Proof Strategy}\label{subsec:illustrate}

To illustrate why Theorem \ref{theorem:asymptotic_consistency_hmm} holds and the idea of its proof, let us consider a special case with $K^* = 1$ and $K = 2$. In this case, we have $\{Y_i, i \geq 1\}$ coming from a single-state HMM with true parameter $\boldsymbol{\phi}^* = (\{1\}; \boldsymbol{\theta}_1)$, which is equivalent to an i.i.d. model with density $f(\cdot|\boldsymbol{\theta}_1)$. We want to show that
\begin{equation}\label{special-case-consist}
\frac{p_2({\bm y}_{1:n})}{p_1({\bm y}_{1:n})} 
\xrightarrow[n \rightarrow \infty]{P^*} 0.
\end{equation}

To approach \eqref{special-case-consist}, we first examine $p_1({\bm y}_{1:n})$, the marginal likelihood of a single-state HMM. Note that for a single-state HMM, the only possible transition matrix is $\{1\}$, so its parameter space $\Phi_1 = \mathcal{Q}_1 \otimes \Theta$ is $d$-dimensional. Also note that the marginal likelihood is the denominator of the posterior distribution, so by Bernstein-von Mises theorem, we have
\begin{equation}\label{special-case-1}
    n^{d/2}\frac{p_1({\bm y_{1:n}})}{p_1({\bm y_{1:n}}|\boldsymbol{\phi}^*)}
=
    \left(
    \frac{p_1({\bm y_{1:n}}|\boldsymbol{\phi}^*)}
    {n^{d/2}p_1({\bm y_{1:n}})}
    \right)^{-1}
    \xrightarrow[n \rightarrow \infty]{P^*} c_1
\end{equation}
for some constant $c_1 > 0$.

We then examine $p_2({\bm y}_{1:n})$, the marginal likelihood of a two-state HMM. Note that $\mathcal{Q}_2$ is two dimensional, so $\Phi_2 = \mathcal{Q}_2 \otimes \Theta \otimes \Theta$ is $2(d+1)$-dimensional. Consider
\begin{equation*}
    \tilde{\boldsymbol{\phi}}^* := \left( 
    \begin{pmatrix}
    1/2 & 1/2 \\
    1/2 & 1/2 
    \end{pmatrix}
    ; \boldsymbol{\theta}_1^*, \boldsymbol{\theta}_1^*
    \right),
\end{equation*}
then, a direct check shows that
\begin{equation}\label{two-equals-one}
    p_2( \cdot |\boldsymbol{\phi}) \equiv p_1( \cdot |\boldsymbol{\phi}^*)
\end{equation}
holds when $\boldsymbol{\phi} = \tilde{\boldsymbol{\phi}}^*$. This means that we can treat $\{Y_i, i \geq 1 \}$ coming from a two-state HMM with true parameter $\tilde{\boldsymbol{\phi}}^*$. It is therefore ``tempted'' to apply the Bernstein-von Mises theorem again to get
\begin{equation}\label{special-case-2}
    n^{2(d+1)/2}\frac{p_2({\bm y_{1:n}})}{p_2({\bm y_{1:n}}|\tilde{\boldsymbol{\phi}}^*)}
    \xrightarrow[n \rightarrow \infty]{P^*} c_2
\end{equation}
for some constant $c_2 > 0$. If we can do this, then combing \eqref{special-case-1} and \eqref{special-case-2} gives
\begin{equation*}
n^{(d+2)/2}\frac{p_2({\bm y_{1:n}})}{p_1({\bm y_{1:n}})}
=
\frac{
    n^{2(d+1)/2}\frac{p_2({\bm y_{1:n}})}{p_2({\bm y_{1:n}}|\tilde{\boldsymbol{\phi}}^*)}
}{
    n^{d/2}\frac{p_1({\bm y_{1:n}})}{p_1({\bm y_{1:n}}|\boldsymbol{\phi}^*)}
}
    \xrightarrow[n \rightarrow \infty]{P^*} \frac{c_2}{c_1},
\end{equation*}
which immediately leads to \eqref{special-case-consist}.

Unfortunately, the argument \textit{does not hold}. The reason is that the argument above uses Bernstein-von Mises theorem, which requires that \eqref{two-equals-one} holds \textit{only} at $\tilde{\boldsymbol{\phi}^*}$; in other words, it requires the \textit{identifiability}. However, this is not the case. To see why, for any two-state transition matrix $Q \in \mathcal{Q}_2$, consider
\begin{equation*}
    \boldsymbol{\phi}_Q^* := (Q; \boldsymbol{\theta}_1^*, \boldsymbol{\theta}_1^*),
\end{equation*}
then, a direct check shows that \eqref{two-equals-one} holds for all $\boldsymbol{\phi}_Q^*$ , meaning that there exists infinitely many $\boldsymbol{\phi}$ that makes \eqref{two-equals-one} hold. In other words, we face \textit{non-identifiability} here, which forbids us to directly apply the Bernstein-von Mises theorem.

To get around this issue, for any $Q \in \mathcal{Q}_2$, consider the sub-parameter space
\begin{equation*}
    \Phi_Q = \{ \boldsymbol{\phi} = (Q; \boldsymbol{\theta}_1, \boldsymbol{\theta}_2): \boldsymbol{\theta}_1 \in \Theta, \boldsymbol{\theta}_2 \in \Theta\}.
\end{equation*}
Note that we can rewrite $p_2({\bm y}_{1:n})$ as
\begin{equation}\label{special-rep}
    p_2({\bm y}_{1:n})
= \int_{\Phi_2} p_2({\bm y}_{1:n}|\boldsymbol{\phi})  p_0(\boldsymbol{\phi}) d\boldsymbol{\phi}
= \int_{\mathcal{Q}_2} \bigg\{ \int_{\Phi_Q} p_2({\bm y}_{1:n}|\boldsymbol{\phi})  p_0(\boldsymbol{\phi}) d\boldsymbol{\phi} \bigg\} dQ.
\end{equation}
In addition, for any \textit{fixed} $Q \in \mathcal{Q}_2$, note that in the $2d$-dimensional set $\Phi_Q = \{Q\} \otimes \Theta \otimes \Theta$, the only $\boldsymbol{\phi}$ that makes \eqref{two-equals-one} hold is $\boldsymbol{\phi}_Q^*$; in other words, the model is identifiable on the subset $\Phi_Q$. Hence, we can apply the Berstein-von Mises theorem on the posterior distribution restricted to $\Phi_Q$, which gives
\begin{equation}\label{special-case-sub}
    \frac{n^{2d/2}\int_{\Phi_Q} p_2({\bm y_{1:n}}|\boldsymbol{\phi}) p_0(\boldsymbol{\phi}) d \boldsymbol{\phi}}{p_1({\bm y_{1:n}}|\boldsymbol{\phi}^*)}
=
    \frac{n^{2d/2}\int_{\Phi_Q} p_2({\bm y_{1:n}}|\boldsymbol{\phi}) p_0(\boldsymbol{\phi}) d \boldsymbol{\phi}}{p_2({\bm y_{1:n}}|\boldsymbol{\phi}_Q^*)}
    \xrightarrow[n \rightarrow \infty]{P^*} c_Q
\end{equation}
for some constant $c_Q \geq 0$. We can further show that the convergence in \eqref{special-case-sub} is uniform over $Q \in \mathcal{Q}_2$, and $\int_{\mathcal{Q}_2} c_Q dQ < \infty$, so we can combine \eqref{special-rep} and \eqref{special-case-sub} to get
\begin{equation}\label{special-rep-2}
    n^{2d/2} \frac{p_2({\bm y}_{1:n})}{p_1({\bm y_{1:n}}|\boldsymbol{\phi}^*)}
= \int_{\mathcal{Q}_2} \left\{ \frac{n^{2d/2} \int_{\Phi_Q} p_2({\bm y}_{1:n}|\boldsymbol{\phi})  p_0(\boldsymbol{\phi}) d \boldsymbol{\phi}}{p_1({\bm y_{1:n}}|\boldsymbol{\phi}^*)} d\boldsymbol{\phi} \right\} dQ
\xrightarrow[n \rightarrow \infty]{P^*} \int_{\mathcal{Q}_2}c_Q dQ.
\end{equation}
Combining \eqref{special-case-1} and \eqref{special-rep-2}, we get
\begin{equation*}
n^{d/2}\frac{p_2({\bm y_{1:n}})}{p_1({\bm y_{1:n}})}
=
\frac{
    n^{2d/2}\frac{p_2({\bm y_{1:n}})}{p_2({\bm y_{1:n}}|\tilde{\boldsymbol{\phi}}^*)}
}{
    n^{d/2}\frac{p_1({\bm y_{1:n}})}{p_1({\bm y_{1:n}}|\boldsymbol{\phi}^*)}
}
    \xrightarrow[n \rightarrow \infty]{P^*} \frac{\int_{\mathcal{Q}_2} c_Q dQ}{c_1},
\end{equation*}
which immediately leads to \eqref{special-case-consist} and completes the argument.

The above argument basically involves four key steps: (i) decompose $\Phi_K$ into a family of subspaces; (ii) show that each of these subspaces is ``identifiable'' so that the Berstein-von Mises theorem can be applied; (iii) show that the aggregate limit ($\int_{\mathcal{Q}_2} c_2 dQ$ in this case) is finite;
and (iv) show that the convergence in ii) is uniform across all subspaces.
The idea of proof for general $K^*$ and $K$ is essentially the same, though the execution is much more complicated, mainly due to the complication required to decompose $\Phi_K$ into identifiable subspaces. See the Appendix for details.

\begin{remark}
One can conceptually interpret the above argument  as follows:
\begin{itemize}
\item For a single-state HMM, the posterior distribution shrinks toward the true parameter $\boldsymbol{\phi}^*$ in a $d$-dimensional space $\Phi_1 = \{1\} \otimes \Theta$, as we need $\boldsymbol{\theta}_1$ converges to $\boldsymbol{\theta}_1^*$. This is shown in \eqref{special-case-1}.
\item For a two-state HMM, even if we forfeit all the degree of freedom in $\mathcal{Q}_2$ by restricting the parameter space to $\Phi_Q$, we still have the posterior distribution shrinks toward the ``true'' parameter $\boldsymbol{\phi}_Q^*$ in a $2d$-dimensional space $\Phi_Q = \{Q\} \otimes \Theta \otimes \Theta$, as we need both $\boldsymbol{\theta}_1$ and $\boldsymbol{\theta}_2$ converge to $\boldsymbol{\theta}_1^*$. This is shown in \eqref{special-case-sub}.
\end{itemize}
In other words, even if we forfeit all additional dimensions in $\mathcal{Q}_2$ compared to $\mathcal{Q}_1$ in order to obtain identifiability, the resulting subspace still has a larger dimension compared to $\Phi_1$, so its marginal likelihood goes to zero faster due to the existence of additional states. This is the conceptual reason why  \eqref{eqn_thm_consist} holds. The same reason holds for general $K > K^*$, as one can see in the Appendix.
\end{remark}

\subsection{Connections with Model Selection of Mixture Models}
\label{subsec:connectwithGM}

In this section, we discuss the connections of the order selection for HMMs with the model selection of mixture models. As mentioned in Section~\ref{sec:gaussianhmmproblem}, the mixture model can be considered as a special case of an HMM, of which the transition matrix has identical rows, i.e. $q_{ij} = s_j$ for all $i,j=1,\cdots K$, where $s_j$ are the proportions of mixture components satisfying $\sum_{j=1}^K s_j = 1$. Consequently,
the model selection of mixture models can follow the same procedure as the order selection for HMMs. Reversely, we can use the model selection of mixture models to determine the order of HMMs. We show that the estimator of the order of an HMM is still consistent if we ``ignore'' the Markov dependency, i.e., regarding it as a mixture model. 

\begin{corollary}
\label{theorem:asymptotic_consistency_gm}
Assume that all the conditions in Theorem \ref{theorem:asymptotic_consistency_hmm} hold, and we restrict $\nu_K(\cdot|\beta_K)$ to be supported on 
\begin{equation*}
\mathcal{Q}_K^{mix} = \{Q_K: q_{1k}=q_{2k}=\cdots=q_{Kk} \mbox{ for all } 1 \leq k \leq K\},
\end{equation*}
i.e., a prior for a finite mixture model without state dependency. Then the result of Theorem \ref{theorem:asymptotic_consistency_hmm} still holds.
\end{corollary}

\begin{remark}
As opposed to Theorem~\ref{theorem:asymptotic_consistency_hmm}, the computational cost required by
Corollary~\ref{theorem:asymptotic_consistency_gm} is much smaller: instead of fitting HMMs, we only need to fit mixture models, which are in lower dimensional spaces with nice independent structures on the latent variables. In both Theorem~\ref{theorem:asymptotic_consistency_hmm} and Corollary~\ref{theorem:asymptotic_consistency_gm}, the convergence rate of the marginal likelihood ratio is $O_{P^*}(n^{-(K-K^*)d/2})$. However, Corollary \ref{theorem:asymptotic_consistency_gm} requires $n$ to be large so that ${\bm y}_{1:n}$ shows a ``mixture model'' behavior through \textit{stability convergence}; see Online Supplement \ref{appendix:proofofconsistencygm}. This leads to a larger constant term in front of $n^{-(K-K^*)d/2}$ for the marginal likelihood ratio of Corollary \ref{theorem:asymptotic_consistency_gm} as compared to Theorem \ref{theorem:asymptotic_consistency_hmm}, especially for the case of nearly diagonal transition matrices.

\end{remark}

\section[Computation]{Computation and Numerical Experiments}
\label{section:computation}

Now we introduce our method of estimating the marginal likelihood and provide numerical results comparing the marginal likelihood method and the BIC. We conclude this section with a brief discussion about the choice of priors.

\subsection{Marginal Likelihood as a Normalizing Constant}

Denote the joint distribution of ${\bm y}_{1:n}$ and $\boldsymbol{\phi}_K$ by $p({\bm y}_{1:n}, \boldsymbol{\phi}_K) = p({\bm y}_{1:n}| \boldsymbol{\phi}_K)\  p_0(\boldsymbol{\phi}_K)$, where 
$p({\bm y}_{1:n}| \boldsymbol{\phi}_K)$, defined in equation (\ref{eqn_full_lik}), is the likelihood, to which we integrate out the hidden states. The marginal likelihood $p_K({\bm y}_{1:n})$ of a $K$-state HMM is 
$\int_{\Phi_K}\ p({\bm y}_{1:n}| \boldsymbol{\phi}_K)\ p_0(\boldsymbol{\phi}_K) \ d\boldsymbol{\phi}_K$ as defined in (\ref{eqn_marg_likelihood}) 

Our strategy to estimate the marginal likelihood is based on the observation that $p_K({\bm y}_{1:n})$, in fact, can be regarded as the normalizing constant of the posterior distribution $p(\boldsymbol{\phi}_K|{\bm y}_{1:n})  = p({\bm y}_{1:n}, \boldsymbol{\phi}_K) / p_K({\bm y}_{1:n})$.
Thus, the estimation of $p_K({\bm y}_{1:n})$ can be recast as the estimation of normalizing constant of this posterior density. 
To do this, note that we can obtain posterior samples from $p(\boldsymbol{\phi}_K|{\bm y}_{1:n})$ using a Markov chain Monte Carlo (MCMC) algorithm (see~\cite{liu2001} and references therein), since the un-normalized posterior likelihood $p({\bm y}_{1:n}, \boldsymbol{\phi}_K)$ can be evaluated at any $\boldsymbol{\phi}_K$ using the forward algorithm {\color{blue}which integrates out the hidden states ~\citep{BW:1966, Baum:1970,xuan2001algorithms}.} Alternatively, we can sample from the augmented space $\Phi_K\times\mathcal{X}_K^n$, i.e., sample model parameters and the hidden states together till convergence. This alternative  approach corresponds to the data augmentation method in~\cite{TW:1987} and has been used for HMM model fitting {\color{blue}by \cite{Ryden:2008}.} Given that we can sample from the posterior distribution, the question becomes: how to estimate the normalizing constant based on the (posterior) samples.

\subsection{Estimation Procedure}\label{subsubsec:EstimationProcedure}

There is a large literature on the estimation of normalization constants, see~\cite{DiCiccio:1997} and references therein. Among them, we want to mention a few that is related to the algorithm we adopt. Methods based on importance sampling and reciprocal importance sampling require knowledge of a ``good'' importance function whose region of interest covers that of the joint posterior to be integrated {\color{blue}\citep{ChenS:1997, GelfandD:1994, Geweke:1989,OhB:1993,  Steele:2006}.}

The importance sampling and reciprocal importance sampling are simple and fast ways of estimating the normalizing constant if an importance function close to the target density can be specified. Since we already have posterior samples from the unnormalized density, it can be utilized as a guide for choosing a good importance function for either the importance sampling or the reciprocal importance sampling. Therefore, our strategy is to use the importance sampling or the reciprocal importance sampling to estimate the normalizing constant $p_K({\bm y}_{1:n})$, where the importance function is chosen based on the posterior sample from $p(\boldsymbol{\phi}_K | {\bm y}_{1:n})$. Since the posterior samples do not necessarily give enough information about the tail of the posterior distribution, the importance function might be a poor approximation of the target posterior distribution in the tail region, which can result in unstable estimators. We use the locally restricted importance sampling or reciprocal importance sampling, which is more robust to the tail behavior of the target posterior distribution $p(\boldsymbol{\phi}_K|{\bm y}_{1:n})$, see~\cite{DiCiccio:1997}. We now give our procedure for estimating the marginal likelihood $p_K({\bm y}_{1:n})$ for each $K \in \{1,2,\ldots,\overline{K}\}$.

\begin{itemize}
\item[1.] {Obtain posterior samples.} Sample from $p(\boldsymbol{\phi}_K |{\bm y}_{1:n})$ using a preferred MCMC algorithm, and denote the samples by $\{\boldsymbol{\phi}_K^{(i)}\}_{i=1}^{N}$ (where $N$ is often a few thousand). \label{step1:postsampling}
\item[2.] {Find a ``good'' importance function.} Fit a mixture model using the samples $\{\boldsymbol{\phi}_K^{(i)}\}_{i=1}^{N}$, where the number of mixing components is given by either (a) any clustering algorithm, or (b) a pre-fixed number which is large enough. Construct the importance function $g(\cdot)$ by fitting a Gaussian mixture or using a heavier-tailed density as the mixture component; for example, a student-$t$ distribution with a small degree of freedom, such as $2$ or $3$, with the same location and scale parameters as the fitted Gaussian mixture components. 
\label{algorithm:step_fit_gm}
\item[3.] {Choose a finite region.} Choose $\Omega_{K}$ to be a bounded subset of the parameter space such that $1/2 < \int_{\Omega_{K}} g( \boldsymbol{\phi}_K ) d\boldsymbol{\phi}_K < 1$. This can be achieved by finding an appropriate finite region for each mixing component of $g(\cdot)$, avoiding the tail parts.
\item[4.] Estimate 
$p_K({\bm y}_{1:n})$ using either way as follows:
\begin{itemize}
\item {Reciprocal importance sampling.} Approximate $p_K({\bm y}_{1:n})$ by
\begin{align}
\label{eqn:reciprocal_importance_weighting_marginal_likelihood}
\hat{p}_K^{(RIS)}({\bm y}_{1:n}) =  \left[ \frac{1}{N \int_{\Omega_{K}} g( \boldsymbol{\phi}_K ) d\boldsymbol{\phi}_K}\sum_{i = 1}^{N} \frac{g(\boldsymbol{\phi}_K^{(i)})}{p({\bm y}_{1:n}, \boldsymbol{\phi}_K^{(i)})} 1_{\Omega_K}(\phi_K^{(i)})\right]^{-1},
\end{align}
where $1_\Omega(x)$ is $1$ if $x \in \Omega$, and zero otherwise.
\item {Importance sampling.}
\begin{itemize}
\item[(a)] Draw $M$ independent samples from $g(\cdot)$, denoted by $\{\boldsymbol{\psi}_{K}^{(j)}\}_{1\leq j\leq M}$. \item[(b)] Approximate $p_K({\bm y}_{1:n})$ by
\begin{equation}
\label{eqn:importance_weighting_marginal_likelihood}
\hat{p}_K^{(IS)}({\bm y}_{1:n}) = \frac{1}{M P_{\Omega}}\sum_{j = 1}^{M} \frac{p({\bm y}_{1:n},  \boldsymbol{\psi}_K^{(j)})}{g(\boldsymbol{\psi}_K^{(j)})} 1_{\Omega_K}(\psi_K^{(i)}),
\end{equation}
where $P_{\Omega} = \#\mathcal{S} / N$ with $\mathcal{S} = \{i: \boldsymbol{\phi}_K^{(i)}\in \Omega_{K}; 1\leq i\leq N\}$.
\end{itemize}
\end{itemize}
\end{itemize}

The purpose of step 2 is to construct a reasonable importance function that covers the mode of the target density $p(\boldsymbol{\phi}_K | {\bm y}_{1:n})$. Thus, the clustering algorithm, if adopted, does not need to be ``optimal'' in any sense. Therefore, a conservative recommendation is to choose overly-fitted Gaussian (or student-$t$) mixtures based on the posterior samples obtained in step 1. Moreover, the heavy-tailed distribution and the truncated regions both serve the purpose of obtaining a robust importance sampling estimator. If reciprocal importance sampling is used, a heavy-tailed distribution is not recommended for the sake of estimation robustness. 

Simulation studies of various target densities (skewed, heavy-tailed, and high-dimensional) with known normalizing constants validate the efficacy of the proposed procedure, regardless of the shape of the target density or the dimension of the parameter space. See Online Supplement \ref{appendix:simulations_normalizing_constant}.

\subsection{Simulation Studies for HMM Order Selection}
\label{sec:simulation_study}

\begin{figure}[tbph]
\centering
\includegraphics[width=0.45\textwidth]{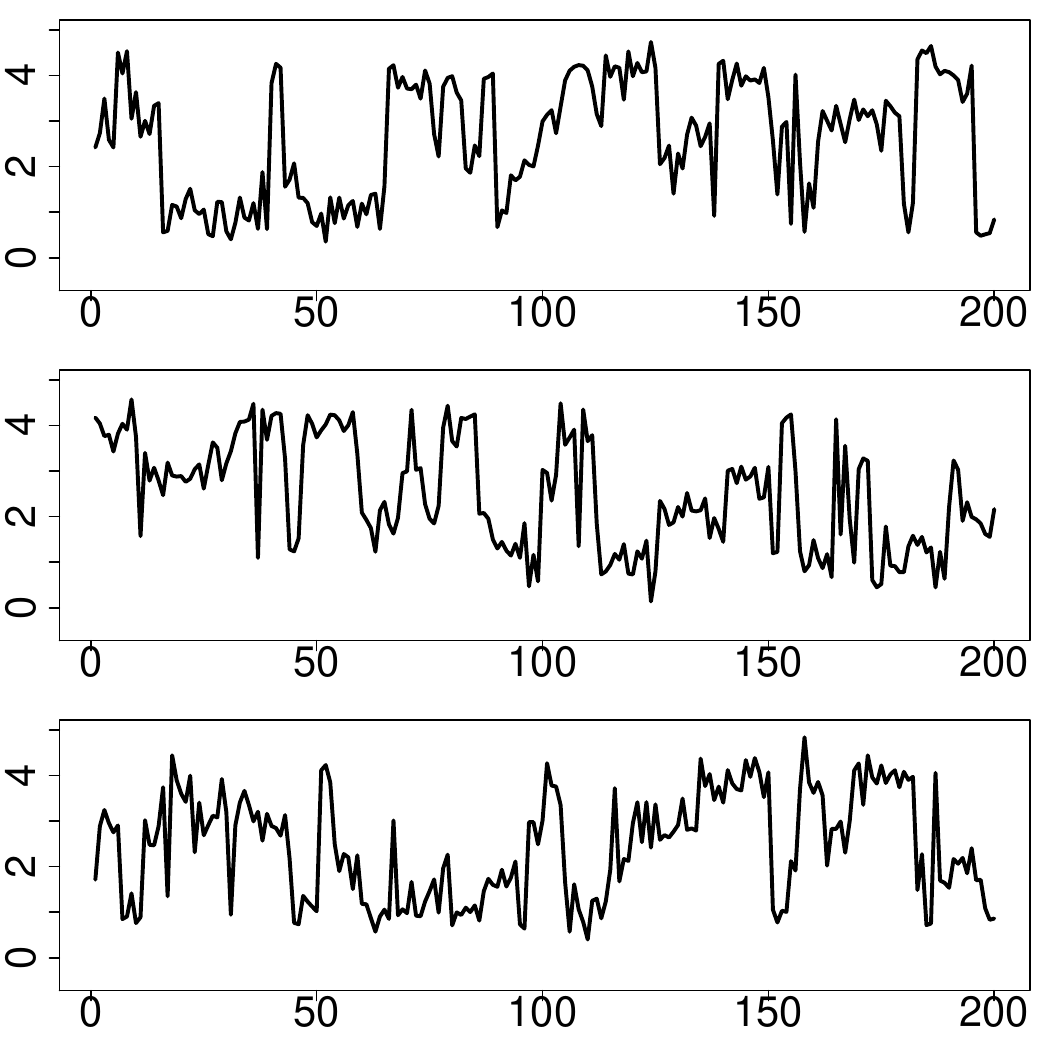}
\includegraphics[width=0.45\textwidth]{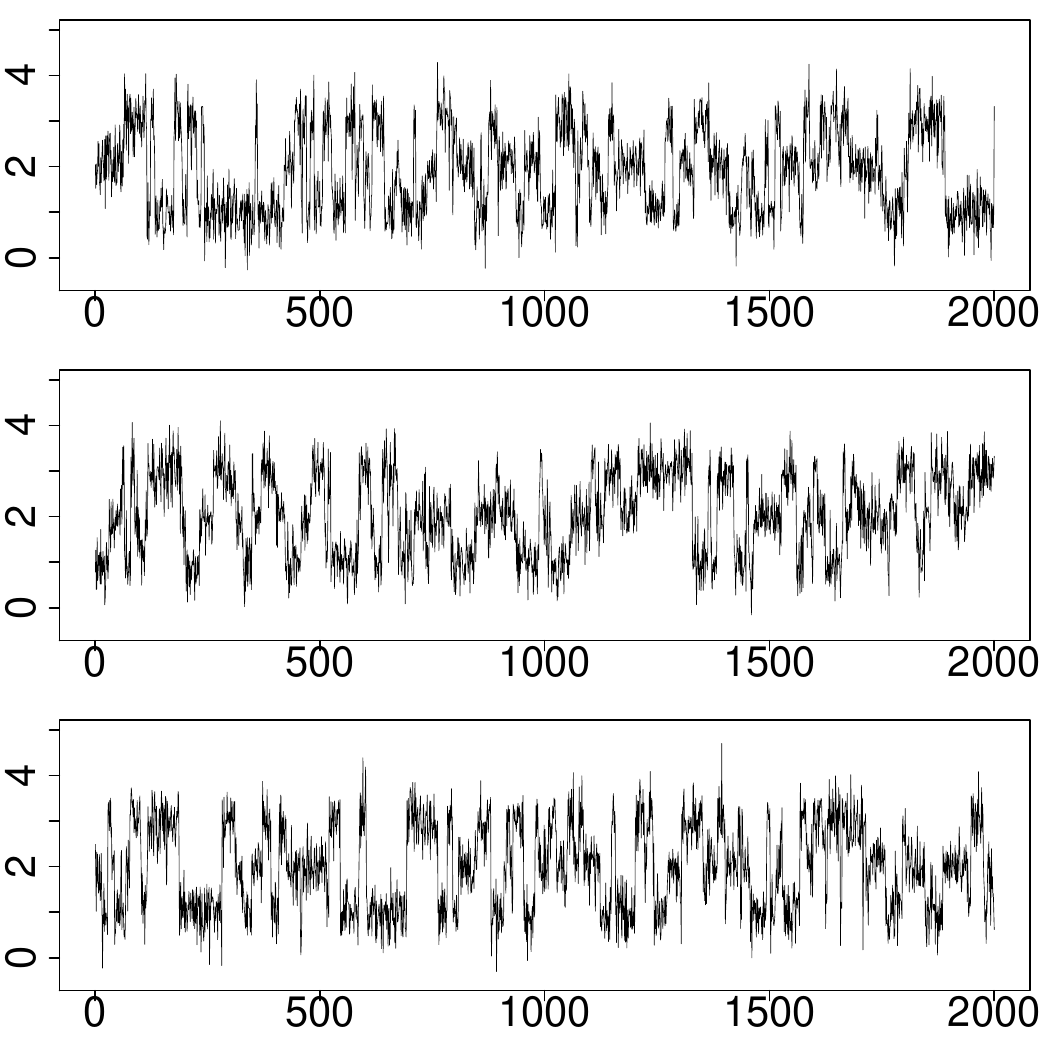}
\caption[Sample HMM traces.]{
\footnotesize{
Sample HMM traces. The top three panels show three simulated HMM traces with $n=200$ observations and $K=4$ hidden states: $\sigma = 0.3$, and the transition matrix is $P_4^{(2)}$. The bottom three panels show three simulated HMM traces with $n=2,000$ observations and $K=3$ hidden states: $\sigma = 0.4$ and the transition matrix is $P_3^{(3)}$.}
\label{fig:demonsimulatedtraces}}
\end{figure}

In the numerical experiments, we fix the mean parameters of a $K$-state HMM to be $\boldsymbol{\mu} = (1, 2,\ldots, K)$ and vary the variances $\boldsymbol{\sigma}^2 = (\sigma^2,\ldots, \sigma^2)$ {\color{blue}in the first set of simulation studies.} The equal variances assumption is adopted here for simplicity of the presentation of the results, but this is not part of the model assumptions. {\color{blue}We will relax this assumption in the second set of simulations.} We consider four kinds of transition matrices, corresponding to flat ($P_K^{(1)}$), moderate and strongly diagonal ($P_K^{(2)}, P_K^{(3)}$) and strongly off-diagonal ($P_K^{(4)}$) cases: 

\begin{align}
\label{eqn:sim_hmm_transmat}
P_K^{(1)} &= \frac{1}{K} E_{K},\ P_K^{(2)} = \left[0.8 - \frac{0.2}{K-1}\right] I_{K} + \frac{0.2}{K-1} E_{K},\\
\label{eqn:sim_hmm_transmat1}
P_K^{(3)} &= \left[0.95 - \frac{0.05}{K-1}\right] I_{K} + \frac{0.05}{K-1} E_{K},\\ P_K^{(4)} &= \frac{0.9}{K-1} E_{K} - \left[\frac{0.9}{K-1} - 0.1\right] I_{K},\label{eqn:sim_hmm_transmat2}
\end{align}
where $E_{K}$ is the $K\times K$ matrix with all elements equal to $1$ and $I_{K}$ is the $K\times K$ identity matrix. The number of observations, $n$, varies from $200$ to $2000$, and the true number of hidden states, $K$, ranges from $3$ to $5$. Figure~\ref{fig:demonsimulatedtraces} illustrates a few simulated HMM traces. 

{\color{blue}
Similar to \cite{Dumont2014}, we compared the proposed method with BIC only. This is because that, among the model selection methods in HMM, (a) AIC is known to be inconsistent, cf. \cite{fuh2024kullback}, and is good for forecasting and not good for estimation. (b) The proposed method has strong connection with BIC as they are both consistent under certain regularity conditions. (c) Other methods such as AICc, MDL, minimum message length (MML), etc. have not been fully explored in HMM model selection literature.}

We conduct $m = 200$ repeated simulations, each of which compares the marginal likelihood method with the BIC as follows. (1) Simulate $n$ observations from the HMM with $K$ states and the specified set of parameters. (2) Apply the Baum-Welch (EM) algorithm with multiple starting points (in our case, $50$ randomly generated starting points) to obtain the (local) maximum likelihood values for $\tilde{K}$-state HMM, thus giving the BIC of HMMs with $\tilde{K}$-states denoted by $BIC_n(\tilde{K})$, $\tilde{K} = 2, 3, 4,\ldots$; let $\hat{K}_{n}^{BIC} = \text{arg max}_{\tilde{K}} BIC_n(\tilde{K})$. (3) Calculate the marginal likelihood of a $\tilde{K}$-state HMM based on the importance sampling procedure detailed in Section~\ref{section:computation}, $\tilde{K} = 2, 3,4,\ldots$; let $\hat{K}_n^{ML} = \text{arg max}_{\tilde{K}} P_{\tilde{K}}({\bm y}_{1:n})$.

\textcolor{blue}{Unless users specify their choices for the hyperparameters in the prior distribution, we set the priors as follows in our R package. And this default choice is the one we adopt for repeated simulation studies. The prior mean for the mean parameters of the $K$ states are set as the $K$ equally spaced quantile levels between $0.05$ and $0.95$ of the observed trajectory, and the corresponding variances are set to be $100^2$. The priors for each row of the transition matrix are flat, that is, Dirichlet with hyperparameters all equal to $1$. For the scaled-inverse chi-square prior for the variance parameters, we set the degree of freedom parameter $\nu=3$ and the scale parameter $s^2$ as:
\begin{equation*}
    s = {\rm quantile}(\boldsymbol{y}, 0.75) - {\rm quantile}(\boldsymbol{y}, 0.25))) / (2 * K),
\end{equation*}
where $\boldsymbol{y}$ is the observed trajectory.
}

Note that although BIC is not well defined for Gaussian HMM as the likelihood is unbounded, in practice, people often use the EM (Baum-Welch) algorithm with multiple starting points to obtain the local maximum of the likelihood and then calculate the BIC based on the (local) maximum from the multiple runs. For a fair comparison, we follow this practice.

Table~\ref{table:simulation_results} summarizes the results from repeated simulations, showing the frequency of correct identification of the true number of hidden states using the marginal likelihood method and the BIC when $n=2000$. Table~\ref{table:simulation_results_continued} in the Appendix gives the results when $n=200$.

\begin{table*}[tbph]
\centering
\caption{Correct identification frequency (percentage) of $K^*$}
\scriptsize
\begin{tabular}{|c|c|c|cc|cc|cc|cc|}
\hline
\hline
\multirow{2}{*}{K} & \multirow{2}{*}{$\sigma$} & \multirow{2}{*}{$n$} & \multicolumn{2}{c|}{$Q_K = P_K^{(1)}$} & \multicolumn{2}{c|}{$Q_K = P_K^{(2)}$} & \multicolumn{2}{c|}{$Q_K = P_K^{(3)}$} & \multicolumn{2}{c|}{$Q_K = P_K^{(4)}$}\\
& & & ML & BIC  & ML & BIC & ML & BIC & ML & BIC\\
\hline
3 & 0.2 & 200 & 99  & 100  & 99   & 100  & 95   & 83   & 100  & 100  \\
3 & 0.3 & 200 & 67  & 10   & 100  & 100  & 96   & 85   & 98.5 & 74.5 \\
3 & 0.4 & 200 & 1.5 & 0.5  & 92   & 65   & 90.5 & 81   & 29   & 2.5  \\
3 & 0.5 & 200 & 0.5 & 0    & 41   & 11.5 & 80   & 63   & 4.5  & 0    \\
\hline
4 & 0.2 & 200 & 98  & 83   & 88.5 & 92   & 85.5 & 56   & 99   & 98   \\
4 & 0.3 & 200 & 6.5 & 0    & 98.5 & 77   & 77   & 51.5 & 30.5 & 0    \\
4 & 0.4 & 200 & 0   & 0    & 50   & 18   & 46.5 & 27   & 0    & 0    \\
4 & 0.5 & 200 & 0   & 0    & 4    & 0    & 11   & 7.5  & 0    & 0    \\
\hline
5 & 0.2 & 200 & 81  & 16.5 & 87   & 66.5 & 67   & 37.5 & 88   & 26.5 \\
5 & 0.3 & 200 & 1   & 0    & 82   & 25   & 42.5 & 19   & 2    & 0    \\
5 & 0.4 & 200 & 0   & 0    & 17   & 2    & 8.5  & 3.5  & 0    & 0    \\
5 & 0.5 & 200 & 0   & 0    & 0.5  & 0    & 0.5  & 0    & 0    & 0 \\ 
\hline
\end{tabular}
\begin{flushleft}
The frequency (in \%) of correct identification of the true number of hidden states, out of $200$ repeated simulations for each entry, using the marginal likelihood method (ML) and the BIC. Here $n=200$ observations are considered (see Online Supplement for $n=2000$ cases). $K(=3,4,5)$ is the true number of hidden states; $\sigma$ is the standard deviation of each hidden state around its mean; $Q_K$ denotes the transition matrix: the matrices $P_K^{(1)}, P_K^{(2)} , P_K^{(3)}, P_K^{(4)}$ are defined in equations (\ref{eqn:sim_hmm_transmat}) to (\ref{eqn:sim_hmm_transmat2}).
\label{table:simulation_results}
\end{flushleft}
\end{table*}

\textcolor{blue}{In the second set of numerical simulations with heterogeneous variances, we set 
\begin{equation*}
    \sigma_{\rm heter} = \sigma_{\rm old} \times (0.5+2 \times runif(K, 0, 1)),
\end{equation*}
where $\sigma_{\rm old}$ refers to the $\sigma$ levels we set for the previous simulation experiments with homogeneous variances. And the $\sigma_{\rm heter}$ are the new heterogeneous variances. The results are in Table~\ref{tab:n2000heter}. The results for $n=200$ are in Table~\ref{tab:n200heter} in the Appendix. }

\begin{table*}[tbph]
\scriptsize
    \centering
    
\textcolor{blue}{
    \caption{Correct identification frequency (percentage) of $K^*$}
\begin{tabular}{|c|c|c|cc|cc|cc|cc|}
\hline
\hline
\multirow{2}{*}{K} & \multirow{2}{*}{$\sigma$} & \multirow{2}{*}{$n$} & \multicolumn{2}{c|}{$Q_K = P_K^{(1)}$} & \multicolumn{2}{c|}{$Q_K = P_K^{(2)}$} & \multicolumn{2}{c|}{$Q_K = P_K^{(3)}$} & \multicolumn{2}{c|}{$Q_K = P_K^{(4)}$}\\
& &  & ML & BIC  & ML & BIC & ML & BIC & ML & BIC\\
            \hline
3&	0.2& 2000&	84.5& 69&		90&	 96&	99&100 &	92.5&		84.5\\
3&	0.3& 2000&	55.5&	29.5&	91.5&76.5&			98.5&94.5&	84 &56.5\\
3&	0.4& 2000&		48&15&		83&52.5&	91&75&	62.5&		33.5\\
3&	0.5& 2000&		42&10&		69.5&34.5&	79.5&	60&	47.5& 16.5\\
4&	0.2& 2000&	56&	26&		81.5&67&		91.5&82.5& 59&	40.5\\
4&	0.3& 2000&		24&6.5&		64&48&		75.5&59.5& 36.5&	16\\
4&	0.4& 2000&	19.5&9&		48&34.5&		60& 56& 23.5&	10\\
4&	0.5& 2000&	10.5&	5& 31&	29&	 54.5& 48.5&	16&	8.5\\
5&	0.2& 2000&		35&12&		56.5& 45.5&		73.5&	63.5& 38.5 & 16\\
5&	0.3& 2000&		20&2.5&	46.5&	38&	 59&	51&	  21& 7\\
5&	0.4& 2000&		10&0.5&	35&	26.5&		45.5&52.5& 6& 	1\\
5&	0.5& 2000&	2.5& 0&		23&12&	34&	40&	 4.5 & 1\\
\hline
    \end{tabular}
    \begin{flushleft}
   The frequency (in \%) of correct identification of the true number of hidden states, out of $200$ repeated simulations for each entry, using the marginal likelihood method (ML) and the BIC. Here $n=2000$ observations are considered (see Online Supplement for $n=200$ cases). $K(=3,4,5)$ is the true number of hidden states; $\sigma$ is the standard deviation of each hidden state around its mean; $Q_K$ denotes the transition matrix: the matrices $P_K^{(1)}, P_K^{(2)} , P_K^{(3)}, P_K^{(4)}$ are defined in equations (\ref{eqn:sim_hmm_transmat}) to (\ref{eqn:sim_hmm_transmat2}). The heterogeneous variances are specified as $ \sigma_{\rm heter} = \sigma \times (0.5+2 \times U$, where $U$ is independently simulated from a uniform distribution on $(0,1)$ for each hidden state.  
    \end{flushleft}}
    \label{tab:n2000heter}
\end{table*}

From the simulation studies, it is evident that the marginal likelihood method outperforms the BIC in several aspects. First, the frequency of correct identification of the number of hidden states using the marginal likelihood method is much higher, especially when the number of observations is small ($200$ as opposed to $2000$). Second, the marginal likelihood method is more robust to low signal-to-noise ratio, which can be seen from Table~\ref{table:simulation_results}. The success rates of the marginal likelihood method and the BIC both drop as the noise level $\sigma$ increases from $0.2$ to $0.5$. However, the success rate of the BIC drops much more as opposed to that of the marginal likelihood. Third, since the number of (unknown) model parameters is quadratic in $K$, given the same number of observations, {\color{blue}the more hidden states there are, the harder the order selection is.} The marginal likelihood method appears more robust to the true number of hidden states than the BIC: the success rate of the marginal likelihood is much higher than that of the BIC when the true number of states is high.

\subsection{Discussion on Choice of Priors}
\label{sec:choice_prior_distribution}

From the asymptotic results in Section~\ref{sec:asymptotic_study}, the influence of priors vanishes as the number of observations goes to infinity. However, in practice, the number of observations is finite, and the choice of priors would have an impact on the results. Now we give our recommendations of the choice of prior distributions based on empirical evidence in running simulation studies. Practitioners should be aware that the \textit{best} prior distribution often comes from incorporating scientific knowledge of the specific problem in the field of study. 

In the simulation studies in Section~\ref{sec:simulation_study}, we choose flat, conjugate priors, and the results look quite promising. The prior for each row of the transition matrix is an independent Dirichlet distribution with parameters all equal to $1$, corresponding to a `flat prior'.  
The priors for the means $\{\mu_k\}_{k=1}^{K}$ are set to be independent Gaussian with means $\{\mu_{0k}\}_{k=1}^{K}$ and large standard deviations, e.g., $10$ or $20$ times the interquartile range of ${\bm y}_{1:n}$. $\{\mu_{0k}\}_{k=1}^K$ is chosen to be data-dependent: the $\mu_{0k}$ are set as the evenly spaced quantiles of the observations ${\bm y}_{1:n}$. The priors for the variances of each hidden state $\{\sigma_k^2\}_{k=1}^K$ are chosen to be independent inverse chi-squared distribution with a degree of freedom $3$, and the scale can be chosen based on empirical estimators of the variability in the data: we can simply take the square root of the scale as the interquartile range of the observations divided by $2K$.

\begin{remark}
Note that the choice of $\nu_K$ is a Dirichlet distribution, which does not have support within $\mathcal{Q}_K^\epsilon$ as in Condition 7) in Section \ref{subsec:consistencytheorem}. However, as discussed in Remark \ref{remark:prior}, it is possible that a $\nu_K$ vanishes quick enough at the boundary of $\mathcal{Q}_K^\epsilon$ is sufficient for the consistency of $\hat{K}_n$. Our simulation studies verify this possibility, although a theoretical proof remains absent for this situation. Note that the Dirichlet prior is also recommended in \cite{Gassiat2014}.
\end{remark}

\section[Applications]{Application to Single-Molecule Experimental Data}
\label{subsec:realdata}

We apply the proposed marginal likelihood method to a set of single-molecule experimental data studied in~\cite{smHHMMfret:2016}, where single-molecule experiments are conducted to study the co-translational protein targeting process, a universal protein transportation mechanism in which proteins are transported to appropriate destinations inside or outside of a cell through the membrane. This process is crucial to the proper functioning of cells, and transportation errors can lead to serious diseases; {\color{blue}see \cite{akopian2013signal}, \cite{saraogi2011molecular} and \cite{shan2016atpase} for more discussion}.

In the biophysics community, each observed time-series trajectory is often modeled as a Gaussian HMM {\color{blue}\citep{SMART, HAMMY:2006}}. The number of hidden states of each HMM corresponds to the number of conformations of a molecular complex, which is of biological significance, as it reveals the dynamics and function of the molecular complex {\color{blue}\citep{Blanco:2010, Yang:2005}.} 

\begin{figure}[tbph]
\centering
\includegraphics[width = 0.6\textwidth]{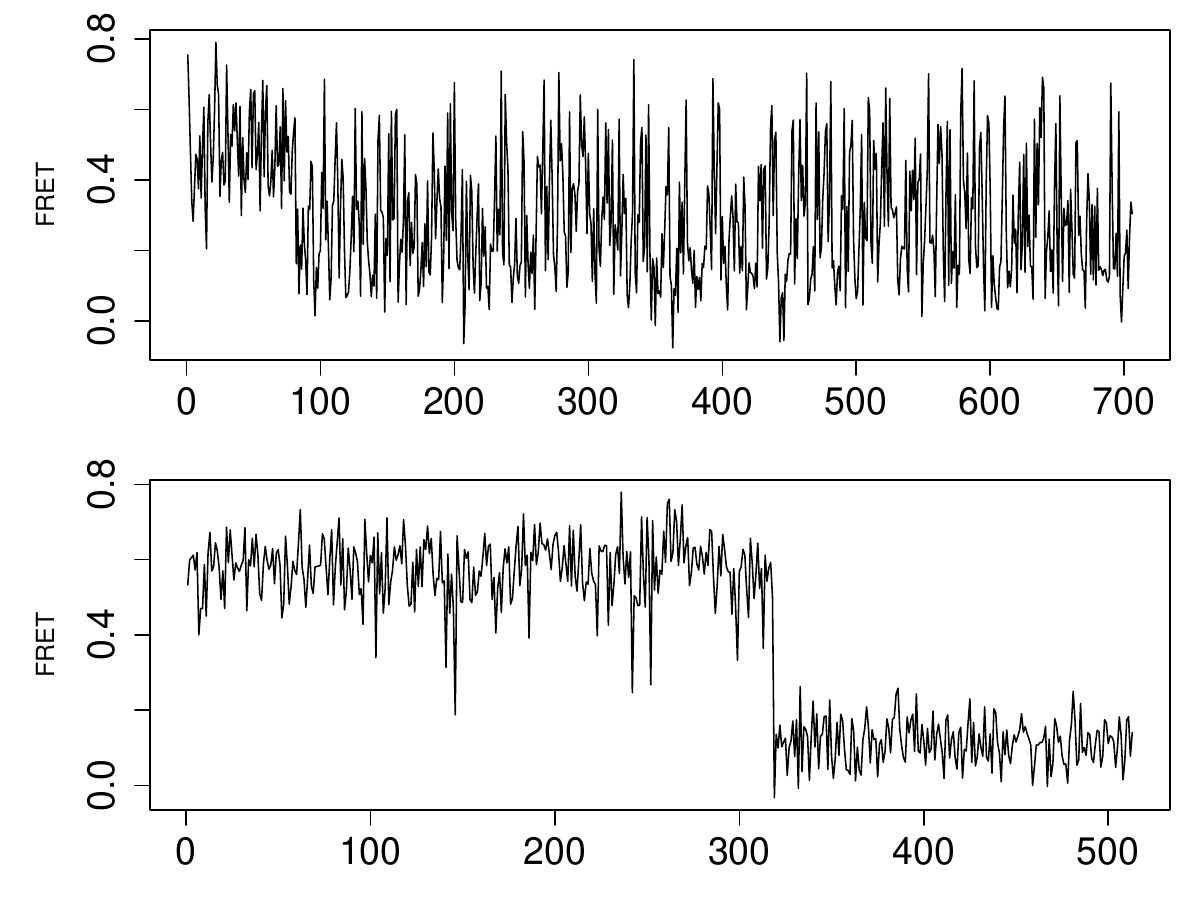}
\caption[Two traces from single-molecule data.]{Two experimental FRET trajectories from the single-molecule data studied in~\cite{smHHMMfret:2016}.}
\label{fig:twotracesTranslocon}
\end{figure}

The marginal likelihood method overall gives similar results as compared to the BIC applied in~\cite{smHHMMfret:2016}. However, for several FRET trajectories, the marginal likelihood method and the BIC give different results. Figure~\ref{fig:twotracesTranslocon} shows two experimental FRET trajectories in which the marginal likelihood method and the BIC give different state-selection: in both cases, the marginal likelihood method gives a selection of three hidden states, whereas the BIC gives a selection of two hidden states. As analyzed and explained in~\cite{smHHMMfret:2016}, these two trajectories are in fact believed to be 3-state trajectories once the information from multiple trajectories was combined together (using a hierarchical HMM) to help identify rarely occurring states. The fact that the marginal likelihood method correctly selects three states in this example indicates that it is more sensitive in detecting rarely occurring states in an HMM. 

\section{Conclusions}
\label{section:conclusions}

In this paper, we use the marginal likelihood to determine the number of hidden states for hidden Markov models. The proposed method is theoretically consistent under mild conditions.
Furthermore, we propose a computation algorithm to robustly estimate the order of an HMM trace through the estimation of normalizing constants. Extensive simulation studies verify our proposed approach and demonstrate its power against the widely adopted approach, the BIC, which lacks theoretical justification. We have provided an R package on CRAN (https://cran.r-project.org/package=HMMmlselect) that implements our proposed method. The package is named HMMmlselect.


\acks{We greatly thank Professor S.C. Kou from Harvard University for his critical comments and assistance on this paper. Y. Chen is supported by  NSF DMS 2113397, NSF PHY 2027555, NASA 22-SWXC22\_2-0005 and 22-SWXC22\_2-0015.}


\begin{appendix}

\section{Proofs of Consistency Theorems}
\label{appendix:proof_theorem_asymptotic_efficiency_hmm}

This section contains five parts. Section \ref{subsec:ProofConcept} briefly explains how we extend the strategy in Section \ref{subsec:illustrate} to general $K^*$ and $K$. Sections \ref{Proof-under} and \ref{Proof-over} prove Theorem \ref{theorem:asymptotic_consistency_hmm}; in particular, Section \ref{Proof-under} proves \eqref{eqn_thm_consist_under}, and Section \ref{Proof-over} proves \eqref{eqn_thm_consist}. The necessary lemmas used in Section \ref{Proof-over}, as well as further theoretical details, are presented in the Online Supplement.

Throughout this section, the notations are consistent with those in the main text. We use $\boldsymbol{\phi}_K$ to denote the parameter under $K$ states, which consists of $Q_K$, the transition matrix, and $\boldsymbol{\theta}_k$, the parameters for $Y_1$ given $X_1=k$, $k=1,\ldots, K$. $\boldsymbol{\theta}_k^*$, $Q^*$ and $\boldsymbol{\phi}^*$ are the true parameters. The priors on $\boldsymbol{\theta}_k$ and $Q_K$ are $\pi$ and $\nu_K$, respectively. For each $K$,  $\Phi_K = \mathcal{Q}_K \times \Theta^K$ is the parameter space, and $\overline{\Phi}_K = \overline{\mathcal{Q}}_K \times \Theta^K$ is the closure of it (for which we use $\overline{S}$ to denote the closure of any set $S$.)

\subsection{Generalizing the Proof Strategy}
\label{subsec:ProofConcept}
In Section \ref{subsec:ProofConcept}, we have demonstrated our proof concept under the special case with $K^*=1$ and $K=2$. In below, we will explain how we extend this strategy to general $K^*$ and $K$, and define necessary notations along the way.

First, for all $K \geq 1$, define $\Delta_K = K(K+d-1)$, which is the number of parameters for a HMM with $K$ states. Then, similar to \eqref{special-case-1}, the Berstein-von Mises theorem for HMM (\cite{HMMCLT}) ensures that the posterior distribution under $K^*$, the true number of states, has asymptotic normality. In other words, there exists a positive definite Fisher information matrix $J^*$ such that
\begin{equation}\label{eqn:CLT-K*}
    \frac{p_{K^*}({\bm y}_{1:n}|\boldsymbol{\phi}^*) p_0(\boldsymbol{\phi}^*)}{n^{\Delta_{K^*}/2} p_{K^*}({\bm y}_{1:n})} \xrightarrow[n \rightarrow \infty]{P^*} (2\pi)^{-\Delta_{K^*}/2} |J^*|^{1/2}.
\end{equation}
This allows us to control the asymptotic behavior for $p_{K^*}({\bm y}_{1:n})$.

Now let us investigate $p_{K}({\bm y}_{1:n})$. In the general case, we have two possibilities:
\begin{itemize}
\item[1)] For $K<K^*$, since a $K$-state HMM can be treated as an under-estimated model for the $K^*$-state HMM, we can prove that there exists $c>0$ such that
\begin{equation}\label{eqn:CLT-under}
    \frac{p_{K}({\bm y}_{1:n})}{p_{K^*}({\bm y}_{1:n}|\boldsymbol{\phi}^*)} = O_{P^*} \left( e^{-cn} \right).
\end{equation}
Combining \eqref{eqn:CLT-K*} and \eqref{eqn:CLT-under}, we have \eqref{eqn_thm_consist_under};
\item[2)] For $K>K^*$, as illustrated in Section \ref{subsec:illustrate}, we will face the non-identifiability issue, so we need to use the strategy in Section \ref{subsec:illustrate} to get
\begin{equation}\label{eqn:CLT-over}
    \frac{n^{\Delta_{K^*}+(K-K^*)d} p_{K}({\bm y}_{1:n})}{p_{K^*}({\bm y}_{1:n}|\boldsymbol{\phi}^*)
    } \xrightarrow[n \rightarrow \infty]{P^*} C
\end{equation}
for some constant $C>0$. Combining \eqref{eqn:CLT-K*} and \eqref{eqn:CLT-over}, we have \eqref{eqn_thm_consist}.
\end{itemize}
Combining the two cases, we prove Theorem \ref{theorem:asymptotic_consistency_hmm}.

Let us be more specific about the strategy in case 2). As illustrated in Section \ref{subsec:illustrate}, the core mechanism in this strategy is to decompose the parameter space into subspaces so that each subspace contains only one ``true'' parameter. Here ``true'' refers to having essentially the same law as we shall see shortly. In the case with $K^*=1$ and $K=2$, as illustrated in \ref{subsec:illustrate}, we use $Q \in \mathcal{Q}_2$ to decompose $\Phi_2$ into $\Phi_Q$, so that each $\Phi_Q$ contains only one ``true'' parameter 
\begin{equation*}
    \boldsymbol{\phi}_Q^* := (Q; \boldsymbol{\theta}_1^*, \boldsymbol{\theta}_1^*).
\end{equation*}
Note that this true parameter can be viewed as we ``split'' the original state into two states, with the corresponding probability ``weight'' determined by $Q$. In other words, in this special case, $Q$ is the index that specifies the ``split'' and ``weight'', while $\boldsymbol{\phi}_Q^*$ is therefore the only ``true'' parameter that satisfies such specification.

This idea of specifying ``split'' and ``weight'' can be extended to general $K > K^*$ as follows. We will reparameterize $\boldsymbol{\phi}_K$ by $(\alpha, \gamma)$ so that $\alpha$ contains all the non-identifiability (including both ``split'' and ``weight''), and $\gamma$ is the identifiable part. As such, through this reparameterization, we can divide the neighborhood of all the  ``true" values into a family of subspaces, with each of them identifiable and having exactly one ``true" value with ergodicity so that the Bernstein-von Mises theorem holds on each of the subspaces. This provides us the subspaces needed to execute steps (i)-(iv) in Section \ref{subsec:illustrate}, which leads us to \eqref{eqn:CLT-over} and completes the proof.

To formally describe this reparameterization approach, we need to first characterize all the ``true'' parameters for a model with $K$ states, following a similar idea in \cite{LerouxM:1992}. To be more precise, for each $K$, recall $\overline{\Phi}_K = \overline{\mathcal{Q}}_K \times \Theta^K$ is the closure of $\Phi_K$. For any $\boldsymbol{\phi}_1 \in \overline{\Phi}_{K_1}$ and $\boldsymbol{\phi}_2 \in \overline{\Phi}_{K_2}$, we define $\boldsymbol{\phi}_1$ and $\boldsymbol{\phi}_2$ to be equivalent, denoting as $\boldsymbol{\phi}_1 \backsim \boldsymbol{\phi}_2$, if and only if there exists initial distribution $\mu(\boldsymbol{\phi}_l) = (\mu_1(\boldsymbol{\phi}_\ell), \cdots, \mu_{K_\ell}(\boldsymbol{\phi}_\ell))$ for $\ell=1, 2$ such that:
\begin{itemize}
    \item[1)] for $\ell=1,2$, if $\{ X_i^\ell: i \geq 0\}$ is a Markov chain with initial distribution $X_0 \sim \mu(\boldsymbol{\phi}_\ell)$ and transition matrix given by $\boldsymbol{\phi}_\ell$, then $\{ \boldsymbol{\theta}_{X_i^\ell}: i \geq 0\}$ is a stationary process;
    \item[2)] the processes  $\{ \boldsymbol{\theta}_{X_i^1}: i \geq 0\}$ and $\{ \boldsymbol{\theta}_{X_i^2}: i \geq 0\}$ follow the same law.
\end{itemize}
Let $\overline{\Phi}_K^* = \{ \boldsymbol{\phi}_K: \boldsymbol{\phi}_K \in \overline{\Phi}_K, \boldsymbol{\phi}_K \backsim \boldsymbol{\phi}^*\}$ be the set of ``true'' parameters under the equivalent class in $\Phi_K$, $\Phi_K^* = \overline{\Phi}_K^* \cap \Phi_K$ be its interior, and $\partial \Phi_K^* = \overline{\Phi}_K^* - \Phi_K^*$ be its boundary.

Now, to obtain the reparametrization $(\alpha, \gamma)$ of $\boldsymbol{\phi}_K$, we first construct a set $\mathcal{A}$, and then construct a set $\Gamma_\alpha$ and a function $\varphi_\alpha$ for each $\alpha \in \mathcal{A}$ so that $\cup_{\alpha \in \mathcal{A}} \varphi_\alpha(\Gamma_\alpha)$ forms a partition of $\Phi_K$. Therefore, we have
\begin{align}\label{change-of-var}
\int_{\Phi_K} p_K({\bm y}_{1:n}|\boldsymbol{\phi}_K) p_0(\boldsymbol{\phi}_K) d\boldsymbol{\phi}_K 
    = \int_{\mathcal{A}} \int_{\Gamma_\alpha} p_K({\bm y}_{1:n}|\varphi_\alpha(\gamma)) p_0(\varphi_\alpha(\gamma)) \vert \Lambda_\alpha(\gamma) \vert d\gamma d\alpha,
\end{align}
in which $\Lambda_\alpha(\gamma)$ is the corresponding Jacobian determinant when changing variables. In addition, we will prove the steps (i)-(iv) in Section \ref{subsec:illustrate} in the following form, respectively:
\begin{enumerate}[(i)]
\item There exists $\mathcal{A}^+ \subset \mathcal{A}$ and $c > 0$ such that, as $n \rightarrow \infty$,
    \begin{align}\label{Apply-Lemma1}
    \frac{ 
    \int_{\mathcal{A}-\mathcal{A}^+} \int_{\Gamma_\alpha} p_K({\bm y}_{1:n}|\boldsymbol{\phi}_K) p_0(\varphi_\alpha(\gamma)) \vert \Lambda_\alpha(\gamma) \vert d\gamma d\alpha}{p_{K^*}({\bm y}_{1:n}|\boldsymbol{\phi}^*)} 
    = O_{P^*}(e^{-c n}).
    \end{align}

\item Define 
\begin{equation}\label{Delta-Gamma}
\Delta_\Gamma = \Delta_{K^*}+(K-K^*)d
\end{equation}
(recall $\Delta_K = K(K+d-1)$). Then, for each $\alpha \in \mathcal{A}^+$, $\Gamma_\alpha$ is $\Delta_\Gamma$-dimensional, and $\overline{\varphi_\alpha(\Gamma_\alpha)} \cap \overline{\Phi}_K^*$
has exactly one component, which is a true value $\boldsymbol{\phi}_\alpha^* \in \Phi_K^*$ with an ergodic transition. In addition, there exists constant $C_\alpha$ such that, as $n \rightarrow \infty$,
    \begin{align}\label{CLT-SubSpace}
    \frac{n^{\frac{\Delta_\Gamma}{2}}  \int_{\Gamma_\alpha} p_K({\bm y}_{1:n}|\varphi_\alpha(\gamma)) p_0(\varphi_\alpha(\gamma)) |\Lambda_\alpha(\gamma)| d\gamma}{p_{K^*}({\bm y}_{1:n}|\boldsymbol{\phi}^*)} 
    = C_\alpha + o_{P^*}(1).
    \end{align}

\item The constant $C_\alpha$ in ii) satisfies 
    \begin{equation}\label{Int-J}
        \int_{\alpha \in \mathcal{A}^
        +} C_\alpha d\alpha < \infty.
    \end{equation}

\item The convergence in \eqref{CLT-SubSpace} is uniform across $\alpha \in \mathcal{A}^+$ so that we have, as $n \rightarrow \infty$, 
    \begin{align}\label{CLT-Int}
    n^{\frac{\Delta_\Gamma}{2}}
    \frac{ 
    \int_{\mathcal{A}^+} \int_{\Gamma_\alpha} p_K({\bm y}_{1:n}|\boldsymbol{\phi}_K) p_0(\varphi_\alpha(\gamma)) \vert \Lambda_\alpha(\gamma) \vert d\gamma d\alpha}{p_{K^*}({\bm y}_{1:n}|\boldsymbol{\phi}^*)}
    = 
    \int_{\alpha \in \mathcal{A}^
        +} C_\alpha d\alpha + o_{P^*}(1).
    \end{align}
\end{enumerate}
Combining \eqref{change-of-var}-\eqref{CLT-Int}, we get \eqref{eqn:CLT-over} and completes the proof.
See Appendix \ref{Proof-over} for the proofs, where the detailed description of the reparameterization $(\alpha, \gamma)$ is given.

\subsection{ Proof of \eqref{eqn_thm_consist_under} ($K<K^*$)}\label{Proof-under}

As discussed in Section \ref{subsec:ProofConcept}, since \eqref{eqn:CLT-K*} holds, we only need to prove \eqref{eqn:CLT-under} holds in order to prove \eqref{eqn_thm_consist_under}. This is provided as follows.

\begin{proof}[Proof of \eqref{eqn:CLT-under} and \eqref{eqn_thm_consist_under}]
Since $K<K^*$, for any $Q_K = \{q_{k\ell}, 1 \leq k, \ell \leq K \} \in \mathcal{Q}_K$, define $\tilde{Q}_{K^*} = \{\tilde{q}_{k\ell}, 1 \leq k, \ell \leq K^* \}$ by
\begin{equation*}
    \tilde{q}_{k\ell} = 
    \begin{cases}
    q_{k\ell} & k,\ell < K \\
    \frac{1}{K^*-K+1}q_{kK} & k < K, \ell \geq K\\ 
    q_{K\ell} & k \geq K, \ell < K\\
    \frac{1}{K^*-K+1}q_{KK} & k, \ell \geq K\\ 
    \end{cases}
\end{equation*}
which is the transition matrix when we ``split" state $K$ into states $K, K+1, \cdots, K^*$ with equal probabilities. In addition, for any $\boldsymbol{\phi}_K = (Q_K; \boldsymbol{\theta}_1, \cdots, \boldsymbol{\theta}_K) \in \Phi_K$, define
\begin{equation*}
    \tilde{\boldsymbol{\phi}}_{K^*} = (\tilde{Q}_{K^*}; \boldsymbol{\theta}_1, \cdots, \boldsymbol{\theta}_{K-1}, \boldsymbol{\theta}_{K}, \boldsymbol{\theta}_{K}, \cdots, \boldsymbol{\theta}_{K}) \in \Phi_{K^*}.
\end{equation*}
Then, a direct computation shows that 
\begin{equation}\label{pf-under-1}
    p_{K}({\bm y}_{1:n}|\boldsymbol{\phi}_K) = p_{K^*}({\bm y}_{1:n}|\tilde{\boldsymbol{\phi}}_{K^*}).
\end{equation}
Moreover, note that
\begin{equation}\label{pf-under-2}
    \tilde{\boldsymbol{\phi}}_K \in \tilde{\Phi}_{K^*} := \left\{ \boldsymbol{\phi}\in \Phi_{K^*}: \boldsymbol{\theta}_K = \boldsymbol{\theta}_{K+1} = \cdots = \boldsymbol{\theta}_{K^*}  \right\},
\end{equation}
and since we assume that the true parameters satisfy $\boldsymbol{\theta}_k^* \neq \boldsymbol{\theta}_\ell^*$ for all $k \neq \ell$, there exists $\delta > 0$ such that
\begin{equation}\label{pf-under-3}
\tilde{\Phi}_{K^*} \subset \{ \boldsymbol{\phi}: \boldsymbol{\phi} \in \Phi_{K^*}, \Vert \boldsymbol{\phi} - \boldsymbol{\phi}^* \Vert \geq \delta \}.
\end{equation}
 Combining \eqref{pf-under-1} to \eqref{pf-under-3}, we have
\begin{align}\label{pf-under-4}
\notag
    p_{K}({\bm y}_{1:n}) = & \int_{\Phi_K} p_{K}({\bm y}_{1:n}|\boldsymbol{\phi}_K) p_0(\boldsymbol{\phi}_K) d\boldsymbol{\phi}_K 
    = \int_{\Phi_K} p_{K^*}({\bm y}_{1:n}|\tilde{\boldsymbol{\phi}}_{K^*}) p_0(\boldsymbol{\phi}_K) d\boldsymbol{\phi}_K \\
    \leq & \sup_{\boldsymbol{\phi} \in \tilde{\Phi}_{K^*}} p_{K^*}({\bm y}_{1:n}|\boldsymbol{\phi}) \int_{\Phi_K}  p_0(\boldsymbol{\phi}_K) d\boldsymbol{\phi}_K 
\leq \sup_{\Vert \boldsymbol{\phi} - \boldsymbol{\phi}^* \Vert \geq \delta } p_{K^*}({\bm y}_{1:n}|\boldsymbol{\phi}).
\end{align}
In addition, Lemma 3.1 in \cite{HMMCLT} states that under conditions 1)-7), for any $\delta > 0$, there exists $\epsilon > 0$ such that
\begin{align}\label{pf-under-5}
P^*\left\{ \sup_{\Vert \boldsymbol{\phi} - \boldsymbol{\phi}^* \Vert \geq \delta} \frac{L_{K^*}({\bm y}_{1:n}|\boldsymbol{\phi}) - L_{K^*}({\bm y}_{1:n}|\boldsymbol{\phi}^*)}{n} \leq -\epsilon \right\} \xrightarrow[n \rightarrow \infty]{} 1,
\end{align}
where $L_K := \log p_K$. Combining \eqref{pf-under-4} and \eqref{pf-under-5}, we have
\begin{align}\label{pf-under-6}
\notag
    & \frac{p_{K}({\bm y}_{1:n})}{p_{K^*}({\bm y}_{1:n}|\boldsymbol{\phi}^*)}
\leq \sup_{\Vert \boldsymbol{\phi} - \boldsymbol{\phi}^* \Vert \geq \delta} \frac{p_{K^*}({\bm y}_{1:n}|\boldsymbol{\phi})}{p_{K^*}({\bm y}_{1:n}|\boldsymbol{\phi}^*)} \\
= & \exp \left\{ n \times \sup_{\Vert \boldsymbol{\phi} - \boldsymbol{\phi}^* \Vert \geq \delta } \frac{L_{K^*}({\bm y}_{1:n}|\boldsymbol{\phi}) - L_{K^*}({\bm y}_{1:n}|\boldsymbol{\phi}^*)}{n} \right\} = O_{P^*}(e^{-\epsilon n})
\end{align}
as $n \rightarrow \infty$. Combining \eqref{eqn:CLT-K*} and \eqref{pf-under-6}, we complete the proof.
\end{proof}

\begin{remark}
One of the main ingredients in the proof above is \eqref{pf-under-5}, which essentially says that, outside of any open ball of $\boldsymbol{\phi}^*$, the likelihood $p_{K^*}({\bm y}_{1:n}|\boldsymbol{\phi})$ decays exponentially compare to $p_{K^*}({\bm y}_{1:n}|\boldsymbol{\phi}^*)$ as $n \rightarrow \infty$. This can actually be extended to $p_K$ with $K > K^*$. We state this as Lemma \ref{lemma:concentrate} in Section \ref{proof:glossary}, which will be used as an intermediate step in the proof of Lemma \ref{Lemma:i}.
\end{remark}

\subsection{ Proof of \eqref{eqn_thm_consist} ($K>K^*$)}\label{Proof-over}

 As discussed in Section \ref{subsec:ProofConcept}, since \eqref{eqn:CLT-K*} holds, we need to prove \eqref{change-of-var} - \eqref{CLT-Int} for the proof of \eqref{eqn_thm_consist}. To this end, we will have a reparameterization $(\alpha, \gamma)$ of $\boldsymbol{\phi}_K$ such that $\alpha$ contains the non-identifiability parameters, and $\gamma$ is the identifiable part. We will first present such reparameterization and then prove \eqref{change-of-var} - \eqref{CLT-Int} accordingly.

Note that, as discussed in Section \ref{subsec:ProofConcept}, for the construction of $(\alpha, \gamma)$ to be a reparameterization with $\alpha$ containing all the non-identifiability, $\alpha$ needs to capture the ``split'' of the states as well as the ``weight'' of the split. To illustrate the idea behind the reparameterization, we first discuss two special cases and then present the construction of the general case.

\begin{example}[Mixture Normal with Two Components]\label{exmp-mixture}
Consider the case when $q_{k\ell} = q_\ell$ for all $\ell=1, 2, \cdots, K$ (so that $X_i$ are i.i.d.), and $Y_i \backsim N({\boldsymbol{\theta}}_{X_i}, 1)$. When the true number of components $K^*=2$, and the true parameter $\boldsymbol{\phi}^*:=((q_1^*, q_2^*), ({\boldsymbol{\theta}}_1^*, {\boldsymbol{\theta}}_2^*))$ is defined as
\begin{equation*}
    (q_1^*, q_2^*) = (\frac{1}{2}, \frac{1}{2}),~~~ ({\boldsymbol{\theta}}_1^*, {\boldsymbol{\theta}}_2^*) = (1,2).
\end{equation*}
Then the likelihood function under $\boldsymbol{\phi}^*$ is
\begin{equation*}
    p_2({\bm y}_{1:n}|\boldsymbol{\phi}^*) = \prod_{i=1}^n \left\{ \frac{1}{2} \phi(y_i - 1) + \frac{1}{2} \phi(y_i-2) \right\},
\end{equation*}
where $\phi(\cdot)$ is the probability density function for the standard normal distribution.

Now, suppose we fit the model with $K = 3$ and parameter $\boldsymbol{\phi}_3 \in \Phi_3$ defined by $(q_1, q_2, q_3)$ and $(\boldsymbol{\theta}_1, \boldsymbol{\theta}_2, \boldsymbol{\theta}_3)$, which corresponds to the likelihood function
\begin{align*}
    p_3({\bm y}_{1:n}|\boldsymbol{\phi}_3) 
    = \prod_{i=1}^n \left\{ q_1 \phi(y_i - {\boldsymbol{\theta}}_1) + q_2 \phi(y_i-{\boldsymbol{\theta}}_2) + q_3 \phi(y_i - {\boldsymbol{\theta}}_3) \right\}.
\end{align*}
This creates identifiability issue since there are multiple $\boldsymbol{\phi}_3$ satisfying $\boldsymbol{\phi}_3 \backsim \boldsymbol{\phi}^*$. For example, the following $\boldsymbol{\phi}_3$ all makes $p_3(\cdot|\boldsymbol{\phi}_3) \equiv p_2(\cdot|\boldsymbol{\phi}^*)$:
\begin{itemize}
    \item $(q_1, q_2, q_3) = (\frac{1}{4}, \frac{1}{4}, \frac{1}{2})$, $({\boldsymbol{\theta}}_1, {\boldsymbol{\theta}}_2, {\boldsymbol{\theta}}_3) = (1, 1, 2)$;
    \item $(q_1, q_2, q_3) = (\frac{1}{3}, \frac{1}{6}, \frac{1}{2})$, $({\boldsymbol{\theta}}_1, {\boldsymbol{\theta}}_2, {\boldsymbol{\theta}}_3) = (1, 1, 2)$;
    \item $(q_1, q_2, q_3) = (\frac{1}{2}, \frac{1}{4}, \frac{1}{4})$, $({\boldsymbol{\theta}}_1, {\boldsymbol{\theta}}_2, {\boldsymbol{\theta}}_3) = (1, 2, 2)$.
\end{itemize}
In fact, any $\boldsymbol{\phi}_3$ satisfying
\begin{equation}\label{exmp-mixture-iden}
    {\boldsymbol{\theta}}_k = \boldsymbol{\theta}_\ell^* = \ell \mbox{ for all } k \in S_\ell,  \mbox{ and } \sum_{k \in S_\ell} q_k = q_\ell^* = \frac{1}{2} 
\end{equation}
for all $\ell = 1, 2,$ for some partition $S=(S_1, S_2)$ of $\{1, 2, 3\}$
will make $\boldsymbol{\phi}_3 \backsim \boldsymbol{\phi}^*$. 

Now, suppose $(\alpha, \gamma)$ is a reparameterization of $\boldsymbol{\phi}_3$, and $\Phi_{3,\alpha}$ is the subspace of $\Phi_3$ conditioning on a given $\alpha$. Then, by \eqref{exmp-mixture-iden}, to make sure that $\alpha$ contains all the non-identifiability, we need that for given $\alpha$, \eqref{exmp-mixture-iden} has at most one solution in $\Phi_{3,\alpha}$. This means that we need to control $\alpha$ by at least the following two matters:
\begin{enumerate}
\item To make sure the $\boldsymbol{\theta}_k$ part in \eqref{exmp-mixture-iden} has at most one solution, $\alpha$ must specify the partition $S=(S_1, S_2)$. Equivalently, this means that $\alpha$ needs to specify, for each true state $\ell$, the collection of states $S_\ell$ that the state $\ell$ is ``split'' into.

\item Given $S$, note that the $q_k$ part in \eqref{exmp-mixture-iden} only specify the sum of $q_k$ across each $S_\ell$. Therefore, to ensure that the $q_k$ part in \eqref{exmp-mixture-iden} has at most one solution, $\alpha$ must also specify
\begin{align*}
    \left( \frac{q_1}{\sum_{k \in S_{s(1)}} q_k}, \frac{q_2}{\sum_{k \in S_{s(2)}} q_k}, \frac{q_3}{\sum_{k \in S_{s(3)}} q_k} \right),
\end{align*}
where
\begin{equation}\label{exmp-s(k)}
    s(k) = 
    \begin{cases}
    1, & if~k \in S_1, \\
    2, & if~k \in S_2.
    \end{cases}
\end{equation}
Equivalently, this means that $\alpha$ also needs to specify the ``weight'' of each $q_i$ relative to $\sum_{k \in S_{s(i)}} q_k$.
\end{enumerate}

In fact, we can construct such $\alpha$ and $\Phi_{3,\alpha}$ as follows. Set $\alpha = (S,W)$, in which $S=(S_1, S_2)$ is a partition of $\{1, 2, 3\}$, and
\begin{equation*}
    W = (W_1, W_2, W_3) \in (0,1)^3.
\end{equation*}
To construct $\Phi_{3,\alpha}$, we first set
\begin{equation}\label{Theta_k_exmp}
    \Theta_1 = (-\infty, 1.5), ~~ \Theta_2 = [1.5, \infty),
\end{equation}
so that
\begin{equation}\label{exmp-mixture-unique}
    \boldsymbol{\theta} \in \{\boldsymbol{\theta}_1^*, \boldsymbol{\theta}_2^*\} \cap \Theta_\ell =  \{1, 2\} \cap \Theta_\ell  \Rightarrow \boldsymbol{\theta} = \ell.
\end{equation}
Then, for any $\alpha = (S,W)$, we set
\begin{align*}
    \Phi_{3,\alpha} = \Big\{ \boldsymbol{\phi}_3: & \mbox{ for all } \ell \in \{1, 2\},  \boldsymbol{\theta}_k \in \Theta_\ell \mbox{ for all } k \in S_\ell \\
    & \mbox{ and for each } i \in \{1, 2, 3\},  \frac{q_i}{\sum_{k \in S_{s(i)}}q_k} = W_i \Big\}.
\end{align*}

To show that \eqref{exmp-mixture-iden} has at most one solution in $\Phi_{3,\alpha}$ for any $\alpha$, suppose that $\boldsymbol{\phi}_3 \in \Phi_{3,\alpha}$ satisfies \eqref{exmp-mixture-iden} for some partition $S'$. Since $\boldsymbol{\phi}_3$ satisfies the first part of \eqref{exmp-mixture-iden}, by \eqref{exmp-mixture-unique}, we must have $S'=S$, so that we have $\boldsymbol{\theta}_k = \boldsymbol{\theta}_\ell^* = \ell$ for all $k \in S_\ell$. In addition, for any $i \in \{1,2,3\}$,
\begin{equation*}
    q_i = \frac{q_i}{\sum_{k \in S_{s(i)}} q_k} \times \sum_{k \in S_{s(i)}} q_k 
    = W_i \times \sum_{k \in S_{s(i)}} q_k
    = W_i \times \frac{1}{2},
\end{equation*}
in which the second equality is because $\boldsymbol{\phi}_3 \in \Phi_{3,\alpha}$, and the last equality is because $\boldsymbol{\phi}_3$ satisfies \eqref{exmp-mixture-iden}. To sum up, given $\alpha$, the solution of \eqref{exmp-mixture-iden} is uniquely determined by $\alpha$; in other words, there is at most one solution in $\Phi_{3,\alpha}$ for each $\alpha$, as desired.
\end{example}

\begin{remark}
The identifiability issue in mixture models has also been discussed in the BIC literature; see \cite{keribin2000consistent} as well as \cite{drton2017bayesian}. See also \cite{watanabe2013widely} for some extensions toward machine learning.
\end{remark}

\begin{example}[Gaussian HMM with Two Components]\label{exmp-GaussianHMM}
Now let us extend the idea in Example \ref{exmp-mixture} to Gaussian HMM. Let $Y_i \backsim N(\boldsymbol{\theta}_{X_i}, 1)$, but without the constraint of $q_{k\ell}=q_\ell$. Consider the case when the true number of states is $K^*=2$, and the true parameter $\boldsymbol{\phi}^*=(Q^*, (\boldsymbol{\theta}_1^*, \boldsymbol{\theta}_2^*))$ is defined by
\begin{equation*}
    Q^* = 
    \begin{pmatrix}
    q_{11}^* & q_{12}^* \\
    q_{21}^* & q_{22}^*
    \end{pmatrix}
    =
    \begin{pmatrix}
        \frac{1}{2} & \frac{1}{2} \\
        \frac{1}{2} & \frac{1}{2}
    \end{pmatrix},~~~
    (\boldsymbol{\theta}_1^*, \boldsymbol{\theta}_2^*) = (1,2).
\end{equation*}
Note that the invariant probability measure corresponding to $Q^*$ is $\mu(\boldsymbol{\phi}^*) = (\mu_1(\boldsymbol{\phi}^*), \mu_2(\boldsymbol{\phi}^*)) = (\frac{1}{2}, \frac{1}{2})$.

Again, we fit the model with $K=3$. Consider $\boldsymbol{\phi}_3 \in \Phi_3$ and its corresponding invariant probability measure $\mu(\boldsymbol{\phi}_3) = (\mu_1(\boldsymbol{\phi}_3), \mu_2(\boldsymbol{\phi}_3),  \mu_3(\boldsymbol{\phi}_3))$. Similar to equation \eqref{exmp-mixture-iden} in Example \ref{exmp-mixture}, we like to find the condition that makes $\boldsymbol{\phi}_3 \backsim \boldsymbol{\phi}^*$. Note that by the definition of the equivalent class, $\boldsymbol{\phi}_3 \backsim \boldsymbol{\phi}^*$ means that, under the corresponding invariant probability, the distributions of $\{ \boldsymbol{\theta}_{X_i}, i \geq 0\}$ are the same; in particular, $(\boldsymbol{\theta}_{X_0}, \boldsymbol{\theta}_{X_1})$ have the same distribution. Let $P_{\mu(\boldsymbol{\phi})}^{\boldsymbol{\phi}}$ denote the probability under ${\boldsymbol{\phi}}$ and invariant probability $\mu(\boldsymbol{\phi})$. Then, for any $i, j = 1, 2$, the probability
\begin{align*}
   P_{\mu(\boldsymbol{\phi}_3)}^{\boldsymbol{\phi}_3} \left\{ (\boldsymbol{\theta}_{X_0}, \boldsymbol{\theta}_{X_1}) = (i,j) \right\} 
   & =
   P_{\mu(\boldsymbol{\phi}_3)}^{\boldsymbol{\phi}_3}
   \left\{ X_0 \in \{ k: \boldsymbol{\theta}_k = i\}, X_1 \in \{ \ell: \boldsymbol{\theta}_\ell = j\} \right\} \\
   & = 
   \sum_{k: \boldsymbol{\theta_k} =i} 
   \sum_{\ell: \boldsymbol{\theta_\ell} =j} 
   \mu_k(\boldsymbol{\phi}_3) q_{k\ell}
\end{align*}
must be equal to
\begin{align*}
   P_{\mu(\boldsymbol{\phi}^*)}^{\boldsymbol{\phi}^*} \left\{ (\boldsymbol{\theta}_{X_0}, \boldsymbol{\theta}_{X_1}) = (i,j) \right\} 
   =
   P_{\mu(\boldsymbol{\phi}^*)}^{\boldsymbol{\phi}^*}
   \left\{ X_0 = i, X_1 = j \right\} 
   = 
   \mu_i(\boldsymbol{\phi}^*) q_{ij}^* = \frac{1}{4}.
\end{align*}
Hence, similar to \eqref{exmp-mixture-iden}, for $\boldsymbol{\phi}_3 \backsim \boldsymbol{\phi}^*$, we have
\begin{align}\label{exmp-GaussianHMM-iden}
\notag
    & {\boldsymbol{\theta}}_k = \boldsymbol{\theta}_i^* = i \mbox{ for all } k \in S_i, \\
    \mbox{ and } & \sum_{k \in S_i} \sum_{\ell \in S_j} \mu_k(\boldsymbol{\phi}_3) q_{k\ell} = \mu_i(\boldsymbol{\phi}^*) q_{ij}^* = \frac{1}{4} 
\end{align}
for all $i,j = 1, 2,$ for some partition $S=(S_1, S_2)$ of $\{1, 2, 3\}$. Note that \eqref{exmp-GaussianHMM-iden} is almost like \eqref{exmp-mixture-iden}, just changing the consideration of $\boldsymbol{\theta}_{X_1}$ to the pair $(\boldsymbol{\theta}_{X_0}, \boldsymbol{\theta}_{X_1})$.

As that in Example \ref{exmp-mixture}, \eqref{exmp-GaussianHMM-iden} suggests that we should have $\alpha = (S,W)$ with $S$ being the partition of $\{1,2,3\}$, and 
\begin{equation}\label{exmp-W}
    W = 
    \begin{pmatrix}
    W_{11} & W_{12} & W_{13} \\
    W_{21} & W_{22} & W_{23} \\
    W_{31} & W_{32} & W_{33} 
    \end{pmatrix}
\end{equation}
such that $W_{ij} \in (0,1)$ for all $1 \leq i, j \leq 3$.

The construction of $\Phi_{3,\alpha}$ for given $\alpha = (S,W)$ is also similar. Recall $s(k)$ defined in \eqref{exmp-s(k)}, and $\Theta_1$ and $\Theta_2$ defined the \eqref{Theta_k_exmp}. For any $\alpha = (S,W)$, set
\begin{align*}
    \Phi_{3,\alpha} = & \Bigg\{ \boldsymbol{\phi}_3: \mbox{ for all } \ell \in \{1, 2\},  \boldsymbol{\theta}_k \in \Theta_\ell \mbox{ for all } k \in S_\ell  \\
    &  \mbox{ and for all } i,j \in \{1, 2, 3\}, \frac{\mu_iq_{ij}}{\sum_{k \in S_{s(i)}}\sum_{\ell \in S_{s(j)}} \mu_k q_{k\ell}} = W_{ij} \Bigg\}.
\end{align*}
By using an argument similar to that in Example \ref{exmp-mixture}, we have that for each $\alpha = (S, W)$, \eqref{exmp-GaussianHMM-iden} has at most one solution in $\Phi_{3,\alpha}$.
\end{example}~\\

We now extend the idea in Example \ref{exmp-GaussianHMM} to a general setting. First, for the construction of $\alpha = (S,W)$, let  
\begin{equation*}
S = (S_1, S_2, \cdots, S_{K^*})
\end{equation*}
be a partition of $\{1, 2, \cdots, K\}$. Given $S$, we extend the definition of $s(k)$ in \eqref{exmp-s(k)} to
\begin{equation}\label{s(k)}
    s(k) = \ell \mbox{ if } k \in S_\ell.
\end{equation}
As for $W$ in \eqref{exmp-W}, we generalize it to
\begin{equation*}
    W = 
    \begin{pmatrix}
    W_{11} & W_{12} & \cdots & W_{1K} \\
    W_{21} & W_{22} & \cdots & W_{2K} \\
    \vdots & \vdots & \ddots & \vdots \\
    W_{K1} & W_{K2} & \cdots & W_{KK} \\
    \end{pmatrix}
\end{equation*}
such that $W_{ij} \in (0,1)$ for all $1 \leq i,j \leq K$.

As for the construction of $\Phi_{K,\alpha}$, similar to \eqref{Theta_k_exmp}, we first decompose $\Theta= \cup_{1 \leq k \leq K^*} \Theta_k$ such that $\Theta_k$ are disjoint, and $\boldsymbol{\theta}_k^*$ is an interior point of $\Theta_k$ for all $k=1, 2, \cdots, K^*$. Similar to \eqref{exmp-mixture-unique},  this guarantees that
\begin{equation}\label{Theta-unique}
\boldsymbol{\theta} \in \left\{\boldsymbol{\theta}_1^*, \cdots, \boldsymbol{\theta}_K^*\right\} \cap \Theta_k \Rightarrow \boldsymbol{\theta} = \boldsymbol{\theta}_k^*.
\end{equation}
Now, given $\alpha = (S, W)$, set
\begin{align}\label{Phi_Kalpha}
    \Phi_{K,\alpha} = & \Bigg\{ \boldsymbol{\phi}_K: \mbox{ for all } \ell \in \{1, 2, \cdots, K^*\},
    \boldsymbol{\theta}_k \in \Theta_\ell \mbox{ for all } k \in S_\ell,  \\
\notag
    & \mbox{ and for all } i,j \in \{1, 2, \cdots, K\},  \frac{\mu_iq_{ij}}{\sum_{k \in S_{s(i)}}\sum_{\ell \in S_{s(j)}} \mu_k q_{k\ell}} = W_{ij} \Bigg\}.
\end{align}

Next using $\alpha$ to construct the reparameterization $(\alpha, \gamma)$ of $\Phi_K$, we first note that by the definition of $\Phi_{K,\alpha}$ in \eqref{Phi_Kalpha}, for each $\boldsymbol{\phi}_K \in \Phi_K$, there exists an unique $\alpha$ such that $\boldsymbol{\phi}_K \in \Phi_{K,\alpha}$. Define
\begin{equation}\label{set-A}
    \mathcal{A} = \left\{\alpha = (S,W): \Phi_K \cap \Phi_{K,\alpha} \neq \emptyset \right\}
\end{equation}
so that $\cup_{\alpha \in \mathcal{A}} \Phi_{K,\alpha}$ forms a partition of $\Phi_K$.

It remains to parameterize $\Phi_{K,\alpha}$ for any $\alpha \in \mathcal{A}$ through $\gamma$ so that we can rewrite the marginal likelihood as an integration over $(\alpha, \gamma)$ in  \eqref{change-of-var}. As $\boldsymbol{\phi}_K = (Q_K; \boldsymbol{\theta}_1, \cdots, \boldsymbol{\theta}_K)$, this means that we need a reparameterization of $Q_K$. In addition, note that for given $\alpha = (S,W)$, we have  
\begin{equation*}
    \frac{\mu_iq_{ij}}{\sum_{k \in S_{s(i)}}\sum_{\ell \in S_{s(j)}} \mu_k q_{k\ell}} = W_{ij} 
    ~~ \Longleftrightarrow ~~
    \mu_iq_{ij} = W_{ij} \times \left\{\sum_{k \in S_{s(i)}}\sum_{\ell \in S_{s(j)}} \mu_k q_{k\ell} \right\}
\end{equation*}
for all $i,j \in \{1, 2, \cdots, K\}$, so $\left\{\sum_{k \in S_{s(i)}}\sum_{\ell \in S_{s(j)}} \mu_k q_{k\ell} \right\}$ basically suggests a reparameterization of $\mu_i q_{ij}$, which can further lead to a reparameterization of $Q_K$.

The details are given below. Let us first set
\begin{equation*}
    R =
    \begin{pmatrix}
    R_{11} & R_{12} & \cdots & R_{1K^*} \\
    R_{21} & R_{22} & \cdots & R_{2K^*} \\
    \vdots & \vdots & \ddots & \vdots \\
    R_{K^*1} & R_{K^*2} & \cdots & R_{K^*K^*}
    \end{pmatrix}
\end{equation*}
such that $R_{ij} \in (0,1)$ for all $1 \leq i,j \leq K^*$. For fixed $\alpha = (S,W)$, consider the mapping $\psi_{\alpha}$ from $R$ to $Q_K$ as
\begin{equation}\label{psi}
    \psi_{\alpha}(R) =
    \begin{pmatrix}
    \psi_{\alpha 11}(R) & \psi_{\alpha 12}(R) & \cdots & \psi_{\alpha 1K}(R)  \\
    \psi_{\alpha 21}(R) & \psi_{\alpha 22}(R) & \cdots & \psi_{\alpha 2K}(R)  \\
    \vdots & \vdots & \ddots & \vdots \\
    \psi_{\alpha K1}(R) & \psi_{\alpha K2}(R) & \cdots & \psi_{\alpha KK}(R)
    \end{pmatrix}.
\end{equation}
Here for all $1 \leq i, j \leq K$,
\begin{equation}\label{psi-ij}
 \psi_{\alpha ij}(R) 
 =  \frac{W_{ij} R_{s(i)s(j)}}{\sum_{\ell=1}^K W_{i\ell} R_{s(i)s(\ell)}}.
\end{equation}
Further set
\begin{equation*}\label{gamma}
\gamma = (R; \boldsymbol{\theta}_1, \cdots, \boldsymbol{\theta}_K),
\end{equation*}
and define the mapping $\varphi_\alpha$ from $\gamma$ to $\boldsymbol{\phi}_K = (Q_K; \boldsymbol{\theta}_1, \cdots, \boldsymbol{\theta}_K)$ as
\begin{equation}\label{varphi}
\varphi_\alpha(\gamma) = (\psi_\alpha(R); \boldsymbol{\theta}_1, \cdots, \boldsymbol{\theta}_K).
\end{equation}

To show that for any $\alpha \in \mathcal{A}$, $\gamma$ indeed forms a reparametrization of $\Phi_{K,\alpha}$ through $\varphi_\alpha$, we need to verify that for any $\alpha \in \mathcal{A}$, $\varphi_\alpha$ is a bijective function. To do so, we directly claim that
\begin{equation}\label{varphi-inv}
    \varphi_\alpha^{-1}(\boldsymbol{\phi}_K) =
    \left( \psi_\alpha^{-1}(Q_K); \boldsymbol{\theta}_1, \cdots, \boldsymbol{\theta}_K \right),
\end{equation}
where 
\begin{equation*}
    \psi_\alpha^{-1}(Q_K) = 
    \begin{pmatrix}
    \psi_{\alpha 11}^{-1}(Q_K) & \psi_{\alpha 12}^{-1}(Q_K) & \cdots & \psi_{\alpha 1K^*}^{-1}(Q_K) \\
    \psi_{\alpha 21}^{-1}(Q_K) & \psi_{\alpha 22}^{-1}(Q_K) & \cdots & \psi_{\alpha 2K^*}^{-1}(Q_K) \\
    \vdots & \vdots & \ddots & \vdots \\
    \psi_{\alpha K^*1}^{-1}(Q_K) & \psi_{\alpha K^*2}^{-1}(Q_K) & \cdots & \psi_{\alpha K^*K^*}^{-1}(Q_K)
    \end{pmatrix}.
\end{equation*}
Here for each $1 \leq i,j \leq K^*$, 
\begin{equation}\label{psi-inv}
    \psi_{\alpha i j}^{-1}(Q_K) = \sum_{k \in S_i} \sum_{\ell \in S_j} \mu_k q_{k \ell},
\end{equation}
which is the denominator appears in \eqref{Phi_Kalpha}. To show that $\varphi_\alpha^{-1}$ is indeed the inverse function of $\varphi_\alpha$, suppose that $\gamma = (R; \boldsymbol{\theta}_1, \cdots, \boldsymbol{\theta}_K) =  \varphi_\alpha^{-1}(\boldsymbol{\phi}_K)$ for some $\boldsymbol{\phi}_K \in \Phi_{K,\alpha}$.
Since $\boldsymbol{\phi}_K \in \Phi_{K,\alpha}$, by \eqref{Phi_Kalpha} and \eqref{psi-inv}, we have
\begin{equation}\label{mu-q-rep}
    \mu_i q_{ij} = W_{ij} \sum_{k \in S_{s(i)}}\sum_{\ell \in S_{s(j)}} \mu_k q_{k \ell} = W_{ij} \psi_{\alpha s(i) s(j)}^{-1}(Q_K),
\end{equation}
which leads to 
\begin{equation}\label{mu-rep}
    \mu_i = \sum_{j=1}^K W_{i j} \psi_{\alpha s(i) s(j)}^{-1}(Q_K)
\end{equation}
since $\sum_{j=1}^K q_{ij} = 1$. Hence, by \eqref{mu-q-rep} and \eqref{mu-rep}, we have
\begin{equation}\label{A.26}
    q_{ij} = \frac{\mu_i q_{ij}}{\mu_i} = 
    \frac{W_{ij} \psi_{\alpha s(i) s(j)}^{-1}(Q_K)}{
    \sum_{\ell=1}^K W_{i \ell} \psi_{\alpha s(i) s(\ell)}^{-1}(Q_K)}.
\end{equation}
Combining \eqref{A.26} with \eqref{psi-ij}, we have that
\begin{equation*}
    \psi_\alpha(\psi_\alpha^{-1}(R)) = R,
\end{equation*}
and thus
\begin{equation*}
    \varphi_\alpha(\varphi_\alpha^{-1}(\gamma)) = \gamma.
\end{equation*}
As such, $\varphi_\alpha$ is invertible, therefore if we set
\begin{equation*}
    \Gamma_\alpha := \varphi_\alpha^{-1}(\Phi_{K,\alpha}),
\end{equation*}
then $\varphi_\alpha: \Gamma_\alpha \rightarrow \Phi_{K,\alpha}$ is a bijective function. This proves that $(\alpha, \gamma)$ indeed forms a reparameterization of $\Phi_K$.

We then show that $\alpha$ indeed captures all the non-identifiability; that is, we want to show that for each $\alpha = (S,W)\in \mathcal{A}$, $\overline{\Phi}_{K,\alpha}$ contains at most one ``true'' value. To show this, consider any true value $\boldsymbol{\phi}_K^* = (Q_K; \boldsymbol{\theta}_1, \cdots, \boldsymbol{\theta}_K) \in \overline{\Phi}_K^*$, and suppose it lies in $\overline{\Phi}_{K,\alpha}$ for some $\alpha \in \mathcal{A}$. Given this $\alpha$, set $\gamma_\alpha^* = \varphi_\alpha^{-1}(\boldsymbol{\phi}_K^*)$. In addition, since $\boldsymbol{\phi}_K^* \backsim \boldsymbol{\phi}^*$, $(\boldsymbol{\theta}_{X_0}, \boldsymbol{\theta}_{X_1})$ have the same law under $\boldsymbol{\phi}_K^*$ and $\boldsymbol{\phi}^*$, by condition 3) we have
\begin{align}\label{id-1}
     \sum_{i,j=1}^{K} \sum_{k \in S_{s(i)}} \sum_{\ell \in S_{s(j)}} \mu_k(\boldsymbol{\phi}_K^*) q_{k\ell} 1_{\boldsymbol{\theta}_k}(\boldsymbol{\theta}_0)
     1_{\boldsymbol{\theta}_\ell}(\boldsymbol{\theta}_1)
\equiv 
     \sum_{i,j=1}^{K^*} \mu_i^*(\boldsymbol{\phi}^*) q_{ij}^* 
     1_{\boldsymbol{\theta_i^*}}(\boldsymbol{\theta}_0)
     1_{\boldsymbol{\theta_j^*}}(\boldsymbol{\theta}_1)
\end{align}
as a function of $(\boldsymbol{\theta}_0, \boldsymbol{\theta}_1)$.

Note that by \eqref{Phi_Kalpha} and \eqref{set-A}, for any $\alpha = (S,W) \in \mathcal{A}$, $W$ has positive entries, so by \eqref{psi-ij}, for any $\gamma = (R; \boldsymbol{\theta}_1, \cdots, \boldsymbol{\theta}_K) \in \overline{\Gamma}_\alpha$, $\psi_\alpha(R)$ is a transition matrix with positive entries. As such, since $\gamma_\alpha^* = \varphi_\alpha^{-1}(\boldsymbol{\phi}_K^*) \in \varphi_\alpha^{-1}(\overline{\Phi}_{K,\alpha}) =  \overline{\Gamma}_\alpha$, we have $\mu_k(\boldsymbol{\phi}_K^*) = \mu_k(\varphi_\alpha(\gamma_\alpha^*)) > 0$ for all $\alpha \in \mathcal{A}$ and $k=1, 2, \cdots, K$. In addition, since $\boldsymbol{\phi}_K^* \in \overline{\Phi}_{K,\alpha}$, by \eqref{Phi_Kalpha} and \eqref{Theta-unique}, we have $\boldsymbol{\theta}_k \in \overline{\Theta}_i$ for all $i = 1, \cdots, K^*$ and $k \in S_{s(i)}$, which means that $\boldsymbol{\theta}_k \in \left\{ \boldsymbol{\theta}_1^*, \cdots, \boldsymbol{\theta}_{K^*}^* \right\}$ if and only if $\boldsymbol{\theta}_k = \boldsymbol{\theta}_i^*$. Hence, the only possibility for \eqref{id-1} to hold is when
\begin{align}\label{theta-equal}
\boldsymbol{\theta}_k = \boldsymbol{\theta}_i^* & \mbox{ for all } i = 1, 2, \cdots, K^* \mbox{ and }  k \in S_{s(i)};
\end{align}
and when this holds, \eqref{psi-inv} and \eqref{id-1} directly lead to
\begin{align}\label{explain-construct}
    R_{ij} &  = \psi_{\alpha i j}^{-1}(Q_K) = \sum_{k \in S_i} \sum_{\ell \in S_j} \mu_{k} q_{k \ell}
     = \mu_i^*(\boldsymbol{\phi}^*) q_{ij}^* := R_{ij}^*.
\end{align}
Therefore, $\boldsymbol{\phi}_K = (Q_K; \boldsymbol{\theta}_1, \cdots, \boldsymbol{\theta}_K)$ is a ``true" value in $\overline{\Phi}_{K,\alpha}$ if and only if 
\begin{align}\label{R*}
\psi_\alpha^{-1}(Q_K) =R^* 
:= 
\begin{pmatrix}
\mu_1^*(\boldsymbol{\phi}^*) q_{11}^* & 
\mu_1^*(\boldsymbol{\phi}^*) q_{12}^* &
\cdots &
\mu_1^*(\boldsymbol{\phi}^*) q_{1K^*}^* \\
\mu_2^*(\boldsymbol{\phi}^*) q_{21}^* & 
\mu_2^*(\boldsymbol{\phi}^*) q_{22}^* &
\cdots &
\mu_2^*(\boldsymbol{\phi}^*) q_{2K^*}^* \\
\vdots & \vdots & \ddots & \vdots \\
\mu_{K^*}^*(\boldsymbol{\phi}^*) q_{K^*1}^* & 
\mu_{K^*}^*(\boldsymbol{\phi}^*) q_{K^*2}^* &
\cdots &
\mu_{K^*}^*(\boldsymbol{\phi}^*) q_{K^*K^*}^*
\end{pmatrix}.
\end{align}
Hence, for any $\alpha \in \mathcal{A}$, there exists at most one ``true'' value in each $\Phi_{K,\alpha}$, as desired.

In fact, note that for $\alpha = (S,W)$ with $S=(S_1, \cdots, S_{K^*})$, the only possibility that there exists no ``true'' value in $\overline{\Phi}_{K,\alpha}$ is when $S_i = \emptyset$ for some $i$. In that case, \eqref{id-1} does not hold since when \eqref{theta-equal} holds, the RHS of \eqref{id-1} has the indicator $1_{\boldsymbol{\theta}_i^*}$, but the LHS has not. Hence, if we define 
\begin{align}\label{set+}
    \mathcal{A}^+ := \left\{ \alpha = (S,W): S = (S_1, \cdots, S_{K^*}) \mbox{ such that } S_k \neq \emptyset \mbox{ for all } k = 1, 2, \cdots, K^* \right\},
\end{align}
then for each $\alpha \in \mathcal{A}^+$, $\overline{\Phi}_{K,\alpha}$ contains exactly one ``true'' value
\begin{equation}\label{phi_alpha*}
    \boldsymbol{\phi}_\alpha^* := \varphi_\alpha(\gamma_\alpha^*),
\end{equation}
with
\begin{equation}\label{gamma_alpha*}
    \gamma_\alpha^* := (R^*; \boldsymbol{\theta}_{s(1)}^*, \boldsymbol{\theta}_{s(2)}^*, \cdots, \boldsymbol{\theta}_{s(K)}^*).
\end{equation}

\begin{remark}
Note that \eqref{id-1} considers the pair $(\boldsymbol{\theta}_0, \boldsymbol{\theta}_1)$ in order to incorporate the information of the transition matrix. This concept has been used in the HMM literature. See, for example, the proof of Lemma 2 in \cite{Leroux:1992}, as well as the composite likelihood method in \cite{chen2016composite}.
\end{remark}

We are now ready to show that this reparametrization satisfies \eqref{change-of-var} - \eqref{CLT-Int}. For \eqref{change-of-var}, by standard change of variables,
\begin{align}\label{change-of-var-full}
    \int_{\Phi_K} p_K({\bm y}_{1:n}|\boldsymbol{\phi}_K)p_0(\boldsymbol{\phi}_K) d\boldsymbol{\phi}_K 
    = \int_{\mathcal{A}} \int_{\Gamma_\alpha} p_K({\bm y}_{1:n}|\varphi_\alpha(\gamma))p_0(\varphi_\alpha(\gamma)) |\Lambda_\alpha(\gamma)| d\gamma d\alpha,   
\end{align}
where $\Lambda_\alpha(\gamma)$ is the Jacobian determinant specified by the mapping $\varphi_\alpha(\gamma)$ (as a function of $(\alpha, \gamma)$.) On the other hand, \eqref{Apply-Lemma1} - \eqref{CLT-Int} are provided by the following four lemmas.

\begin{lemma}\label{Lemma:i} Suppose conditions 1)-5) hold. Then  \eqref{Apply-Lemma1} holds for the $\mathcal{A}^+$ defined in \eqref{set+}.
\end{lemma}

\begin{lemma}\label{Lemma:ii} Suppose conditions 1)-7) hold, and recall $\Delta_\Gamma = \Delta_{K^*} + (K-K^*)d$ defined in \eqref{Delta-Gamma}. Then, for each $\alpha \in \mathcal{A}^+$, 
\begin{enumerate}[a)]
\item $\Gamma_\alpha$ is $\Delta_\Gamma$-dimensional, then the dimension of $\Gamma_\alpha$ remains the same for all $\alpha \in \mathcal{A}^+$;
\item $\overline{\varphi_\alpha(\Gamma_\alpha)} \cap \overline{\Phi}_K^*$ contains exactly one component, which is a true value $\varphi_\alpha(\gamma_\alpha^*) = \boldsymbol{\phi}_\alpha^* \in \Phi_K^*$ with ergodic transition matrix;
\item $\gamma_\alpha^*$ is an interior point of $\Gamma_\alpha$;
\item there exists a constant $C_\alpha$ such that \eqref{CLT-SubSpace} holds.
\end{enumerate}
\end{lemma}

\begin{lemma}\label{Lemma:iv} Suppose conditions 1)-7) hold. Then \eqref{Int-J} holds for the $C_\alpha$ in Lemma \ref{Lemma:ii}.
\end{lemma}

\begin{lemma}\label{Lemma:iii} Suppose conditions 1)-7) hold. Then \eqref{CLT-Int} holds.
\end{lemma}

Now we are ready to prove \eqref{eqn_thm_consist} in Theorem
\ref{theorem:asymptotic_consistency_hmm}.\\

\begin{proof}[Proof of \eqref{eqn_thm_consist}]
Owing to Lemmas \ref{Lemma:i} - \ref{Lemma:iii}, \eqref{Apply-Lemma1} - \eqref{CLT-Int} hold. \eqref{eqn_thm_consist} follows by combining them with \eqref{eqn:CLT-K*} and \eqref{change-of-var-full}.
\end{proof}

It remains to prove Lemmas \ref{Lemma:i}-\ref{Lemma:iii}, which are presented in the Online Supplement.\\

\end{appendix}

\vskip 0.2in
\bibliography{HMMmarginal}

\newpage

\begin{center}
\Large{\bf Online Supplement}
\end{center}

The Online Supplement below contains three parts. Section \ref{appendix:detail_proof} gives the proofs of the key lemmas in Appendix \ref{Proof-over} as well as other theoretical details related to the consistency results. Section \ref{appendix:simulations_normalizing_constant} presents simulation studies for estimating normalizing constants. Section \ref{appendix:simulationrobustness} gives proof of the robustness of the computational algorithm.

\setcounter{section}{0}
\renewcommand{\thesection}{S\arabic{section}}

\section{Detail Proofs for the Consistency Theorems}\label{appendix:detail_proof}

This section contains three parts. Section \ref{appendix:key_lemmas} proves Lemmas \ref{Lemma:i} to \ref{Lemma:iii} in the Appendix, which completes the proof of Theorem \ref{theorem:asymptotic_consistency_hmm}. To simplify the presentation, some of the intermediate steps in the proofs of these Lemmas are deferred to Section \ref{proof:glossary}. Finally, Section \ref{appendix:proofofconsistencygm} proves Corollary \ref{theorem:asymptotic_consistency_gm}.

\subsection{Proofs for the Key Lemmas in Appendix \ref{appendix:proof_theorem_asymptotic_efficiency_hmm}}\label{appendix:key_lemmas}

\begin{proof}[Proof of Lemma \ref{Lemma:i}]
Recall $\Phi_{K,\alpha} = \varphi_\alpha(\Gamma_\alpha)$, so we have
\begin{align}
\notag
    & \int_{\mathcal{A}-\mathcal{A}^+} \int_{\Gamma_\alpha} p_K({\bm y}_{1:n}|\varphi_\alpha(\gamma)) p_0(\varphi_\alpha(\gamma))|\Lambda_\alpha(\gamma)|d\gamma d\alpha \\
\notag
    = & \int_{\cup_{\alpha \in \mathcal{A}-\mathcal{A}^+} \Phi_{K,\alpha}}  p_K({\bm y}_{1:n}|\boldsymbol{\phi}_K) p_0(\boldsymbol{\phi}_K) d\boldsymbol{\phi}_K \\
\label{pf:Lemma-i-1}
     = & \int_{\Phi_K - \cup_{\alpha \in \mathcal{A}^+} \Phi_{K,\alpha}}  p_K({\bm y}_{1:n}|\boldsymbol{\phi}_K) p_0(\boldsymbol{\phi}_K) d\boldsymbol{\phi}_K.
\end{align}
In addition, recall $\Phi_K^*$ and $\overline{\Phi}_K^*$ defined in Section \ref{subsec:ProofConcept}. For any $\delta > 0$ and $\boldsymbol{\phi}_K \in \overline{\Phi}_K$, let 
\begin{equation}\label{B-delta}
\mathcal{B}_\delta(\boldsymbol{\phi}_K) := \left\{ \boldsymbol{\phi}_K: \boldsymbol{\phi}_K \in \Phi_K, \Vert \boldsymbol{\phi} - \boldsymbol{\phi}_K \Vert < \delta \right\},
\end{equation}
and $\Phi_{K,\delta}^* = \bigcup_{\boldsymbol{\phi}_K \in \overline{\Phi}_K^*} \mathcal{B}_\delta(\boldsymbol{\phi}_K)$. In other words, $\mathcal{B}_\delta(\boldsymbol{\phi}_K)$ is a neighborhood of the ``true" values in $\overline{\Phi}_K^*$. Then, by Lemma \ref{lemma:concentrate} provided in Section \ref{proof:glossary}, we have
\begin{equation}\label{pf:Lemma-i-2}
    \frac{\int_{\Phi_K - \Phi_{K,\delta}^*} p_{K}({\bm y}_{1:n} |\boldsymbol{\phi}_K) p_0(\boldsymbol{\phi}_K) d\boldsymbol{\phi}_K }{p_{K^*}({\bm y}_{1:n}|\boldsymbol{\phi}^*)} = O_{P^*} \left( e^{-\epsilon n} \right)
\end{equation}
as $n \rightarrow \infty$. Hence, by \eqref{pf:Lemma-i-1} and \eqref{pf:Lemma-i-2}, we only need to show that $\Phi_{K,\delta}^* \subset \cup_{\alpha \in \mathcal{A}^+} \Phi_{K,\alpha}$ for some $\delta > 0$, and \eqref{Apply-Lemma1} immediately follows.

Recall the $\backsim$ relationship defined in Section \ref{subsec:ProofConcept}. Consider any $\boldsymbol{\phi}_K  = (Q_K; \boldsymbol{\theta}_1, \cdots, \boldsymbol{\theta}_K) \in \overline{\Phi}_K^*$. Let $\{ X_n^1: n \geq 0\}$ be a Markov chain under the law govern by $\boldsymbol{\phi}_K$ and initial distribution $\mu(\boldsymbol{\phi}_K)$, and $\{ X_n^2: n \geq 0\}$ be a Markov chain under the law govern by $\boldsymbol{\phi}^*$ and initial distribution $\mu(\boldsymbol{\phi}^*)$. Since $\boldsymbol{\phi}_K \backsim \boldsymbol{\phi}^*$, we have $\boldsymbol{\theta}_{X_1^1}$ follows the same law as $\boldsymbol{\theta}_{X_1^2}^*$. But
$\boldsymbol{\theta}_{X_1^1}^*$ takes value on ${\bm T} := \{  \boldsymbol{\theta}_k: 1 \leq k \leq K \}$ with $P\{ \boldsymbol{\theta}_{X_1^1} =  \boldsymbol{\theta}_k \} = \mu_k(\boldsymbol{\phi})$ for all $k$, and $\boldsymbol{\theta}_{X_1^2}^*$ takes value on ${\bm T^*} := \{ \boldsymbol{\theta}_i^*: 1 \leq i \leq K^* \}$ with $P\{ \boldsymbol{\theta}_{X_1^2}^* =  \boldsymbol{\theta}_i^* \} = \mu_i(\boldsymbol{\phi}^*) > 0$ for all $i$. So, for $\boldsymbol{\theta}_{X_1^1}$ to follow the same law of $\boldsymbol{\theta}_{X_1^2}^*$, we must have ${\bm T}^* = {\bm T}$, which means that for any $\boldsymbol{\phi}_K$ in $\overline{\Phi}_K^*$ and any $1 \leq i \leq K^*$, there is a $1 \leq \ell_i \leq K$ such that $\boldsymbol{\theta}_{\ell_i} = \boldsymbol{\theta}_i^*$.

In addition, by the construction of $\Theta_i$, we know that there exists $\delta > 0$ such that $\mathcal{B}_\delta(\boldsymbol{\theta}_i^*) \subset \Theta_i$ for all $1 \leq i \leq K^*$. As such, for any  $\boldsymbol{\phi}_K' = (Q_K'; \boldsymbol{\theta}_1', \cdots, \boldsymbol{\theta}_K') \in \mathcal{B}_\delta(\boldsymbol{\phi}_K)$, we have
\begin{equation*}
    \boldsymbol{\theta}_{\ell_i}' \in \mathcal{B}_\delta(\boldsymbol{\theta}_{\ell_i}) = \mathcal{B}_\delta(\boldsymbol{\theta}_i^*) \subset \Theta_i
\end{equation*}
for any $ 1 \leq i \leq K^*$. However, by the construction of  $\Phi_{K,\alpha}$ defined in \eqref{Phi_Kalpha}, this means that $\boldsymbol{\phi}_K'$ must belongs to $\Phi_{K,\alpha}$ for some $\alpha = (S, W)$ with $S = (S_1, S_2, \cdots, S_{K^*})$ satisfying
\begin{equation*}
    \ell_i \in S_i \mbox{ for all } i = 1, 2, \cdots, K^*,
\end{equation*}
which means that $S_i \neq \emptyset$, so that $\alpha \in \mathcal{A}^+$ as defined in \eqref{set+}. As a consequence, we have $\Phi_{K,\delta}^* \subset \cup_{\alpha \in \mathcal{A}^+} \Phi_{K,\alpha}$
as desired. The proof is completed.
\end{proof}

\begin{proof}[Proof of Lemma \ref{Lemma:ii}] a) For any $\alpha = (S,W) \in \mathcal{A}$, it is straightforward to check that $\Phi_{K,\alpha} = \varphi_\alpha(\Gamma_\alpha)$ and $\Gamma_\alpha$ has the same dimension. In addition, by \eqref{varphi-inv}, it is straightforward to show that
\begin{equation}\label{Gamma-decomp}
 \Gamma_\alpha = \mathcal{R}_\alpha \times \Theta_{s(1)} \times \cdots \times \Theta_{s(K)},
\end{equation}
where
\begin{equation}\label{set-R}
\mathcal{R}_\alpha = \{ R: \varphi_\alpha(R; \boldsymbol{\theta}_{s(1)}^*, \cdots, \boldsymbol{\theta}_{s(K)}^*) \in \Phi_{K,\alpha} \};
\end{equation}
namely, the projection space of $\varphi_\alpha^{-1}(\Phi_{K,\alpha})$ onto the space of $R$.

Now, each $\Theta_{s(k)}$ has dimension $d$. As for the dimension of $\mathcal{R}_\alpha$, note that, by \eqref{psi-inv}, $\psi_\alpha^{-1}$ maps $Q_K$ to a $K^*$-by-$K^*$ matrix satisfying
\begin{align}\label{set-R-dim}
    \sum_{j=1}^{K^*} \psi_{\alpha i j}^{-1}(Q_K) &
    = \sum_{j=1}^{K^*} \sum_{k \in S_i} \sum_{\ell \in S_j} \mu_k q_{k \ell} 
    = \sum_{k \in S_i} \sum_{\ell \in \cup_{j=1}^{K^*} S_j} \mu_k q_{k \ell} 
     = \sum_{k \in S_i} \sum_{\ell = 1}^K \mu_k q_{k \ell} = \sum_{k \in S_i} \mu_k,
\end{align}
since $S=(S_1, \cdots, S_{K^*})$ is a partition of $\{1, 2, \cdots, K\}$, and $\sum_{\ell=1}^K q_{k\ell} = 1$. Hence, $\mathcal{R}_\alpha$ is the space of $K^*$-by-$K^*$ matrices subject to $K^*$ constraints given by \eqref{set-R-dim}, so it has dimension $(K^*)^2 - K^* = K^*(K^*-1)$. Thus, by \eqref{Gamma-decomp}, the dimension of $\Gamma_\alpha$ is
\begin{align*}
    K^*(K^*-1) + dK = K^*(K^*-1)+dK^* + d(K-K^*) 
    = \Delta_{K^*} + d(K-K^*) = \Delta_\Gamma,
\end{align*}
and so does $\Phi_{K,\alpha} = \varphi_\alpha(\Gamma_\alpha)$.

b) By the argument leading up to \eqref{Phi_Kalpha} and \eqref{gamma_alpha*}, $\overline{\varphi_\alpha(\Gamma_\alpha)} \cap \overline{\Phi}_K^* = \left\{ \varphi_\alpha(\gamma_\alpha^*) \right\}$ defined as in \eqref{gamma_alpha*}. In addition, a direct check shows that $\varphi_\alpha(\gamma_\alpha^*) = \boldsymbol{\phi}_\alpha^* = (Q_K; \boldsymbol{\theta}_1, \cdots, \boldsymbol{\theta}_K)$ with $Q_K \in \mathcal{Q}_K$, which is an ergodic transition matrix.

c) Recall $\gamma_\alpha^* = (R^*; \boldsymbol{\theta}_{s(1)}, \cdots, \boldsymbol{\theta}_{s(K)})$. A direct check shows that $\psi_\alpha(R^*)$ is a transition matrix with positive entries, and $\boldsymbol{\theta}_i^*$ is an interior of $\Theta_i$ for all $i=1,2,\cdots,K^*$. As such, $\varphi_\alpha(\gamma_\alpha^*)$ is an interior point of $\Phi_{K,\alpha}$, and therefore, $\gamma_\alpha^*$ is an interior point of $\Gamma_\alpha$.

d) Note that since $\boldsymbol{\phi}^* \backsim \boldsymbol{\phi}_\alpha^* = \varphi_\alpha(\gamma_\alpha^*)$, by Lemma 2 in \cite{LerouxM:1992}, $p_{K^*}({\bm y}_{1:n} | \boldsymbol{\phi}^*) = p_K({\bm y}_{1:n} | \varphi_\alpha(\gamma_\alpha^*))$ and hence,
\begin{align}\label{pf-lemma-ii-d-1}
     \mbox{LHS of \eqref{CLT-SubSpace}} 
    =  \frac{n^{\frac{\Delta_\Gamma}{2}}  \int_{\Gamma_\alpha} p_K({\bm y}_{1:n}|\varphi_\alpha(\gamma)) p_0(\varphi_\alpha(\gamma)) |\Lambda_\alpha(\gamma)| d\gamma}{p_K({\bm y}_{1:n}|\varphi_\alpha(\gamma_\alpha^*))}.
\end{align}
In addition, since $\varphi_\alpha(\gamma)$ is bijective, the inverse function theorem ensures that $|\Lambda_\alpha(\gamma)| > 0$ for all $\alpha \in \mathcal{A}^+$ and $\gamma \in \Gamma_\alpha$. Thus, if $p_0(\varphi_\alpha(\gamma_\alpha^*)) > 0$, then
\begin{equation}\label{p-alpha}
p_\alpha(\cdot) := \frac{p_0(\varphi_\alpha(\cdot))|\Lambda_\alpha(\cdot)|}
{\int_{\Gamma_\alpha}  p_0(\varphi_\alpha(\gamma))|\Lambda_\alpha(\gamma)| d\gamma}
\end{equation}
is a prior distribution on $\Gamma_\alpha$ with positive density at $\gamma_\alpha^*$, and therefore
\begin{align}\label{pf-lemma-ii-d-2}
    \mbox{RHS of \eqref{pf-lemma-ii-d-1}} 
    = \frac{n^{\frac{\Delta_\Gamma}{2}}  \int_{\Gamma_\alpha} p_K({\bm y}_{1:n}|\varphi_\alpha(\gamma)) p_\alpha(\gamma) d\gamma}{p_K({\bm y}_{1:n}|\varphi_\alpha(\gamma_\alpha^*)) p_\alpha(\gamma_\alpha^*)} \times p_0(\varphi_\alpha(\gamma_\alpha^*)) |\Lambda_\alpha(\gamma_\alpha^*)|.
\end{align}
Combining \eqref{pf-lemma-ii-d-2} with the Bernstein-von Mises theorem  in Lemma \ref{Lemma:gBvM} below to have
\begin{align*}
    \mbox{RHS of \eqref{pf-lemma-ii-d-2}} 
    = \frac{(2\pi)^{\Delta_\Gamma/2}}{|J_\alpha|^{1/2}} \times p_0(\varphi_\alpha(\gamma_\alpha^*)) |\Lambda_\alpha(\gamma_\alpha^*)| + o_{P^*}(1) = C_\alpha + o_{P^*}(1),
\end{align*}
where $J_\alpha$ is the Fisher information to be defined in Lemma \ref{Lemma:gBvM}, and
\begin{equation}\label{C_alpha}
    C_\alpha := \frac{(2\pi)^{\Delta_\Gamma/2}}{|J_\alpha|^{1/2}} \times p_0(\varphi_\alpha(\gamma_\alpha^*)) |\Lambda_\alpha(\gamma_\alpha^*)|.
\end{equation}

It remains to handle the case when $p_0(\varphi_\alpha(\gamma_\alpha^*))=0$. For this purpose, let $u_\alpha$ be any continuous prior distribution on $\Gamma_\alpha$ such that $u_\alpha(\gamma_\alpha^*) > 0$. Then, for any $\epsilon > 0$, \begin{equation}\label{p-alpha-epsilon}
p_\alpha^\epsilon(\cdot) := \frac{p_0(\varphi_\alpha(\cdot))|\Lambda_\alpha(\cdot)| + \epsilon u_\alpha(\cdot) }
{\int_{\Gamma_\alpha}  \left( p_0(\varphi_\alpha(\gamma))|\Lambda_\alpha(\gamma)| + \varphi_\alpha(\gamma) \right) d\gamma}
\end{equation}
is a continuous prior distribution on $\Gamma_\alpha$ with positive density at $\gamma_\alpha^*$. Hence, by Lemma \ref{Lemma:gBvM},
\begin{align}
\notag
    \mbox{RHS of \eqref{pf-lemma-ii-d-1}} 
\leq & \frac{n^{\frac{\Delta_\Gamma}{2}}  \int_{\Gamma_\alpha} p_K({\bm y}_{1:n}|\varphi_\alpha(\gamma)) \left( p_0(\varphi_\alpha(\gamma)) |\Lambda_\alpha(\gamma)| + \varphi_\alpha(\gamma) \right) d\gamma}{p_K({\bm y}_{1:n}|\varphi_\alpha(\gamma_\alpha^*))} \\
\notag
    \leq & \frac{n^{\frac{\Delta_\Gamma}{2}}  \int_{\Gamma_\alpha} p_K({\bm y}_{1:n}|\varphi_\alpha(\gamma)) p_\alpha^\epsilon(\gamma) d\gamma}{p_K({\bm y}_{1:n}|\varphi_\alpha(\gamma_\alpha^*)) p_\alpha^\epsilon(\gamma_\alpha^*)} \times \varphi_\alpha(\gamma_\alpha^*) |\Lambda_\alpha(\gamma_\alpha^*)|  \\
    & \xrightarrow[n \rightarrow \infty]{P^*} \frac{(2\pi)^{\Delta_\Gamma/2}}{|J_\alpha|^{1/2}} \times \varphi_\alpha(\gamma_\alpha^*) |\Lambda_\alpha(\gamma_\alpha^*)|.
\end{align}
Since $\epsilon > 0$ can be arbitrarily small, we get
\begin{equation*}
    \mbox{RHS of \eqref{pf-lemma-ii-d-1}} = 0+o_p(1) = C_\alpha + o_p(1),
\end{equation*}
since $C_\alpha = 0$ due to \eqref{C_alpha} and the fact that $p_0(\varphi_\alpha(\gamma_\alpha^*)) = 0$ in this case. The proof is completed.
\end{proof}

\begin{proof}[Proof of Lemma \ref{Lemma:iv}] 
Recall $\mathcal{Q}_K^\epsilon$ in condition 7),   $\psi_\alpha$ defined in \eqref{psi}, and $\mathcal{R}_\alpha$ defined in \eqref{set-R}. Let 
\begin{equation}\label{A_ep}
\mathcal{A}_\epsilon^+ := \left\{ \alpha \in \mathcal{A}^+: \psi_\alpha(\mathcal{R}_\alpha) \cap Q_K^\epsilon \neq \emptyset \right\},
\end{equation}
which is the collection of $\alpha$ such that there exists some $\boldsymbol{\phi}_K = (Q_K; \boldsymbol{\theta}_1, \cdots, \boldsymbol{\theta}_K) \in \Phi_{K,\alpha}$ with $Q_K \in \mathcal{Q}_K^\epsilon$. Note that for any $\alpha \in \mathcal{A}^+ - \mathcal{A}_\epsilon^+$, we have
\begin{equation}\label{varphi-R-out}
\psi_\alpha(\mathcal{R}_\alpha) \subset \mathcal{Q}_K - \mathcal{Q}_K^\epsilon.
\end{equation}
In addition, by condtion 7), $p_0(\boldsymbol{\phi}_K)=0$ for any $\boldsymbol{\phi}_K = (Q_K; \boldsymbol{\theta}_1, \cdots, \boldsymbol{\theta}_K)$ with $Q_K \in \mathcal{Q}_K - \mathcal{Q}_K^\epsilon$. So by \eqref{varphi-R-out}, 
\begin{equation}\label{prior-out}
\psi_\alpha(\gamma_\alpha^*) \in \psi_\alpha(\mathcal{R}_\alpha) \subset \mathcal{Q}_K - \mathcal{Q}_K^\epsilon \Rightarrow p_0(\varphi_\alpha(\gamma_\alpha^*)) = 0.
\end{equation}
Hence, by \eqref{C_alpha}, we have
\begin{equation}\label{C_alpha_zero}
    C_\alpha = 0 \mbox{ for any } \alpha \in \mathcal{A}^+ - \mathcal{A}_\epsilon^+,
\end{equation}
and hence
\begin{equation*}
    \int_{\alpha \in \mathcal{A}^+} C_\alpha d\alpha = \int_{\alpha \in \mathcal{A}_\epsilon^+} C_\alpha  d\alpha.
\end{equation*}
In addition, since $\mathcal{Q}_K^\epsilon$ is compact and $\psi_\alpha(\gamma)$ is continuous, $\mathcal{A}_\epsilon^+$ is compact. Furthermore since $C_\alpha $ is continuous in $\alpha$, we have
\begin{equation*}
    \int_{\alpha \in \mathcal{A}_\epsilon^+} C_\alpha  d\alpha < \infty,
\end{equation*}
which completes the proof.
\end{proof}

\begin{remark}\label{remark:prior-conj}
As one can see from above, to ensure the finiteness in \eqref{Int-J}, one does not actually need a stronger condition 7) but only needs to make the prior $p_0$ vanishes quick enough when it approaches the boundary of $\mathcal{Q}_K$. {\color{blue}To be more specific, if one can estimate the rate of $|J_\alpha|$ as $\alpha$ approaches the boundary, then one can choose a $p_0$ that vanishes quick enough near the boundary to make \eqref{Int-J} holds. Potential choices of priors include the Dirichlet distribution with a certain order or priors that vanish exponentially fast near the boundary. See \cite{Gassiat2014} for examples.

However, {\color{blue}priors with non-zero values near the boundary} would raise a question of uniform convergence in Lemma \ref{Lemma:iii}. See Remark \ref{remark-Lemma:L''} for details.}
\end{remark}

\begin{proof}[Proof of Lemma \ref{Lemma:iii}] 
Recall $\mathcal{A}_\epsilon^+$ defined in \eqref{A_ep}. By \eqref{C_alpha_zero}, we have $C_\alpha=0$ for any $\alpha \in \mathcal{A}^+ - \mathcal{A}_\epsilon^+$. In addition, through the same argument leading up to \eqref{prior-out}, for any  $\alpha \in \mathcal{A}^+ - \mathcal{A}_\epsilon^+$ and any $\gamma \in \Gamma_\alpha$, we have
\begin{equation*}
    \varphi_\alpha(\gamma) \notin \Phi_K^\epsilon \Rightarrow p_0(\varphi_\alpha(\gamma)) = 0,
\end{equation*}
which means that 
\begin{equation*}
    \frac{\int_{\Gamma_\alpha} p_K({\bm y}_{1:n}|\varphi_\alpha(\gamma)) p_0(\varphi_\alpha(\gamma))|\Lambda_\alpha(\gamma)|d\gamma}{p_{K^*}({\bm y}_{1:n}|\boldsymbol{\phi}^*)} = 0 = C_\alpha.
\end{equation*}
In other words, as soon as we show that the convergence in \eqref{CLT-SubSpace} is uniform across $\mathcal{A}_\epsilon^+$, we immediately know that the convergence is also uniform across $\mathcal{A}^+$, and hence, \eqref{CLT-Int} holds.

To see why the convergence is uniform for $\alpha \in \mathcal{A}_\epsilon^+$, let us first review how the convergence holds for any fixed $\alpha$. Recall $L_K({\bm y}_{1:n}|\boldsymbol{\phi}_K) = \log p_K({\bm y}_{1:n}|\boldsymbol{\phi}_K)$, and $p_\alpha$ defined in \eqref{p-alpha}. Note that the convergence in \eqref{CLT-SubSpace} comes from Lemma \ref{Lemma:gBvM} below, which uses the Berstein-von Mises theorem in \cite{HMMCLT}, Theorem 3.1. This theorem is obtained through the following steps:
\begin{enumerate}[a)]
    \item Use the law of large numbers (LLN) for $L_K({\bm y}_{1:n}|\varphi_\alpha(\gamma))$ for any $\gamma \in \Gamma_\alpha$ to show that, for any $\delta > 0$, as $n \rightarrow \infty$,
    \begin{align}\label{HMMCLT-L}
        \frac{\int_{\Gamma_\alpha} p_K({\bm y}_{1:n}|\varphi_\alpha(\gamma)) p_\alpha(\gamma) d\gamma}{p_K({\bm y}_{1:n}|\varphi_\alpha(\gamma_\alpha^*))} 
        =
        \frac{\int_{ \mathcal{B}_{\delta,\alpha}(\gamma_\alpha^*)} p_K({\bm y}_{1:n}|\varphi_\alpha(\gamma)) p_\alpha(\gamma) d\gamma}{p_K({\bm y}_{1:n}|\varphi_\alpha(\gamma_\alpha^*))}  + o_{P^*}(1),
    \end{align}
    where $\mathcal{B}_{\delta,\alpha}(\gamma) = \left\{ \gamma' \in \Gamma_\alpha: \Vert \gamma' - \gamma \Vert < \delta \right\}$.
    \item On $\mathcal{B}_{\delta,\alpha}(\gamma_\alpha^*)$, by Taylor expansion of $L_K$ to have
    \begin{align}\label{K-Taylor}
    \notag
    L_{K} \left( {\bm y}_{1:n} \bigg\vert \varphi_\alpha(\gamma_\alpha^*) + \frac{v}{\sqrt{n}} \right) 
    = & L_{K}({\bm y}_{1:n}|\varphi_\alpha(\gamma_\alpha^*)) + v^t \frac{D_\gamma L_{K}({\bm y}_{1:n}|\varphi_\alpha(\gamma_\alpha^*))}{\sqrt{n}} \\
    & + v^t \frac{D_\gamma^2L_{K}({\bm y}_{1:n}|\varphi_\alpha(\tilde{\gamma}_{\alpha,n}))}{2n} v,
    \end{align}
    for some $\tilde{\gamma}_{\alpha,n}$ between $\varphi_{\alpha}(\gamma_\alpha^*)$ and $\varphi_{\alpha}(\gamma_\alpha^*) + \frac{v}{\sqrt{n}}$, with $D_\gamma f$ and $D_\gamma^2 f$ denote the gradient and the Hessian matrix of a function $f$ with respect to $\gamma$, respectively.
    
    \item Use the central limit theorem (CLT) for $D_\gamma L_K({\bm y}_{1:n}|\varphi_\alpha(\gamma_\alpha^*))$ to control the first derivative term in \eqref{K-Taylor}.
    \item Use the LLN for $D_\gamma^2 L_K({\bm y}_{1:n}|\varphi_\alpha(\tilde{\gamma}_{\alpha,n})$ with $\tilde{\gamma}_{\alpha,n} \xrightarrow{P^*} \gamma_\alpha^*$ to control the second derivative term in \eqref{K-Taylor}. Note that a) ensures the convergence of $\tilde{\gamma}_{\alpha,n}$.
\end{enumerate}

Thus, if we want to show that the convergence in \eqref{CLT-SubSpace} is uniform across $\mathcal{A}_\epsilon^+$, we need the following three steps:
\begin{enumerate}[i)]
    \item show that for any sufficiently small $\delta>0$, \eqref{HMMCLT-L} holds uniformly for all $\alpha \in \mathcal{A}_\epsilon^+$;
    
    \item show that the CLT in step c) is uniform for all $\alpha \in \mathcal{A}_\epsilon^+$. More precisely, we need to show that 
    \begin{equation}\label{eqn:L'}
    \frac{1}{\sqrt{n}} \sup_{\alpha \in \mathcal{A}_\epsilon^+} \Big\Vert D_\gamma L_K({\bm y}_{1:n}|\varphi_\alpha(\gamma_\alpha^*)) J_\alpha^{-1/2} \Big\Vert = O_{P^*}(1),
    \end{equation}
    where $J_\alpha$ will be defined in Lemma \ref{Lemma:gBvM}.
    
    \item show that the LLN in step d) is uniform for all $\alpha \in \mathcal{A}_\epsilon^+$. That is, we need to show that, for any $\gamma_{\alpha,n}$ with $ \sup_{\alpha \in \mathcal{A}_\epsilon^+} \Vert \gamma_{\alpha,n} - \gamma_\alpha^* \Vert \rightarrow 0$,
    \begin{equation}\label{eqn:L''}
        \sup_{\alpha \in \mathcal{A}_\epsilon^+} \Big\Vert \frac{1}{n} D_\gamma^2 L_K({\bm y}_{1:n}|\varphi_\alpha(\gamma_{\alpha,n})) + J_\alpha \Big\Vert \xrightarrow[n \rightarrow \infty]{P^*} 0.
    \end{equation}
\end{enumerate}

We will prove i)-iii) through Lemma \ref{Lemma:L} - \ref{lemma:L'}, respectively, which leads to the uniform convergence of \eqref{CLT-SubSpace} and completes the proof.
\end{proof}

\subsection{Lemmas for Intermediate Steps}\label{proof:glossary}

\begin{lemma}[Likelihood Concentration Used in the Proof of Lemma \ref{Lemma:i}]\label{lemma:concentrate} Assume that conditions 1)-5) hold. Then, for any $K \geq K^*$  and $\delta >0$, there exists $\epsilon > 0$ such that
\begin{equation}\label{L-concentration}
P^*\left\{ \sup_{\boldsymbol{\phi} \in \Phi_K - \Phi_{K,\delta}^*} \frac{1}{n}(L_K({\bm y}_{1:n}|\boldsymbol{\phi}) - L_{K^*}({\bm y}_{1:n}|\boldsymbol{\phi}^*)) \leq -\epsilon \right\}  \xrightarrow[n \rightarrow \infty]{} 1,
\end{equation}
where $\Phi_{K,\delta}^*$ is defined right after \eqref{B-delta}. Consequently, as $n \rightarrow \infty$,
\begin{equation}\label{eqn:concentrate}
    \frac{\int_{\Phi_K - \Phi_{K,\delta}^*} p_{K}({\bm y}_{1:n} |\boldsymbol{\phi}_K) p_0(\boldsymbol{\phi}_K) d\boldsymbol{\phi}_K }{p_{K^*}({\bm y}_{1:n}|\boldsymbol{\phi}^*)} = O_{P^*} \left( e^{-\epsilon n} \right).
\end{equation}
\end{lemma}
Lemma \ref{lemma:concentrate} is used to prove \eqref{Apply-Lemma1} as we mentioned in Section \ref{Proof-over}.

\begin{proof}[Proof of Lemma~\ref{lemma:concentrate}]

Define $\tilde{Q}_K^* = \left\{ \tilde{q}_{k\ell}^*, 1 \leq k,\ell \leq K \right\}$ by
\begin{equation*}
    \tilde{q}_{k\ell}^* = 
    \begin{cases}
    q_{k\ell}^* & k,\ell < K^* \\
    \frac{1}{K-K^*+1}q_{kK^*}^* & k < K^*, \ell \geq K^*\\ 
    q_{K\ell}^* & k \geq K^*, \ell < K^*\\
    \frac{1}{K-K^*+1}q_{KK^*}^* & k, \ell \geq K^*\\ 
    \end{cases}
\end{equation*}
which is the transition matrix when we ``split" state $K^*$ into states $K^*, K^*+1, \cdots, K$ with equal probabilities. Further, define
\begin{equation*}
    \tilde{\boldsymbol{\phi}}_K^* = (\tilde{Q}_K^*; \boldsymbol{\theta}_1, \cdots, \boldsymbol{\theta}_{K^*-1}, \boldsymbol{\theta}_{K^*}, \boldsymbol{\theta}_{K^*}, \cdots, \boldsymbol{\theta}_{K^*}) \in \Phi_K.
\end{equation*}
Then, a direct computation shows that 
\begin{equation*}
    p_{K}({\bm y}_{1:n}|\tilde{\boldsymbol{\phi}}_K^*) = p_{K^*}({\bm y}_{1:n}|\boldsymbol{\phi}^*),
\end{equation*}
and therefore, $L_{K^*}({\bm y}_{1:n} | \boldsymbol{\phi}^*) = L_{K}({\bm y}_{1:n} | \tilde{\boldsymbol{\phi}}_K^*)$. In addition, let $\tilde{P}_K^*$ denotes the probability when $\{X_i, i \geq 0\}$ is a Markov chain on $\mathcal{X}_K$, governed by $\tilde{Q}_K^*$, with $X_0$ follows the invariant measure under $\tilde{Q}_K^*$. Then, a direct computation shows that $\{Y_i, i \geq 1\}$ has the same law under $P^*$ and $\tilde{P}_K^*$. As such, for any $\delta > 0$ and $\epsilon > 0$, we have
\begin{align}\label{pf-concentration-2}
\notag
& P^*\left\{ \sup_{\boldsymbol{\phi} \in \Phi_K - \Phi_{K,\delta}^*} \frac{L_K({\bm y}_{1:n}|\boldsymbol{\phi}) - L_{K^*}({\bm y}_{1:n}|\boldsymbol{\phi}^*)}{n} \leq -\epsilon \right\} \\
= &
\tilde{P}_K^*\left\{ \sup_{\boldsymbol{\phi} \in \Phi_K - \Phi_{K,\delta}^*} \frac{L_K({\bm y}_{1:n}|\boldsymbol{\phi}) - L_K({\bm y}_{1:n}|\tilde{\boldsymbol{\phi}}_K^*)}{n} \leq -\epsilon \right\}. 
\end{align}
However, by using the argument in the proof of Theorem 3 in \cite{LerouxM:1992}, we have that under conditions 1)-5), for any $\delta > 0$, there exists $\epsilon > 0$ such that the second line in \eqref{pf-concentration-2} goes to one as $n \rightarrow \infty$, so \eqref{L-concentration} is proven. \eqref{eqn:concentrate} immediately follows since, by \eqref{L-concentration}, with probability approaching to one, we have, as $n \rightarrow \infty$,
\begin{align*}
    \notag
    & \frac{\int_{\Phi_K - \Phi_{K,\delta}^*} p_{K}({\bm y}_{1:n} |\boldsymbol{\phi}_K) p_0(\boldsymbol{\phi}_K) d\boldsymbol{\phi}_K }{p_{K^*}({\bm y}_{1:n}|\boldsymbol{\phi}^*)}
    = \frac{ \int_{\Phi_K - \Phi_{K,\delta}^*} \exp \{ L_{K}({\bm y}_{1:n} |\boldsymbol{\phi}_K) \} p_0(\boldsymbol{\phi}_K) d\boldsymbol{\phi}_K }{ \exp \{L_{K^*}({\bm y}_{1:n}|\boldsymbol{\phi}^*) \}} \\
    \leq &
       \exp \left\{ \sup_{\boldsymbol{\phi}_K \in \Phi_K - \Phi_{K,\delta}^*} \frac{L_K({\bm y}_{1:n}|\boldsymbol{\phi}_K) - L_{K^*}({\bm y}_{1:n}|\boldsymbol{\phi}^*)}{n} \times n \right\} 
    = O_{P^*}\left( e^{- \epsilon n} \right).
\end{align*}
\end{proof}

\begin{remark}
Note that the results in \cite{LerouxM:1992} do not require the probability measure to be under the true number of states $K^*$. It also does not require the true parameters to be unique in the parameter space, as its results are with respect to the quotient topology of the equivalence class. See the last paragraph on page 142 of \cite{LerouxM:1992}. Also note that its definition of equivalent class $\boldsymbol{\phi} \backsim \boldsymbol{\phi}'$ does not require $\boldsymbol{\phi}$ and $ \boldsymbol{\phi}'$ to correspond to the same number of states; it only requires that $\{ \boldsymbol{\theta}_{X_i}: i \geq 0\}$ follow the same law under $\boldsymbol{\phi}$ and $\boldsymbol{\phi}'$.
\end{remark}

\begin{lemma}[Generalized Bernstein–von Mises theorem for HMM]\label{Lemma:gBvM}
Assume conditions 1)-5) holds and $K \geq K^*$. Then, for any $\alpha \in \mathcal{A}^+$ and any continuous prior distribution $p_\alpha$ on $\Gamma_\alpha$ with $p_\alpha(\gamma_\alpha^*) > 0$,
\begin{align}\label{eqn:Lemma:gBVM}
     \frac{n^{\frac{\Delta_\Gamma}{2}}  \int_{\Gamma_\alpha} p_K({\bm y}_{1:n}|\varphi_\alpha(\gamma)) p_\alpha(\gamma) d\gamma}{p_K({\bm y}_{1:n}|\varphi_\alpha(\gamma_\alpha^*)) p_\alpha(\gamma_\alpha^*)} 
     \xrightarrow[n \rightarrow \infty]{P^*} \frac{(2\pi)^{\Delta_\Gamma/2}}{|J_\alpha|^{1/2}},
\end{align}
where
\begin{equation*}
    J_\alpha = \lim_{n \rightarrow \infty} \frac{1}{n} E^*\left[ D_\gamma^2 \log p_K ({\bm y}_{1:n} | \varphi_\alpha(\gamma_\alpha^*)) \right].
\end{equation*}
\end{lemma}

\begin{proof}
Let $P_\alpha^*$ and $E_\alpha^*$ denote the probability and expectation when $\{Y_i, i \geq 1 \}$ is govern by the $K$-state HMM under $\boldsymbol{\phi}_\alpha^* = \varphi_\alpha(\gamma_\alpha^*)$, respectively. Since $\boldsymbol{\phi}_\alpha^* \backsim \boldsymbol{\phi}^*$,  by Lemma 2 in \cite{LerouxM:1992}, $\{Y_i, i \geq 1 \}$ has the same probability law under $P^*$ and $ P_\alpha^*$, which implies that
\begin{align*}
    \lim_{n \rightarrow \infty} \frac{1}{n} E^*\left[ D_\gamma^2 \log p_K ({\bm y}_{1:n} | \varphi_\alpha(\gamma_\alpha^*) \right] 
    = 
    \lim_{n \rightarrow \infty} \frac{1}{n} E_\alpha^*\left[ D_\gamma^2 \log p_K ({\bm y}_{1:n} | \varphi_\alpha(\gamma_\alpha^*) \right],
\end{align*}
with the limit in the last line exists due to Lemma 2 of \cite{Bickel:1998} under conditions 1)-5). Hence, $J_\alpha$ is well defined.

Now, we would like to show that the Berstein-von Mises theorem for HMM in Theorem 3.1 of \cite{HMMCLT} holds on $\Gamma_\alpha$. To do so, we need to check the conditions (B1)-(B6) in \cite{HMMCLT}, which are true since
\begin{itemize}
    \item (B1) holds since, by part c) of Lemma \ref{Lemma:ii}, $\gamma_\alpha^*$ is an interior point of $\Gamma_\alpha$;
    \item (B2)-(B5) are implied by conditions 1)-5);
    \item (B6) holds since, by part b) of Lemma \ref{Lemma:ii}, $\overline{\varphi_\alpha(\Gamma_\alpha)} \cap \overline{\Phi}_K^*$ contains exactly one point $\boldsymbol{\phi}_\alpha^* = \varphi_\alpha(\gamma_\alpha^*)$.
\end{itemize}
Thus, Theorem 3.1 in \cite{HMMCLT} gives 
\begin{align*}
     \frac{n^{\frac{\Delta_\Gamma}{2}}  \int_{\Gamma_\alpha} p_K({\bm y}_{1:n}|\varphi_\alpha(\gamma)) p_\alpha(\gamma) d\gamma}{p_K({\bm y}_{1:n}|\varphi_\alpha(\gamma_\alpha^*)) p_\alpha(\gamma_\alpha^*)} 
     \xrightarrow[n \rightarrow \infty]{P_\alpha^*} \frac{(2\pi)^{\Delta_\Gamma/2}}{|J_\alpha|^{1/2}},
\end{align*}
which immediately implies \eqref{eqn:Lemma:gBVM} since $\{Y_i, i \geq 1\}$ has the same law under $P^*$ and $ P_\alpha^*$. The proof is completed.
\end{proof}

\begin{remark}
For the results cited in the proof of Lemma \ref{Lemma:gBvM}:
\begin{itemize}
    \item \cite{LerouxM:1992} studies consistency of the MLE for HMM. Its Lemma 2 shows that two probability laws agree if and only if their corresponding parameters are in the same equivalent class.
    \item \cite{Bickel:1998} studies asymptotic normality of the MLE for HMM. Its Lemma 2 shows that the corresponding Fisher information exists.
    \item Theorem 3.1 in \cite{HMMCLT} is the Berstein-von Mises theorem for HMM under the identifiability condition.
\end{itemize}
\end{remark}

\begin{lemma}[Uniform Concentration of $L_K$]\label{Lemma:L} Assume conditions 1)--5). Then for any $\delta > 0$, \eqref{HMMCLT-L} holds uniformly for all $\alpha \in \mathcal{A}_\epsilon^+$.
\end{lemma}

\begin{proof}
As in the proof of Lemma \ref{Lemma:ii}, part d), we have
$L_K({\bm y}_{1:n}|\varphi_\alpha(\gamma_\alpha^*)) = L_{K^*}({\bm y}_{1:n}|\boldsymbol{\phi}^*)$. For any $\delta > 0$, define
\begin{equation*}
    \mathcal{V}_\delta = \cup_{\alpha \in \mathcal{A}_\epsilon^+} \varphi_\alpha \left( \Gamma_\alpha - \mathcal{B}_{\delta, \alpha}(\gamma_\alpha^*) \right).
\end{equation*}
If we can show that, for any $\delta > 0$, $\mathcal{V}_\delta$ is bounded away from $\overline{\Phi}_K^*$ by some constant $c_\delta > 0$, then we have
\begin{equation}\label{V-subset}
\mathcal{V}_{\delta} \subset \Phi_K - \Phi_{K,c_\delta}^*,
\end{equation}
with $\Phi_{K,\delta}^*$ being defined in the paragraph right before Lemma \ref{lemma:concentrate}. In this case, Lemma \ref{lemma:concentrate} shows that there exists $c > 0$ such that
\begin{align*}
    & P^* \left\{ \sup_{\boldsymbol{\phi} \in \mathcal{V}_\delta}
    \frac{L_K({\bm y}_{1:n}|\boldsymbol{\phi}) - L_K({\bm y}_{1:n}|\varphi_\alpha(\gamma_\alpha^*))}{n} \leq -c
    \right\} \\
    = & P^* \left\{ \sup_{\boldsymbol{\phi} \in \mathcal{V}_\delta}
    \frac{L_K({\bm y}_{1:n}|\boldsymbol{\phi}) - L_{K^*}({\bm y}_{1:n}|\boldsymbol{\phi}^*)}{n} \leq -c
    \right\} \\
    \geq & P^* \left\{ \sup_{\boldsymbol{\phi} \in \Phi_K 
    - \Phi_{K,c_\delta}^*}
    \frac{L_K({\bm y}_{1:n}|\boldsymbol{\phi}) - L_{K^*}({\bm y}_{1:n}|\boldsymbol{\phi}^*)}{n} \leq -c
    \right\} \xrightarrow[n \rightarrow \infty]{} 1,
\end{align*}
and the uniformity of \eqref{HMMCLT-L} immediately follows through an argument similar to that in Lemma \ref{lemma:concentrate}.

To see why such $c_\delta >0$ exists, first note that, as shown in the proof of Lemma 4, $\mathcal{A}_\epsilon^+$ is compact. Further note that
\begin{equation}\label{Phi_K*_decomp}
    \overline{\Phi}_K^* = \Phi_K^* \cup \left( \overline{\Phi}_K^* - \Phi_K^* \right).
\end{equation}
For $\Phi_K^*$, note that for each $\alpha \in \mathcal{A}_\epsilon^+$ and $\gamma \in \Gamma_\alpha - \mathcal{B}_{\delta,\alpha}(\gamma_\alpha^*)$, the closest $\boldsymbol{\phi} \in \Phi_K^*$ to $\varphi_\alpha(\gamma)$ is $\boldsymbol{\phi}_\alpha^* = \varphi_\alpha(\gamma_\alpha^*)$. Hence, $\mathcal{V}_\delta$ is bounded away from $\Phi_K^*$ by
\begin{equation}\label{c_gamma}
    c_\delta' := \inf \left\{ \Vert \varphi_\alpha(\gamma) - \varphi_\alpha(\gamma_\alpha^*) \Vert: \alpha \in  \mathcal{A}_\epsilon^+, \gamma \in \Gamma_\alpha - \mathcal{B}_{\delta,\alpha}(\gamma_\alpha^*) \right\}.
\end{equation}
As for $\overline{\Phi}_K^* - \Phi_K^*$, note that
\begin{equation*}
    \overline{\Phi}_K^* - \Phi_K \subset \partial \mathcal{Q}_K \times \Theta^K.
\end{equation*}
On the other hand, since $\mathcal{A}_\epsilon^+$ is compact, $\cup_{\alpha \in \mathcal{A}_\epsilon^+} \psi_\alpha(\Gamma)$ is a compact subset of the open set $\mathcal{Q}_K$, which means that $\cup_{\alpha \in \mathcal{A}_\epsilon^+} \psi_\alpha(\Gamma)$ is bounded away from $\partial \mathcal{Q}_K$ by some $c_Q > 0$. Combining with \eqref{Phi_K*_decomp} and \eqref{c_gamma}, we see that $\mathcal{V}_\delta$ is bounded away from $\overline{\Phi}_K^*$ by $c_\delta = \min \{ c_\delta', c_Q\} > 0$, so \eqref{V-subset} holds. The proof is completed.
\end{proof}

\begin{remark}\label{remark-Lemma:L}
Note that this argument will not work when we replace $\mathcal{A}_\epsilon^+$ by $\mathcal{A}^+$, as $c_Q = 0$ in this case, meaning that we can have $\alpha$ such that the likelihood in $\Gamma_\alpha$ concentrating to $\gamma_\alpha^*$ arbitrarily slow. This shows the importance of condition 7).
\end{remark}

\begin{lemma}[Uniform Convergence of $D^2 L_K$ in \eqref{eqn:L''}]\label{lemma:L''} Assume conditions 1)--5). Then \eqref{eqn:L''} holds for any $\gamma_{\alpha,n}$ with $ \sup_{\alpha \in \mathcal{A}_\epsilon^+} \Vert \gamma_{\alpha,n} - \gamma_\alpha^* \Vert \rightarrow 0$ as $n \rightarrow \infty$.
\end{lemma}

\begin{proof}
Fix such a family of sequence $\gamma_{\alpha,n}$. Follow the steps 
in \cite{Bickel:1998}, we extend $\{(X_i, Y_i)\}_{i=1}^\infty$ to $\{(X_i, Y_i)\}_{i=-\infty}^\infty$, and let $p_K({\bm Y}_n|{\bm Y}_{-m:(n-1)}, \boldsymbol{\phi}_K)$ be the conditional likelihood of ${\bm Y}_n$ given ${\bm Y}_{-m:(n-1)}$ and parameter $\boldsymbol{\phi}_K$. 
By \cite{Bickel:1998}, Lemma 10, under conditions 1)-5), for $\alpha \in \mathcal{A}$, there exists an ergodic sequence $\zeta_{\alpha,n}$ such that
\begin{equation*}
    D_\gamma^2 \log p_K({\bm Y}_n|{\bm Y}_{-m:(n-1)}, \varphi_\alpha(\gamma_\alpha^*)) \xrightarrow[m \rightarrow \infty]{L^1}  \zeta_{\alpha,n}.
\end{equation*}
The proof of \cite{Bickel:1998}, Lemma 2, shows that, for each $\alpha \in \mathcal{A}$,
\begin{equation}\label{eqn:LemmaD''-1}
\bigg\Vert \frac{D_\gamma^2 L_K({\bm y}_{1:n}|\varphi_\alpha(\gamma_{\alpha,n}))}{n}
-
\frac{\sum_{i=1}^n \zeta_{\alpha,i}
}{n} \bigg\Vert
\xrightarrow[n \rightarrow \infty]{P^*} 0,
\end{equation}
and
\begin{equation}\label{eqn:LemmaD''-2}
\frac{\sum_{i=1}^n \zeta_{\alpha,i}
}{n} \xrightarrow[n \rightarrow \infty]{} J_\alpha
\end{equation}
with probability one for some non-singular matrix $J_\alpha$.

To further extend \eqref{eqn:LemmaD''-2} to an uniform version, we first show that 
\begin{equation}\label{eqn:LemmaD''-3}
E^* \left[ \sup_{\alpha \in \mathcal{A}_\epsilon^+} \Vert \zeta_{\alpha,n}\Vert \right] < \infty
\end{equation}
for each $n$. As discussed in the proof of Lemma \ref{Lemma:iv}, $W$ is bounded away from zero on $\mathcal{A}_\epsilon^+$, so that $\cup_{\alpha \in \mathcal{A}_\epsilon^+} \Phi_{K,\alpha}$ is bounded away from $\partial \Phi_K$. Hence, by conditions 4) and 5), for any $\alpha' \in \mathcal{A}_\epsilon^+$, there exists an open neighborhood $\mathcal{O}(\alpha')$ such that
\begin{equation*}
    E^*\left[ \sup_{\overline{\alpha} \in \mathcal{O}(\alpha')} \Vert \zeta_{\overline{\alpha},n}\Vert \right]
    \leq 
    E^*\left[ \Vert \zeta_{\alpha',n} \Vert \right] + \epsilon.
\end{equation*}
Since $\{\mathcal{O}(\alpha'): \alpha' \in \mathcal{A}_\epsilon^+\}$ forms an open covering of $\mathcal{A}_\epsilon^+$, which is compact (as discussed in the proof of Lemma \ref{Lemma:L}), we can choose finitely many $\{ \alpha_m, 1 \leq m \leq M \}$ so that  $\{\mathcal{O}(\alpha_m), m =1,\cdots,M\}$ forms an open covering of $\mathcal{A}_\epsilon^+$. As such, (\ref{eqn:LemmaD''-3}) holds since
\begin{align*}
    E^*\left[ \sup_{\alpha \in \mathcal{A}_\epsilon^+} \Vert \zeta_{\alpha', n}) \Vert \right]
    & \leq 
    \max_{1 \leq m \leq M} E^*\left[ \sup_{\overline{\alpha} \in \mathcal{O}(\alpha_m)} \Vert \zeta_{\overline{\alpha},n}\Vert \right]
\\
    & \leq  \max_{1 \leq m \leq M}
    E^*\left[ \Vert \zeta_{\alpha_m,n} \Vert \right] + \epsilon
    < \infty.
\end{align*}

We then show that 
\begin{equation}\label{eqn:LemmaD''-4}
\sup_{\alpha \in \mathcal{A}_\epsilon^+} \bigg\Vert \frac{\sum_{i=1}^n \zeta_{\alpha,i}
}{n} - J_\alpha \bigg\Vert \xrightarrow[n \rightarrow \infty]{P^*} 0.
\end{equation}
Notice that:
\begin{itemize}
\item[a)] $\mathcal{A}_\epsilon^+$ is compact;
\item[b)] by condition 5), $\zeta_{\alpha,1}$ is continuous in $\alpha$;
\item[c)] for each $\alpha$, by \eqref{eqn:LemmaD''-2}, $\frac{1}{n} \sum_{i=1}^n \zeta_{\alpha,i} \xrightarrow{n \rightarrow \infty} E^*\zeta_{\alpha,1}$ with probability one;
\item[d)] by \eqref{eqn:LemmaD''-3},  $\sup_{\alpha \in \mathcal{A}_\epsilon^+} \Vert \zeta_{\alpha,n} \Vert \in L^1$.
\end{itemize}
Therefore, by the uniform law of large numbers provided in Lemma \ref{ULLN}, (\ref{eqn:LemmaD''-4}) holds.

Finally, by using a similar open covering argument as proving \eqref{eqn:LemmaD''-3}, we can show that 
\begin{equation}\label{eqn:LemmaD''-5}
\sup_{\alpha \in \mathcal{A}_\epsilon^+}\bigg\Vert \frac{D_\gamma^2 L_K({\bm y}_{1:n}|\varphi_\alpha(\gamma_{\alpha,n}))}{n}
-
\frac{\sum_{i=1}^n \zeta_{\alpha,i}
}{n} \bigg\Vert
\xrightarrow[n \rightarrow \infty]{P^*} 0.
\end{equation}
The proof is completed by combining (\ref{eqn:LemmaD''-4}) and (\ref{eqn:LemmaD''-5}).
\end{proof}

\begin{remark}\label{remark-Lemma:L''}
The above argument uses the fact that $\cup_{\alpha \in \mathcal{A}_\epsilon^+} \Phi_{K,\alpha}$ is bounded away from $\partial \Phi_K$, which ensures an uniform ergodicity across all $\boldsymbol{\phi}_\alpha^*$ for all $\alpha \in \mathcal{A}_\epsilon^+$. If we replace $\mathcal{A}_\epsilon^+$ to $\mathcal{A}^+$, then $\Phi_{K,\alpha}$ can be arbitrarily close to $\partial \Phi_K$, meaning that $\boldsymbol{\phi}_\alpha^*$ can be arbitrarily close to non-ergodic. This is another reason why the condition 7) is important. See \cite{Gassiat2014} for related discussions. The same issues happen in Lemma \ref{lemma:L'}.
\end{remark}

\begin{lemma}[Uniform Convergence of $D_\gamma L_K$ in \eqref{eqn:L'}]\label{lemma:L'} Assume conditions 1)--6). Then \eqref{eqn:L'} holds.
\end{lemma}

\begin{proof}
By the proof of \cite{Bickel:1998}, Lemma 1, under conditions 1)-5), for any $\alpha \in \mathcal{A}$, there exists a sequence $\{ \eta_{\alpha,i}, i \geq 1\}$ such that
\begin{equation}\label{eqn:LemmaD'-1}
\Bigg\Vert
\frac{D_\gamma L_K({\bm y}_{1:n}|\varphi_\alpha(\gamma_\alpha^*))}{\sqrt{n}}
-
\frac{\sum_{i=1}^n \eta_{\alpha,i}
}{\sqrt{n}}
\Bigg\Vert \xrightarrow[n \rightarrow \infty]{P^*} 0,
\end{equation}
and
\begin{equation*}
\lim_{n \rightarrow \infty} \frac{\sum_{i=1}^n  \eta_{\alpha,i}}{\sqrt{n}} \xrightarrow[n \rightarrow \infty]{d} \mathcal{N}(\vec{0}, J_\alpha).
\end{equation*}
Define
\begin{equation*}
    G_{s:t} = \sup_{\alpha \in \mathcal{A}_\epsilon^+} \sum_{i=s+1}^t \eta_{\alpha,i} J_\alpha^{-1/2}.
\end{equation*}
Suppose we can prove that, as $n \rightarrow \infty$,
\begin{equation}\label{eqn:LemmaD'-2}
    \lim_{n \rightarrow \infty} \sup_{\alpha \in \mathcal{A}_\epsilon^+} \frac{\sum_{i=1}^n  \eta_{\alpha,i} J_\alpha^{-1/2}}{\sqrt{n}} = 
    \lim_{n \rightarrow \infty}  \frac{G_{0:n}}{\sqrt{n}} =   O_{P^*}(1),
\end{equation}
then by (\ref{eqn:LemmaD'-1}), (\ref{eqn:LemmaD'-2}) and a similar open covering argument as in the proof of Lemma \ref{lemma:L''}, we prove Lemma \ref{lemma:L'}.

To prove (\ref{eqn:LemmaD'-2}), we use the central limit theorem of subadditive processes provided by \cite{Ishitani1977}, Theorem 1. To check the conditions, first note that
\begin{align*}
    G_{s:t} = & \sup_{\alpha \in \mathcal{A}_\epsilon^+} \sum_{i=s+1}^t \eta_{\alpha,i} J_\alpha^{-1/2} 
    \leq \sup_{\alpha \in \mathcal{A}_\epsilon^+} \sum_{i=s+1}^u \eta_{\alpha,i} J_\alpha^{-1/2}  +
\sup_{\alpha \in \mathcal{A}_\epsilon^+} \sum_{i=u+1}^t \eta_{\alpha,i} J_\alpha^{-1/2}
= G_{s:u} + G_{u:t},
\end{align*}
so $G_{0:n}$ is subadditive. The condition (2) in \cite{Ishitani1977}, Theorem 1, holds due to Lemma \ref{MixtureCondition}. In addition, note that the proof of \cite{Bickel:1998}, Lemma 1, shows that $\eta_{\alpha,n}$ is a stationary, ergodic martingale increment sequence, so conditions (1), (3) and (4) hold. To sum up, all conditions in \cite{Ishitani1977}, Theorem 1, are satisfied, so \eqref{eqn:LemmaD'-2} holds. The proof is completed.
\end{proof}

\begin{lemma}[Uniform Law of Large Numbers used in Lemma \ref{lemma:L''}]\label{ULLN}
Suppose, for each $\alpha \in \overline{\mathcal{A}}$, $\zeta_{\alpha,i}$ is a sequence of random variables satisfying
\begin{itemize}
    \item[(a)] $\overline{\mathcal{A}}$ is compact;
    \item[(b)] $\zeta_{\alpha,1}$ is continuous in $\alpha$ at each $\alpha \in \overline{\mathcal{A}}$ with probability one;
    \item[(c)] for each $\alpha$, $\frac{1}{n} \sum_{i=1}^n \zeta_{\alpha,i} \xrightarrow{n \rightarrow \infty} E^*\zeta_{\alpha,1}$ with probability one;
    \item[(d)] $E^* \left[ \sup_{\alpha \in \overline{\mathcal{A}}} \Vert \zeta_{\alpha,1} \Vert \right]< \infty$.
\end{itemize}
Then, we have
\begin{equation*}
    \sup_{\alpha \in \overline{\mathcal{A}}} \Bigg\Vert \frac{1}{n}\sum_{i=1}^n \zeta_{\alpha,i} - E^* \zeta_{\alpha,1} \Bigg\Vert \xrightarrow[n \rightarrow \infty]{P^*} 0.
\end{equation*}
\end{lemma}

\begin{proof}
We will only prove the case with one dimensional $\zeta_{\alpha,i}$, as the higher dimensional case can be proved similarly. Throughout this proof, for each $\alpha \in \overline{\mathcal{A}}$ and $\delta > 0$, let $\mathcal{B}_\delta(\alpha) = \left\{ \alpha': \Vert \alpha' - \alpha \Vert < \delta \right\}$.

The proof is basically based on the proof of \cite{jennrich1969asymptotic}, Theorem 2, with minor adjustments. First, for any $\alpha \in \overline{\mathcal{A}}$, define 
$$\Delta_1^\alpha(\delta) = \sup_{\alpha' \in \mathcal{B}_\delta(\alpha)} \zeta_{\alpha',1} - \inf_{\alpha' \in \mathcal{B}_\delta(\alpha)} \zeta_{\alpha',1}.$$
Note that for each $\alpha \in \overline{\mathcal{A}}$, (i) by condition (b), we have $\Delta_1^\alpha (\delta)
\downarrow 0$ with probability one as $\delta \downarrow 0$; (ii) $\Delta_1^\alpha (\delta) \leq 2 \sup_{\alpha \in \overline{\mathcal{A}}} |\zeta_{\alpha,1}|$; (iii) by condition (d), $E^*\left[ \sup_{\alpha \in \overline{\mathcal{A}}} |\zeta_{\alpha,1}| \right]< \infty$. Therefore by (i)-(iii) and the dominated convergence theorem, we have $E^*\left[ \Delta_1^\alpha (\delta) \right] \downarrow 0$ as $\delta \downarrow 0$. Hence, for any $\epsilon > 0$, for all $\alpha \in \overline{\mathcal{A}}$, the exists $\delta_\epsilon(\alpha)>0$ such that $E^*\left[ \Delta_1^\alpha \left( \delta_\epsilon(\alpha) \right) \right] < \epsilon$.

Now, by (a), $\overline{\mathcal{A}}$ is compact, and is covered by $\{\mathcal{B}_{\delta_\epsilon(\alpha)}(\alpha): \alpha \in \overline{\mathcal{A}}\}$, so we can choose finitely many $\alpha_k \in \overline{\mathcal{A}}$ ($k = 1, \cdots, K$) such that $\overline{\mathcal{A}} \subset \bigcup_{k=1}^K \mathcal{B}_{\delta_\epsilon(\alpha_k)}(\alpha_k)$. Hence, we have 
\begin{align}
\notag
&
    \sup_{\alpha \in \overline{\mathcal{A}}} \left\{ \frac{1}{n}\sum_{i=1}^n \zeta_{\alpha,i} - E^*\zeta_{\alpha,1} \right\} 
= 
    \max_k \sup_{\alpha \in \mathcal{B}_{\delta_\epsilon(\alpha_k)}(\alpha_k)}
    \left\{ \frac{1}{n}\sum_{i=1}^n \zeta_{\alpha,i} - E^*\zeta_{\alpha,1} \right\}  \\
\label{ULLN-1}
\leq &
    \max_k 
    \left\{ \frac{1}{n}\sum_{i=1}^n \sup_{\alpha \in \mathcal{B}_{\delta_\epsilon(\alpha_k)}(\alpha_k)} \zeta_{\alpha,i} -  E^*\left[ \inf_{\alpha \in \mathcal{B}_{\delta_\epsilon(\alpha_k)}(\alpha_k)} \zeta_{\alpha,1} \right] \right\}.
\end{align}
Since $\zeta_{\alpha,i}$ are ergodic (as stated in (c)), so are $\sup_{\alpha \in \mathcal{B}_{\delta_\epsilon(\alpha_k)}(\alpha_k)} \zeta_{\alpha,i}$. Hence, by (d) and the weak law of large numbers of ergodic sequence, as $n \rightarrow \infty$,
\begin{equation}\label{ULLN-2}
\frac{1}{n}\sum_{i=1}^n \sup_{\alpha \in \mathcal{B}_{\delta_\epsilon(\alpha_k)}(\alpha_k)} \zeta_{\alpha,i} = E^*\left[
\sup_{\alpha \in \mathcal{B}_{\delta_\epsilon(\alpha_k)}(\alpha_k)} \zeta_{\alpha,1} \right] + o_P(1).
\end{equation}
Combining \eqref{ULLN-1} and \eqref{ULLN-2}, we have
\begin{align*}
    \sup_{\alpha \in \overline{\mathcal{A}}} \left[ \frac{1}{n}\sum_{i=1}^n \zeta_{\alpha,i} - E^*\zeta_{\alpha,1} \right] 
\leq &
    \max_k 
    \left\{ E^*\left[ \sup_{\alpha \in \mathcal{B}_{\delta_\epsilon(\alpha_k)}(\alpha_k)} \zeta_{\alpha,1} \right] -  E^*\left[ \inf_{\alpha \in \mathcal{B}_{\delta_\epsilon(\alpha_k)}(\alpha_k)} \zeta_{\alpha,1} \right] \right\} 
    + o_P(1) \\
= & \max_k E^*[\Delta_1^\alpha(\delta_\epsilon)] + o_P(1)
\leq \epsilon + o_P(1).
\end{align*}
A similar argument shows that
\begin{equation*}
    \inf_{\alpha \in \overline{\mathcal{A}}} \left[ \frac{1}{n}\sum_{i=1}^n \zeta_{\alpha,i} - E^*\zeta_{\alpha,1} \right] \geq - \epsilon + o_{P^*}(1).
\end{equation*}
Since $\epsilon$ can be arbitrarily small, the proof is completed.
\end{proof}

\begin{lemma}[Mixing Condition used in Lemma \ref{lemma:L'}]\label{MixtureCondition} Recall $P^*$ is the probability measure of $\{(X_i, Y_i), i \geq 1\}$ under the true parameter $\boldsymbol{\phi}^*$. Suppose that under $P^*$, $\{ X_i, i \geq 1\}$ is ergodic (irreducible, aperiodic, and positive recurrent.) Then there exist constant $C>0$ and $\rho \in (0,1)$ such that, 
\begin{align*}
    \Big\vert P_\nu^* \left\{ Y_{1:n} \in A \right\} P_\nu^* \left\{Y_{n+m:\infty} \in B \right\} - P_\nu^* \left\{ Y_{1:n} \in A , Y_{n+m:\infty} \in B \right\} \Big\vert 
    \leq C \rho^m \times P_\nu^* \{ Y_{1:n} \in A \},
\end{align*}
for any $A$, $B$, $m \geq 0$ and initial measure $\nu = (\nu_1, \cdots, \nu_{K^*})$. Here, 
\begin{equation*}
    P_\nu^*\{ \cdot \} = \sum_{k=1}^{K^*} \nu_k P^*\{ \cdot | X_0 = k \}. 
\end{equation*}
\end{lemma}

\begin{proof}[Proof of Lemma \ref{MixtureCondition}]
Let $\mu = (\mu_1, \cdots, \mu_{K^*})$ be the invariant probability of $X_n$ (under true parameter $\boldsymbol{\phi}^*$.) Since $X_n$ is ergodic under $P^*$, there exist $R>0$ and $\rho \in (0,1)$ such that for any $n, m,$
\begin{equation}
    \max_{1 \leq k,l \leq K^*} \Big\vert P^*\left\{ X_{n+m} = l \big\vert X_n = k\right\} - \mu_l \Big\vert \leq R \rho^m.
\end{equation}
As such, we have
\begin{align*}
& 
P_\nu^* \left\{ Y_{1:n} \in A, Y_{n+m:\infty} \in B \right\}  \\
= & \sum_{1 \leq k,l \leq K^*} P_\nu^* \left\{ Y_{1:n} \in A, X_n=k, X_{n+m}=l, Y_{n+m:\infty} \in B \right\} \\
= &
\sum_{1 \leq k,l \leq K^*} P_\nu^* \left\{ Y_{1:n} \in A, X_n=k \right\} \times  P^* \left\{ X_{n+m}=l | X_n = k\right\} \times P^* \left\{ Y_{n+m:\infty} \in B | X_{n+m} = l \right\} \\
\leq &
\sum_{1 \leq k,l \leq K^*} P_\nu^* \left\{ Y_{1:n} \in A, X_n=k \right\} \times (\mu_l + R \rho^m) \times P^* \left\{ Y_{n+m:\infty} \in B | X_{n+m} = l \right\} \\
= &
\sum_{1 \leq k,l \leq K^*} P_\nu^* \left\{ Y_{1:n} \in A, X_n=k \right\} \times \mu_l \times P^* \left\{ Y_{n+m:\infty} \in B | X_{n+m} = l \right\}
\\
& +
\sum_{1 \leq k,l \leq K^*} P_\nu^* \left\{ Y_{1:n} \in A, X_n=k \right\} \times R \rho^m \times P^* \left\{ Y_{n+m:\infty} \in B | X_{n+m} = l \right\} \\
\leq &
\sum_{1 \leq k \leq K^*} P_\nu^* \left\{ Y_{1:n} \in A, X_n=k \right\} \times \sum_{1 \leq l \leq K^*} \mu_l P^* \left\{ Y_{n+m:\infty} \in B | X_{n+m} = l \right\} \\
&
+ \sum_{1 \leq k,l \leq K^*} P_\nu^* \left\{ Y_{1:n} \in A, X_n=k \right\} \times R \rho^m \\
= &
P_\nu^* \{ Y_{1:n} \in A \} \times P_\mu^* \{ Y_{n+m:\infty} \in B\} + K^* R \rho^m \times P_\nu^* \{ Y_{1:n} \in A \} .
\end{align*}
A similar argument leads to the lower bound, and so we have
\begin{align}\label{Part1}
\notag
    & \Big\vert P_\nu^* \left\{ Y_{1:n} \in A, Y_{n+m:\infty} \in B \right\} - P_\nu^* \left\{ Y_{1:n} \in A \right\} P_\mu^* \left\{ Y_{n+m:\infty} \in B \right\} \Big\vert \\
    \leq &
    K^* R \rho^m \times P_\nu^* \{ Y_{1:n} \in A \}.
\end{align}
On the other hand, we have 
\begin{align*}
    & P_\nu^*\left\{ Y_{1:n} \in A \right\} P_\nu^*\left\{ Y_{n+m:\infty} \in B \right\} \\
=  & P_\nu^*\left\{ Y_{1:n} \in A \right\} \sum_{1 \leq l \leq K^*} P_\nu^*\{ X_{n+m}=l\} P^*\left\{ Y_{n+m:\infty} \in B \big\vert X_{n+m}=l \right\} \\
\leq &
P_\nu^*\left\{ Y_{1:n} \in A \right\} \sum_{1 \leq l \leq K^*} (\mu_l + R \rho^{n+m}) P^*\left\{ Y_{n+m:\infty} \in B \big\vert X_{n+m}=l \right\}
\\
= &
P_\nu^*\left\{ Y_{1:n} \in A \right\} P_\mu^*\left\{ Y_{n+m:\infty} \in B \right\} \\
& +
P_\nu^*\left\{ Y_{1:n} \in A \right\} \sum_{1 \leq l \leq K^*} R \rho^{n+m} P^*\left\{ Y_{n+m:\infty} \in B \big\vert X_{n+m}=l \right\} \\
\leq &
P_\nu^*\left\{ Y_{1:n} \in A \right\} P_\mu^*\left\{ Y_{n+m:\infty} \in B \right\} + K^* R \rho^{n+m} \times P_\nu^* \{ Y_{1:n} \in A \} \\
\leq &
P_\nu^*\left\{ Y_{1:n} \in A \right\} P_\mu^*\left\{ Y_{n+m:\infty} \in B \right\} + K^* R \rho^m \times P_\nu^* \{ Y_{1:n} \in A \}.
\end{align*}
A similar argument leads to the lower bound, and so we have
\begin{align}\label{Part2}
\notag
& \Big\vert P_\nu^* \left\{ Y_{1:n} \in A \right\} P_\nu^* \left\{Y_{n+m:\infty} \in B \right\} - P_\nu^* \left\{ Y_{1:n} \in A \right\} P_\mu^* \left\{ Y_{n+m:\infty} \in B \right\} \Big\vert \\
\leq & K^* R \rho^m \times P_\nu^* \{ Y_{1:n} \in A \}.
\end{align}
Combining \eqref{Part1} and \eqref{Part2}, we complete the proof for Lemma \ref{MixtureCondition}.
\end{proof}

\subsection{ Proof of Corollary \ref{theorem:asymptotic_consistency_gm}}
\label{appendix:proofofconsistencygm}

Recall that $\boldsymbol{\phi}^*$ is the true parameter, $\mu(\boldsymbol{\phi}^*) = \left( \mu_1(\boldsymbol{\phi}^*), \cdots, \mu_{K^*}(\boldsymbol{\phi}^*) \right)$ is the corresponding invariant probability measure, and $P^*$ is the probability law of $\left\{ (X_i, Y_i), i \geq 1 \right\}$ under $\boldsymbol{\phi}^*$ and initial distribution $\mu(\boldsymbol{\phi}^*)$. In addition, for any $K$, recall
\begin{equation*}
\mathcal{Q}_K^{mix} = \{Q_K: q_{1k}=q_{2k}=\cdots=q_{Kk} \mbox{ for all } 1 \leq k \leq K\},
\end{equation*}
and define $\Phi_K^{mix} = \left\{ \boldsymbol{\phi}_K = (Q_K; \boldsymbol{\theta}_1, \cdots, \boldsymbol{\theta}_K): Q_K \in \mathcal{Q}_K^{mix} \right\}$. Note that for any $\boldsymbol{\phi}_K \in \Phi_K^{mix}$ with $q_{ik} = q_k$ for all $1 \leq k \leq K$, we have
\begin{equation*}
    p({\bm y}_{1:n}|\boldsymbol{\phi}_K) = \prod_{i=1}^n \sum_{x_i=1}^K q_{x_i} f(y_i|\boldsymbol{\theta}_{x_i}),
\end{equation*}
and hence, the log-likelihood
\begin{equation*}
    L({\bm y}_{1:n}|\boldsymbol{\phi}_K) = \sum_{i=1}^n \log\left( \sum_{x_i=1}^K q_{x_i} f(y_i|\boldsymbol{\theta}_{x_i}) \right),
\end{equation*}
which has an additive form. This ensures the asymptotic behavior of $L$ as well as its derivatives at any $\boldsymbol{\phi}_K \in \Phi_K^{mix}$ under $P^*$.

Recall that we use the following three steps to prove \eqref{eqn_thm_consist}:
\begin{enumerate}[(a)]
    \item For $K^*$, the posterior distribution over $\Phi_{K^*}$ is asymptotically normal, centering at $\boldsymbol{\phi}^*$;
    \item For $K < K^*$, the log-likelihood decays exponentially due to the fact that the entire $\Phi_K$ can be viewed as a subset of $\Phi_{K^*}$, bounded away from $\boldsymbol{\phi}^*$;
    \item For $K > K^*$, we decompose $\Phi_K$ into a family of $\Phi_{K,\alpha}$ so that the posterior distribution over $\Phi_{K,\alpha}$ is asymptotically normal, centering at the unique ``true'' value $\boldsymbol{\phi}_\alpha^*$ in $\Phi_{K,\alpha}$.
\end{enumerate}
Note that a reason for (a)-(c) hold is $\boldsymbol{\phi}^*$ belongs to $\Phi_{K^*}$. However, in the case of Corollary \ref{theorem:asymptotic_consistency_gm}, $\boldsymbol{\phi}^*$ is generally not in $\Phi_{K^*}^{mix}$, so there is no way to have a posterior distribution concentrated at $\boldsymbol{\phi}^*$.

To overcome this difficulty, we follow the argument in \cite{bunke1998asymptotic} to  have that the posterior distribution being asymptotically normal centering at some ``pseudo-true'' value in $\Phi_{K^*}^{mix}$ (note that the log-likelihood follows an additive structure in this case.) Hence, by replacing all ``true'' values in the original proof of Theorem \ref{theorem:asymptotic_consistency_hmm} with the corresponding ``pseudo-true'' values, we get the proof. Details are as follows.

\begin{proof}[Part (a) under Corollary \ref{theorem:asymptotic_consistency_gm}]
Define 
\begin{equation*}
    Q_{mix}^* = 
    \begin{pmatrix}
    \mu_1(\boldsymbol{\phi}^*) & \mu_2(\boldsymbol{\phi}^*) & \cdots & \mu_{K^*} (\boldsymbol{\phi}^*) \\
    \mu_1(\boldsymbol{\phi}^*) & \mu_2(\boldsymbol{\phi}^*) & \cdots & \mu_{K^*} (\boldsymbol{\phi}^*) \\
    \vdots & \vdots & \ddots & \vdots \\
    \mu_1(\boldsymbol{\phi}^*) & \mu_2(\boldsymbol{\phi}^*) & \cdots & \mu_{K^*} (\boldsymbol{\phi}^*)
    \end{pmatrix},
\end{equation*}
and set $\boldsymbol{\phi}_{mix}^* = (Q_{mix}^*; \boldsymbol{\theta}_1^*, \cdots, \boldsymbol{\theta}_K^*)$. A direct check shows that $\boldsymbol{\phi}_{mix}^*$ is  ``pseudo-true'' value in the sense that it is the unique point
in $\Phi_{K^*}^{mix}$ that obtains the minimum Kullback-Leibler divergence to $\boldsymbol{\phi}^*$. Hence, by the argument in \cite{bunke1998asymptotic} (with all results within related to i.i.d. data being replaced by the corresponding HMM version provided by \cite{Bickel:1998}, \cite{HMMCLT} and \cite{LerouxM:1992}), the posterior distribution over $\Phi_{K^*}^{mix}$ is asymptotically normal, centering at $\boldsymbol{\phi}_{mix}^*$, that is,
\begin{equation}\label{eqn:CLT-K*-mix}
    \frac{p_{K^*}({\bm y}_{1:n}|\boldsymbol{\phi}_{mix}^*) p_0(\boldsymbol{\phi}_{mix}^*)}{n^{\Delta_{K^*}^{mix}/2} p_{K^*}({\bm y}_{1:n})} = O_{P^*}(1)
\end{equation}
as $n \rightarrow \infty$, where 
\begin{equation}\label{Delta-mix}
\Delta_K^{mix} := Kd + (K-1)
\end{equation}
is the dimension of $\Phi_{K}^{mix}$.
\end{proof}

\begin{proof}[Part (b) under Corollary \ref{theorem:asymptotic_consistency_gm}] Note that when $K < K^*$, the entire $\Phi_K^{mix}$ is a subset of $\Phi_{K^*}^{mix}$, and is bounded away from $\boldsymbol{\phi}_{mix}^*$. As such, we have
\begin{equation*}
    \inf_{\boldsymbol{\phi}_K \in \Phi_K^{mix}} KL(\boldsymbol{\phi}_K, \boldsymbol{\phi}^*) > KL(\boldsymbol{\phi}_{mix}^*, \boldsymbol{\phi}^*)
\end{equation*}
with $KL(\cdot, \cdot)$ denotes the Kullback-Leibler divergence.
As such, a similar argument as in Section \ref{Proof-under} gives
\begin{equation}\label{eqn:CLT-under-mix}
    \frac{p_{K}({\bm y}_{1:n})}{p_{K^*}({\bm y}_{1:n}|\boldsymbol{\phi}_{mix}^*)} = O_{P^*} \left( e^{-cn} \right)
\end{equation}
as $n \rightarrow \infty$ for some $c>0$.
\end{proof}

\begin{proof}[Part (c) under Corollary \ref{theorem:asymptotic_consistency_gm}] In Theorem \ref{theorem:asymptotic_consistency_hmm}, when $K>K^*$, we decompose $\Phi_K$ into $\cup_{\alpha} \Phi_{K,\alpha}$ so that each $\Phi_{K,\alpha}$ contains at most one ``true'' value. Similarly, here we will need to decompose $\Phi_K^{mix}$ into $\cup_\alpha \Phi_{K,\alpha}^{mix}$ so that each $\Phi_{K,\alpha}^{mix}$ contains at most one ``pseudo-true'' value. To do so, set $\alpha = (S,W)$ with $S=(S_1, \cdots, S_{K^*})$ being a partition of $\{1, 2, \cdots, K\}$, and $W = (W_1, \cdots, W_K) \in (0,1)^K$. Further recall the decomposition $\Theta = \cup_{k=1}^{K^*} \Theta_k$ defined before \eqref{Theta-unique} so that $\boldsymbol{\theta}_k^*$ is an interior point for $\Theta_k$ for all $k = 1, 2, \cdots, K^*$, and the function $s(k)$ defined in \eqref{s(k)}. Then, similar to Example \ref{exmp-mixture}, for any $\alpha = (S,W)$, we set
\begin{align*}
    \Phi_{K,\alpha}^{mix} = & \Bigg\{ \boldsymbol{\phi}_K \in \Phi_K^{mix}: \mbox{ for all } \ell \in \{1, 2, \cdots, K^*\},  \boldsymbol{\theta}_k \in \Theta_\ell \mbox{ for all } k \in S_\ell \\
    & \mbox{ and for each } i,j  \in \{1, 2, \cdots, K\},  \frac{q_{ij}}{\sum_{k \in S_{s(i)}}q_{ik}} = W_j \Bigg\}.
\end{align*}
Then, a similar argument as in Section \ref{Proof-over} shows that $\Phi_{K,\alpha}^{mix}$ forms a partition of $\Phi_K^{mix}$, each with dimension $Kd + (K^*-1) = \Delta_{K^*}^{mix} + (K-K^*)d$, and contains at most one ``pseudo-true'' value $\boldsymbol{\phi}_{mix,\alpha}^* := (Q_{mix,\alpha}^*; \boldsymbol{\theta}_{s(1)}, \cdots, \boldsymbol{\theta}_{s(K)})$ with
\begin{equation*}
    Q_{mix,\alpha}^* := 
    \begin{pmatrix}
    \mu_{s(1)}(\boldsymbol{\phi}^*)W_1 & \mu_{s(2)}(\boldsymbol{\phi}^*)W_2 & \cdots & \mu_{s(K)} (\boldsymbol{\phi}^*)W_K \\
    \mu_{s(1)}(\boldsymbol{\phi}^*)W_1 & \mu_{s(2)}(\boldsymbol{\phi}^*)W_2 & \cdots & \mu_{s(K)} (\boldsymbol{\phi}^*)W_K \\
    \vdots & \vdots & \ddots & \vdots \\
    \mu_{s(1)}(\boldsymbol{\phi}^*)W_1 & \mu_{s(2)}(\boldsymbol{\phi}^*)W_2 & \cdots & \mu_{s(K)} (\boldsymbol{\phi}^*)W_K
    \end{pmatrix},
\end{equation*}
in the sense that $\boldsymbol{\phi}_{mix,\alpha}^*$ obtains the minimum Kullback-Leibler divergence to $\boldsymbol{\phi}^*$ on $\Phi_{K,\alpha}^{mix}$. Hence, similar to part (a), by \cite{bunke1998asymptotic} (with all necessary HMM replacement), we have the asymptotic normality of the posterior distribution over $\Phi_{K,\alpha}^{mix}$. Thus, through a similar procedure in Section \ref{Proof-over}, we have, as $n \rightarrow \infty$,
\begin{equation}\label{eqn:CLT-over-mix}
    n^{(\Delta_{K^*}^{mix} + (K-K^*)d)/2} \frac{p_K({\bm y}_{1:n})}{p_{K^*}({\bm y}_{1:n}|\boldsymbol{\phi}_{mix}^*) } = O_{P^*}(1)
\end{equation}
\end{proof}

\begin{proof}[Proof of Corollary \ref{theorem:asymptotic_consistency_gm}]
For $K < K^*$, combining \eqref{eqn:CLT-K*-mix} and \eqref{eqn:CLT-under-mix}, we have \eqref{eqn_thm_consist_under} holds. As for $K > K^*$, combining \eqref{eqn:CLT-K*-mix} and \eqref{eqn:CLT-over-mix}, we have \eqref{eqn_thm_consist} holds. Hence, all results in Theorem \ref{theorem:asymptotic_consistency_hmm} hold under the scenario in Corollary \ref{theorem:asymptotic_consistency_gm}, which completes the proof.
\end{proof}

\section{ Simulation Studies for Estimation of Normalizing Constant}
\label{appendix:simulations_normalizing_constant}

We use models with known normalizing constants to test the performance of our estimating of normalizing constants proposed in Section~\ref{subsubsec:EstimationProcedure}. The first family of models is $d$-dimensional Gaussian mixture models with three distinct components whose covariance matrices are diagonal with diagonal elements all equal to $0.1$, $2\leq d\leq 30$; the normalizing constant is set to be $C_1 = \exp(10)$. The second family of models, with normalizing constant $C_2 = \exp(2)$, has three independent dimensions: the first dimension is Gaussian with mean $1$ and variance $1$, the second dimension is a student-t distribution with a degree of freedom $2$, and the third dimension is Gamma distribution with shape parameter $6$ and scale parameter $2$.

We perform repeated simulations on the two families of models as follows: first, simulate $N_{sim}$ samples independently from the model, apply the importance sampling algorithm mentioned in Section~\ref{subsubsec:EstimationProcedure} with $N_{is}$ samples from the fitted importance function, with the Gaussian tail and the t-tail. For comparison, we also use the reciprocal importance sampling with fitted Gaussian mixture and multivariate t-distribution mixture as the importance functions, respectively. The results are summarized in Table~\ref{table:simulation_results_normalizing_const}. \textcolor{blue}{It shows that the importance sampling method, despite the more computational cost, gives better estimates of the normalizing constant as compared to the reciprocal importance sampling, regardless of the dimension of the parameter space and the complexity of the posterior, e.g. heavy tail or multi-mode. The superior performance of the importance sampling that we see here is more due to our adaptive choice of the importance function based on posterior samples, such that the importance function covers the ``regions of interest'' very well. }

\begin{table*}[tbph]
\caption{Simulation results of the algorithm in Section~\ref{subsubsec:EstimationProcedure}}
\begin{tabular}{cccccccc}
\hline
\hline
$\mathcal{M}$ & $D$ & $N_{sim}$ & $N_{is}$ & $CI_{1}$ & $CI_{2}$ & $CI_{3}$ & $CI_{4}$\\
\hline
2 & 3 & 2,000 & 4,000 & [-0.042, 0.023] & [-0.035, 0.028] & [-0.378, -0.237] & [0.091, 0.265] \\
2 & 3 & 10,000 & 10,000 & [-0.014, 0.013] & [-0.017, 0.016] & [-0.366, -0.282] & [0.112, 0.210] \\ 
\hline 
1 & 4 & 2,000 & 4,000 & [-0.055, 0.042] & [-0.069, 0.043] & [-0.056, 0.035] & [0.548, 0.643] \\
1 & 6 & 2,000 & 4,000 & [-0.063, 0.046] & [-0.070, 0.046] & [-0.076, 0.014] & [0.645, 0.738]\\
1 & 8 & 2,000 & 4,000 & [-0.073, 0.021] & [-0.076, 0.026] & [-0.106, -0.011] & [0.672, 0.791] \\
1 & 10 & 2,000 & 4,000 & [-0.075, 0.018] & [-0.087, 0.045] & [-0.108, -0.023] & [0.659, 0.833] \\
\hline
1 & 10 & 10,000 & 10,000 & [-0.027, 0.020] & [-0.035, 0.023] & [-0.030, 0.001] & [0.774, 0.847]\\
1 & 15 & 10,000 & 10,000 & [-0.039, 0.012] & [-0.047, 0.026] & [-0.049, -0.015] & [0.777, 0.924]\\
1 & 20 & 10,000 & 10,000 & [-0.040, 0.012] & [-0.052, 0.028] & [-0.073, -0.032] & [0.525, 0.975]\\
1 & 25 & 10,000 & 10,000 & [-0.042, 0.003] & [-0.069, 0.024] & [-0.109, -0.067] & [0.598, 1.010]\\
1 & 30 & 10,000 & 10,000 & [-0.050, 0.003] & [-0.064, 0.016] & [-0.143, -0.104] & [0.122, 1.069]\\
\hline
\end{tabular}
\begin{flushleft}
\footnotesize
Simulation results of estimating normalizing constants of models 1 and 2 ($\mathcal{M} = 1, 2$ in column 1) using the algorithm in Section~\ref{subsubsec:EstimationProcedure}. the last four columns are the 95\% confidence intervals of $\log(\hat{C}/C)$, where $\hat{C}$ is the estimator and $C$ is the true value, in 100 repeated simulations using the importance sampling with Gaussian tail ($CI_1$), the importance sampling with t tail with a degree of freedom 2 ($CI_2$), the reciprocal importance sampling with Gaussian tail ($CI_3$) and the reciprocal importance sampling with t tail degree of freedom 2 ($CI_4$). $N_{sim}$ is the number of observations and $N_{is}$ is the number of samples from the importance function; $D$ is the dimension of the space.
\label{table:simulation_results_normalizing_const}
\end{flushleft}
\end{table*}

\begin{table*}[tbph]
\centering
\caption{Table \ref{table:simulation_results} in the main text continued.}
\footnotesize
\begin{tabular}{|c|c|c|cc|cc|cc|cc|}
\hline
\hline
\multirow{2}{*}{K} & \multirow{2}{*}{$\sigma$} & \multirow{2}{*}{$n$} & \multicolumn{2}{|c|}{$Q_K = P_K^{(1)}$} & \multicolumn{2}{|c|}{$Q_K = P_K^{(2)}$} & \multicolumn{2}{|c|}{$Q_K = P_K^{(3)}$} & \multicolumn{2}{|c|}{$Q_K = P_K^{(4)}$}\\
& & & ML & BIC  & ML & BIC & ML & BIC & ML & BIC\\
\hline
3 & 0.2 & 2000 & 98   & 100  & 96   & 100 & 97   & 100  & 100  & 100  \\
3 & 0.3 & 2000 & 100  & 100  & 96   & 100 & 96.5 & 100  & 100  & 100  \\
3 & 0.4 & 2000 & 99   & 46   & 100  & 100 & 99.5 & 100  & 100  & 100  \\
3 & 0.5 & 2000 & 3.5  & 0    & 100  & 100 & 99   & 100  & 99.5 & 84   \\
\hline
4 & 0.2 & 2000 & 93.5 & 100  & 97   & 100 & 100  & 99   & 99   & 100  \\
4 & 0.3 & 2000 & 88.5 & 91.5 & 90.5 & 100 & 97   & 100  & 97.5 & 100  \\
4 & 0.4 & 2000 & 6    & 0    & 98.5 & 100 & 100  & 100  & 85.5 & 11.5 \\
4 & 0.5 & 2000 & 0    & 0    & 98.5 & 100 & 98   & 99   & 1    & 0    \\
\hline
5 & 0.2 & 2000 & 93   & 100  & 98.5 & 100 & 98   & 89.5 & 95   & 100  \\
5 & 0.3 & 2000 & 68   & 28.5 & 99   & 100 & 99.5 & 95   & 74.5 & 86.5 \\
5 & 0.4 & 2000 & 0    & 0    & 97.5 & 100 & 100  & 94.5 & 0    & 0    \\
5 & 0.5 & 2000 & 0    & 0    & 98.5 & 100 & 99   & 100  & 0    & 0   
\\
\hline
\end{tabular}
\begin{flushleft}
{$n=2000$ observations are considered.}
\end{flushleft}
\label{table:simulation_results_continued}
\end{table*}

\begin{table}[tbph]
    \centering
    \caption{Table \ref{tab:n2000heter} in the main text continued.}
    \begin{tabular}{|c|c|c|cc|cc|cc|cc|}
\hline
\hline
\multirow{2}{*}{K} & \multirow{2}{*}{$\sigma$} & \multirow{2}{*}{$n$} & \multicolumn{2}{|c|}{$Q_K = P_K^{(1)}$} & \multicolumn{2}{|c|}{$Q_K = P_K^{(2)}$} & \multicolumn{2}{|c|}{$Q_K = P_K^{(3)}$} & \multicolumn{2}{|c|}{$Q_K = P_K^{(4)}$}\\
& & & ML & BIC  & ML & BIC & ML & BIC & ML & BIC\\
\hline
3&	0.2& 200&	48.5&	23.5&	82&56&		78&66&	66 & 34.5\\
3&	0.3&200&	26&11&	57&29.5	& 57&	39.5& 41.5 & 18.5\\
3&	0.4&200&	13.5& 2&32.5& 12&		41&27.5& 20.5 &	9.5\\
3&	0.5&200&	8.5& 2&		21.5& 6&	32& 16.5&	12.5 & 5.5\\
4&	0.2&200&	19.5& 5&	49&	29.5&	34&23.5&			32.5 & 8\\
4&	0.3&200&	5& 0 &	20.5& 10&	16.5& 11.5&	7.5 & 1\\
4&	0.4&200&	0.5& 0&	6&1.5&	5&	3&	1 & 0\\
4&	0.5&200&	0&	0&		2.5& 1&	2.5& 	0.5&	0&	0\\
5&	0.2 &200& 0&	0&		27& 15&		13.5& 6& 17 &	3\\
5&	0.3&200&	0.5& 0&	7.5 & 1.5	&	6&2.5 & 4&	0.5\\
5&	0.4&200&	0.5& 0&	0&	0 & 1 &	0.5&	1 & 0\\
5&	0.5&200&	0&	0&	0.5&	0&	1&	0.5&	0&	0\\
\hline
    \end{tabular}
   \begin{flushleft}$n=200$ observations with heterogeneous variances.
   \end{flushleft}
    \label{tab:n200heter}
\end{table}

\section{ Simulation Robustness}
\label{appendix:simulationrobustness}

In Section \ref{section:computation}, we propose an estimation procedure to approximate marginal likelihoods, which results from viewing it as a problem of estimating normalization constants. Note that the consistency theorem in the main text is for exact marginal likelihoods, whereas in practice, we can only approximate the marginal likelihoods using our proposed estimators. We now present the robustness of this estimation procedure for marginal likelihoods, towards estimating the number of hidden states. We first prove a general property for estimating normalization constants in Section \ref{subsec:Robust_NormalConst} (not limited to HMM). We then show the robustness of the adopted estimator in Section \ref{subsec:AdoptedEst} under the HMM setting in the main text.

\subsection{ Property of the Normalizing Constant Estimator}\label{subsec:Robust_NormalConst}

The following is a general property of the normalizing constant estimator using the locally restricted version, which is not limited to HMMs.

Let $\tilde{p}(z) = p(z)/C$ be a pdf on $\Omega$, for which we only know the un-normalized density $p(z)$; the normalizing constant $C = \int p(z) dz$ is unknown. Let $Z_{1:N} = \{ Z_1, \cdots, Z_N \}$ be independent random samples from $\tilde{p}$, and $\{g(\cdot|t): t \in \mathcal{T}\}$ be a family of densities indexed by parameter $t \in\mathcal{T}$. Given $Z_{1:N} = z_{1:N}$, define
\begin{equation*}
    \hat{T}_N = \hat{T}_N(z_{1:N}) := \mbox{argmax}_T  \prod_{i=1}^N g(z_i|T).
\end{equation*}
Assume that $V_{1:M} = \{V_1, \cdots, V_M\}$ are independent samples from $g(\cdot | \hat{T}_N)$. In Section \ref{section:computation} we have shown that, given $V_{1:M} = v_{1:M}$, $C$ can be approximated by the locally restricted importance sampling estimator $\hat{C}_{n, N}^{loc}$ defined as
\begin{align}\label{C_loc}
\hat{C}_{M, N}^{loc} = \hat{C}_{M, N}^{loc}(z_{1:N}, v_{1:M}) := \frac{1}{M \hat{P}_{\Omega_r}} \bigg[ \sum_{j=1}^{M} \frac{p(v_j)}{g(v_j|\hat{T}_N)} 1_{\Omega_r}(v_j) \bigg],
\end{align}
where $\Omega_r \subset \Omega$ and $\hat{P}_{\Omega_r} =  \frac{1}{N} \sum_{i=1}^N 1_{\Omega_r}(z_i)$.

\begin{lemma}
\label{lemma:simrobustlemma}
Let $\Omega_r$ and $\hat{C}_{M,N}^{loc}$ be defined as in (\ref{C_loc}). Assume that for any $t \in \mathcal{T}$, $\frac{1}{2} < \int_{\Omega_r} g(z|t) < 1$ and $\tilde{c}_l \leq \frac{\tilde{p}(\cdot)}{g(\cdot|t)} \leq \tilde{c}_u$ on $\Omega_r$ for some constants $\tilde{c}_u > \tilde{c}_l > 0$. Then, there exists $L_M, L_N > 0$ such that for all $M > L_M$ and $N > L_N$, 
\begin{equation}
\label{eqn:simrobustlemma}
      \frac{Var[\hat{C}_{M, N}^{loc}(Z_{1:N}, V_{1:M})]}{C^2} \leq \left[\frac{1}{M}+\frac{1}{N}\right] \max\left\{ \frac{2\tilde{c}_u}{\tilde{c}_l}, \frac{2}{\tilde{c}_l}\right\}.
\end{equation}
Here $C=\int p(z)dz$ is the normalizing constant.
\end{lemma}

\begin{proof}[Proof of Lemma \ref{lemma:simrobustlemma}]
Note that the variance of $\hat{C}_{M, N}^{loc}$ is the ratio of two independent random variables $\frac{1}{M}\sum_{j=1}^{M} \frac{p(Y_j)}{g(V_j|\hat{T}_N)} 1_{\Omega_r}(V_j)$ and $\frac{1}{N}\sum_{i=1}^N 1_{\Omega_{r}}(Z_i)$. The expectations and variances of them can be computed as follows.
\begin{equation*}
E\left(\frac{1}{N}\sum_{i=1}^N 1_{\Omega_{r}}(Z_i)\right) = \int_{\Omega_r} \tilde{p}(z) dz,
\end{equation*}
\begin{align*}
Var\left[\frac{1}{N}\sum_{i=1}^N 1_{\Omega_{r}}(Z_i)\right] = \frac{1}{N}  \int_{\Omega_r} \tilde{p}(z) dz \left[ 1 - \int_{\Omega_r} \tilde{p}(z) dz \right];
\end{align*}
\begin{align*}
E\left(\frac{1}{M}\sum_{j=1}^{M} \frac{p(V_j)}{g(V_j|\hat{T}_N)} 1_{\Omega_r}(V_j)\right) 
= & E\left[ E\left(\frac{1}{M}\sum_{j=1}^{M} \frac{p(V_j)}{g(V_j|\hat{T}_N)} 1_{\Omega_r}(V_j) \bigg| \hat{T}_N\right) \right] \\
= & \int_{\Omega_r} p(z) dz = C \int_{\Omega_r} \tilde{p}(z) dz;
\end{align*}
\begin{align*}
Var\left[\frac{1}{M}\sum_{j=1}^{M} \frac{p(y_j)}{g(V_j|\hat{T}_N)} 1_{\Omega_r}(V_j)\right]
= & Var\left[ E\left(\frac{1}{M}\sum_{j=1}^{M} \frac{p(V_j)}{g(V_j|\hat{T}_N)} 1_{\Omega_r}(V_j) \bigg| \hat{T}_N\right) \right] \\ 
 & + E\left[ Var\left(\frac{1}{M}\sum_{j=1}^{M} \frac{p(V_j)}{g(V_j|\hat{T}_N)} 1_{\Omega_r}(V_j) \bigg| \hat{T}_N\right) \right]\\
= & E\left[ \frac{1}{M} \int_{\Omega_r} \frac{p^2(z)}{g(z|\hat{T}_N)} dz \right] - \frac{1}{M} \left(\int_{\Omega_r} p(z) dz\right)^2\\ = & C^2 E\left[ \frac{1}{M} \int_{\Omega_r} \frac{\tilde{p}^2(z)}{g(z|\hat{T}_N)} dz \right] - C^2\frac{1}{M} \left(\int_{\Omega_r} \tilde{p}(z) dz\right)^2.
\end{align*}
By delta method, we have 
\begin{align}
\notag
 Var(\hat{C}_{M, N}^{loc})
= & \frac{Var\left[\frac{1}{M}\sum_{j=1}^{M} \frac{p(V_j)}{g(V_j|\hat{T}_N)} 1_{\Omega_r}(V_j)\right]}{\left[E\left(\frac{1}{N}\sum_{i=1}^N 1_{\Omega_{r}}(Z_i)\right)\right]^2}\\ 
\notag
&+ \frac{Var\left[\frac{1}{N}\sum_{i=1}^N 1_{\Omega_{r}}(Z_i)\right]\left[E\left(\frac{1}{M}\sum_{j=1}^{M} \frac{p(V_j)}{g(V_j|\hat{T}_N)} 1_{\Omega_r}(V_j)\right)\right]^2}{\left[E\left(\frac{1}{N}\sum_{i=1}^N 1_{\Omega_{r}(Z_i)}\right)\right]^4} \\
\notag
& + o(M^{-1}) + o(N^{-1})\\
\notag
= & \frac{C^2}{M}\left\{ \frac{E\left[ \int_{\Omega_r} \frac{\tilde{p}^2(y)}{g(z|\hat{T}_N)} dz \right]}{\left(\int_{\Omega_r} \tilde{p}(z) dz\right)^2} - 1 \right\}  + \frac{C^2}{N} \left[\left( \int_{\Omega_r} \tilde{p}(z) dz \right)^{-1} - 1\right] \\
& + o(M^{-1}) + o(N^{-1}).
\label{eqn_pf_NormalConst}
\end{align}
Now, by the condition on $g(\cdot|t)$, we have
\begin{align*}
    &\tilde{c}_l \int_{\Omega_r} {\tilde{p}(z)} dz\leq \int_{\Omega_r} \frac{\tilde{p}^2(z)}{g(z|\hat{T}_N)} dz \leq \tilde{c}_u \int_{\Omega_r} {\tilde{p}(y)} dy, \\
    & \frac{1}{\tilde{c}_u}\leq \left[\int_{\Omega_r} {\tilde{p}(y)} dy\right]^{-1} \leq \frac{2}{\tilde{c}_l}.
\end{align*}
Plugging these back into (\ref{eqn_pf_NormalConst}) and taking $M$ and $N$ sufficiently large, we complete the proof.
\end{proof}

\subsection{ Robustness of Adopted Estimation Procedure}\label{subsec:AdoptedEst}

After having the result for general settings in Section \ref{subsec:Robust_NormalConst}, we now show the robustness result under the HMM setting in the main text.
Recall that for any number of state $K$, given ${\bm Y}_{1:n}={\bm y}_{1:n}$,  $p_K({\bm y}_{1:n})$ denotes the corresponding marginal likelihood. In Section \ref{subsubsec:EstimationProcedure}, we propose a way of estimating $p_K({\bm y}_{1:n})$ via the importance sampling. Let $N(n)$ and $M(n)$ be positive integers. Using the notations in Section \ref{subsec:Robust_NormalConst}, we first draw $Z_{1:N(n)}^K$ from $p_K( \cdot |{\bm y}_{1:n})$, the posterior distribution of $\boldsymbol{\phi}_K$ given the observed trajectory, and choose an importance function $g_K(\cdot|\hat{T}_{N(n)}^K)$ indexed by some parameter $\hat{T}_{N(n)}^K = \hat{T}_{N(n)}^K (Z_{1:N(n)}^K)$. Then, we draw $V_{1:M(n)}^K$ from $g_K(\cdot|\hat{T}_{N(n)}^K)$. Given $V_{1:M(n)}^K = v_{1:M(n)}^K$, Section \ref{subsubsec:EstimationProcedure} gives the estimated marginal likelihood
\begin{align*}
    \hat{p}_K({\bm y}_{1:n}, {\bm v}_{1:M(n)}^K) 
    = \frac{1}{M(n)\hat{P}_{\Omega_{n,r}}} \sum_{j=1}^{M(n)} \frac{p_K({\bm y}_{1:n},v_j^K)}{g_K(v_j^K|\hat{T}_{N(n)}^K)} 1_{\Omega_{n,K,r}}(v_j),
\end{align*}
where $\Omega_{n,K,r} \subset \Phi_K$, the parameter space for a $K$-state HMM.

We claim that this estimated marginal likelihood can approximate the marginal likelihood reasonably well so that the consistency of $\hat{K}_n$ still holds if we replace the marginal likelihood with its estimated one. Let $P$ be the probability measure such that, for any $n$ and any $K$,
\begin{align}
\notag
    & P\{ {\bm X}_{1:n} = {\bm x}_{1:n}, {\bm Y}_{1:n} \in d{\bm y}_{1:n}, Z_{1:N(n)}^K \in dz_{1:N(n)}^K, V_{1:M(n)}^K \in dv_{1:M(n)}^K \} \\
\notag
= & P^*\left\{ {\bm X}_{1:n} = {\bm x}_{1:n}, {\bm Y}_{1:n} \in d{\bm y}_{1:n} \right\} \times \prod_{i=1}^{N(n)} p_K(z_i^K|{\bm y}_{1:n}) \times \prod_{j=1}^{M(n)} g_K(v_j^K|\hat{T}_{N(n)}^K),
\end{align}
i.e., $P$ is defined by the product of three measures (1) the probability measure under the true parameter for $({\bm X}_{1:n}, {\bm Y}_{1:n})$, which is denoted by $P^*$, (2) given ${\bm Y}_{1:n}={\bm y}_{1:n}$, the probability measure for the vector $Z_{1:N(n)}^K$ with i.i.d. components with density function $p_K(\cdot | {\bm y}_{1:n})$ and (3) given $Z_{1:N(n)}^K=z_{1:N(n)}^K$, the probability measure for the vector $V_{1:M(n)}^K$ with i.i.d. components with density function $g_K(\cdot | \hat{T}_{N(n)}^K)$.

\begin{theorem}
\label{theorem:simrobustnessnc}
Assume $K \neq K^*$ and that the following conditions hold.
\begin{itemize}
    \item [(i)] Under $P^*$, $p_K({\bm Y}_{1:n})/p_{K^*}({\bm Y}_{1:n}) \xrightarrow{n \rightarrow \infty} 0$ in probability.
    \item [(ii)] There exist $\tilde{c}_u>\tilde{c}_l>0$ and $D(n) \rightarrow \infty$ such that, for any $N(n) \geq D(n)$, we have $P\{\mathcal{E}_{n,K}\} \rightarrow 1$ and $P\{\mathcal{E}_{n,K^*}\} \rightarrow 1$, where
        \begin{align*}
        \mathcal{E}_{n,K} := \left\{ \forall \phi \in \Omega_{n,K,r}, \frac{p_K(\phi|{\bm y}_{1:n})}{g(\phi|\hat{T}_{N(n)}^K)} \in (\tilde{c}_l, \tilde{c}_u) \mbox{ and } \int_{\Omega_{n,K,r}} g(\phi|\hat{T}_{N(n)}^K) d\phi > 1/2 \right\}.
    \end{align*}
\end{itemize}
Then, there exists $A(n) \rightarrow \infty$ such that, for any sequences of $M(n)$ and $N(n)$ satisfying $M(n) \geq A(n)$ and $N(n) \geq A(n)$, we have, under $P$, $$\hat{p}_K({\bm Y}_{1:n}, V_{1:M(n)}^K)/\hat{p}_{K^*}({\bm Y}_{1:n}, V_{1:M(n)}^{K^*}) \xrightarrow{n \rightarrow \infty} 0 \mbox{in probability.}$$
\end{theorem}

\begin{proof}[Proof of Theorem \ref{theorem:simrobustnessnc}]
Let us first focus on $K$. For any $n$, let $L_{M,n}$ and $L_{N,n}$ be the corresponding bounds in Lemma \ref{lemma:simrobustlemma} so that (\ref{eqn:simrobustlemma}) holds for all $M(n) > L_{M,n}$ and $N(n) > L_{N,n}$. Consider $M(n) > L_{M,n}$, $N(n) > \max\{ L_{N,n}, D(n)\}$. Since $N(n) > D(n)$, we have $P\{ \mathcal{E}_{n,K}^c \} \rightarrow 0$. Let $\hat{\sigma}_K^2({\bm y}_{1:n})$ be the variance of $\hat{p}_K({\bm Y}_{1:n}, V_{1:M(n)}^K)$ conditioning on ${\bm Y}_{1:n} = y_{1:n}$. Set $c_{n,K} \rightarrow 0$ satisfying
\begin{equation*}
 c_{n,K} \left[\frac{1}{M}+\frac{1}{N}\right] \max\left\{ \frac{2\tilde{c}_u}{\tilde{c}_l}, \frac{2}{\tilde{c}_l}\right\} \rightarrow \infty.
\end{equation*}
Then, since $M(n) > L_{M,n}$ and $N(n) > L_{N,n}$, by Lemma \ref{lemma:simrobustlemma}, we have
\begin{equation}
\label{eqn:fromlemmavarnormalizingconst}
P \left\{ \hat{\sigma}_K({\bm Y}_{1:n}) \leq c_{n,K}^{1/2} p_K({\bm Y}_{1:n}); \mathcal{E}_{n,K} \right\} \xrightarrow{n \rightarrow \infty} 1.
\end{equation}
Let $\{b_{n, K}\}_{n\geq 1}$ be a sequence of positive constants such that $\lim_{n\rightarrow\infty}b_{n,K} c_{n,K} =  \gamma< 1$ and $\lim_{n\rightarrow\infty} b_{n,K} = \infty$. For any $n$ and any ${\bm y}_{1:n}$, from the Chebychev's inequality,

\begin{equation*}
    P\left\{ \vert \hat{p}_K({\bm Y}_{1:n}, V_{1:M(n)}^K) - p_K({\bm Y}_{1:n}) \vert > b_{n,K}^{1/2} \hat{\sigma}_K({\bm Y}_{1:n}); \mathcal{E}_{n,K} \Big\vert {\bm Y}_{1:n} = {\bm y}_{1:n} \right\} \leq b_{n,K}^{-1},
\end{equation*}
and since this holds for all ${\bm y}_{1:n}$, we have
\begin{equation*}
    P\left\{ \vert \hat{p}_K({\bm Y}_{1:n}, V_{1:M(n)}^K) - p_K({\bm Y}_{1:n}) \vert > b_{n,K}^{1/2} \hat{\sigma}_K({\bm Y}_{1:n}); \mathcal{E}_{n,K} \Big\vert {\bm Y}_{1:n} = {\bm y}_{1:n} \right\} \leq b_{n,K}^{-1},
\end{equation*}
As a consequence, 
\begin{align*}
& P  \left\{ 1 - b_{n,K}^{1/2} \frac{\hat{\sigma}_K({\bm Y}_{1:n})}{p_K({\bm Y}_{1:n})} \leq \frac{\hat{p}_K({\bm Y}_{1:n}, V_{1:M(n)}^K)}{p_K({\bm Y}_{1:n})} \leq 1 + b_{n,K}^{1/2} \frac{\hat{\sigma}_K({\bm Y}_{1:n})}{p_K({\bm Y}_{1:n}); \mathcal{E}_{n,k}} ; \mathcal{E}_{n,K}\right\}\\
= & 1- P\left\{ \vert \hat{p}_K({\bm Y}_{1:n}, \boldsymbol{\varphi}_K) - p_K({\bm Y}_{1:n}) \vert > b_{n,K}^{1/2} \hat{\sigma}_K({\bm Y}_{1:n}) ; \mathcal{E}_{n,k} \right\} \\
\geq  & 1-b_{n,K}^{-1}\xrightarrow{n \rightarrow \infty} 1.
\end{align*}
Together with Equation~(\ref{eqn:fromlemmavarnormalizingconst}) and $P\{ \mathcal{E}_{n,K}^c\} \rightarrow 0$, we have
\begin{equation}
\label{eqn:bncnineq1}
P \left\{ 1 - b_{n,K}^{1/2} c_{n,K}^{1/2}\leq \frac{\hat{p}_K({\bm Y}_{1:n}, V_{1:M(n)}^K)}{p_K({\bm Y}_{1:n})} \leq 1 +  b_{n,K}^{1/2} c_{n,K}^{1/2} \right\} \xrightarrow{ n \rightarrow \infty} 1.
\end{equation}
Similarly, we have
\begin{equation*}
P \left\{ 1 - b_{n,K^*}^{1/2} c_{n,K^*}^{1/2}\leq \frac{\hat{p}_{K^*}({\bm Y}_{1:n}, V_{1:M(n)}^{K^*})}{p_{K^*}({\bm Y}_{1:n})} \leq 1 +  b_{n,K^*}^{1/2} c_{n,K^*}^{1/2} \right\} \xrightarrow{ n \rightarrow \infty} 1.
\end{equation*}
We complete the proof by noticing that (i) the condition implies that $$p_{K}({\bm Y}_{1:n}) / p_{K^*}({\bm Y}_{1:n}) \xrightarrow{n \rightarrow \infty} 0 \mbox{ in probability (under $P$),}$$ 
(ii) $\lim_{n\rightarrow\infty}b_{n,K} c_{n,K}\rightarrow \gamma < 1$, $\lim_{n\rightarrow\infty}b_{n,K^*} c_{n,K^*}\rightarrow \gamma < 1$. 
\end{proof}

\end{document}